\documentclass{m2anBjoern}

\usepackage{accents}
\usepackage[toc,page]{appendix}
\usepackage{enumerate}

\newcommand{\myim}[1]{{\mathtt{im}( {#1})}} 
\newcommand{\myker}[1]{{\mathtt{ker}( {#1})}} 
\newcommand{\mydim}{\mathtt{dim}} 
\newcommand{\bv}[1]{{\mathbf #1}} 
\newcommand{\bt}[1]{{\mathbf{#1}}}
\newcommand{\scalco}[1]{{\mathfrak{#1}}} 

\newcommand{\tfa}{q}
\newcommand{\tfb}{w}
\newcommand{\sfa}{\rho}
\newcommand{\sfb}{m}

\newcommand{\Massa}{\bt{Q}}
\newcommand{\Massb}{\bt{W}}
\newcommand{\MassBoth}{\bt{M}}

\newcommand{\funSpace}[1]{\mathcal{#1}} 
\newcommand{\Vsa}{\funSpace{Q}}
\newcommand{\Vsb}{\funSpace{W}}

\newcommand{\Pfun}{\Pi}
\newcommand{\Pvec}{\bt{\Pi}}

\newcommand{\vwei}{\boldsymbol{\xi}}
\newcommand{\wei}{\xi}

\newcommand{\ubv}[1]{{\underaccent{\bar}{\mathbf{#1}}}} 
\newcommand{\Lscal}{\langle}	\newcommand{\Rscal}{\rangle} 
\newcommand{\scps}{{\ast}}

\newcommand{\HamPDE}{\mathcal{H}}

 \newcommand{\blsVertex}{{{\nu}}} 
 \newcommand{\TraceOp}{{\mathcal{T}}}
 \newcommand{\Onepipe}{e}
 
\usepackage{tikz}
\usetikzlibrary{calc,positioning,shapes,arrows,decorations.shapes} 
\tikzstyle{rect}=[
   rectangle, draw=black, font=\small\sffamily\bfseries, inner sep=9pt]

\newtheorem{system}[thrm]{System}
\newtheorem{assumption}[thrm]{Assumption}

\begin{document}
\title[]{On snapshot-based model reduction under compatibility conditions for a nonlinear flow problem on networks}\thanks{The support of the German Federal Ministry of Education and Research (BMBF) via the project \textit{EiFer} is acknowledged. Moreover, we thank for the support of the DFG research
training group 2126 on algorithmic optimization.}
\author{Bj{\"o}rn Liljegren-Sailer}\address{Universit{\"a}t Trier, FB IV - Mathematik, Lehrstuhl Modellierung und Numerik, D-54286 Trier, Germany
\email{Corresponding author: bjoern.sailer@uni-trier.de}}
\author{Nicole Marheineke}\sameaddress{1}
%
%
\begin{abstract}
This paper is on the construction of structure-preserving, online-efficient reduced models for the barotropic Euler equations with a friction term on networks. The nonlinear flow problem finds broad application in the context of gas distribution networks. We propose a snapshot-based reduction approach that consists of a mixed variational Galerkin approximation combined with quadrature-type complexity reduction. Its main feature is that certain compatibility conditions are assured during the training phase, which make our approach structure-preserving. The resulting reduced models are locally mass conservative and inherit an energy-bound and port-Hamiltonian structure. We also derive a well-posedness result for them. In the training phase, the compatibility conditions pose challenges, we face constrained data approximation problems as opposed to the unconstrained training problems in the conventional reduction methods. The training of our model order reduction consists of a principal component analysis under a compatibility constraint and, notably, yields reduced models that fulfill an optimality condition for the snapshot data. The training of our quadrature-type complexity reduction involves a semi-definite program with combinatorial aspects, which we approach by a greedy procedure.
Efficient algorithmic implementations are presented. The robustness and good performance of our structure-preserving reduced models are showcased at the example of gas network simulations.
\end{abstract}


\maketitle

\section{Introduction}

Conventional discretization methods, such as finite difference or finite element methods, are powerful tools for the numerical investigation of engineering applications. Nonetheless, they may reach their limits due to the high computational demand in many-query tasks of complex applications. Under certain circumstances, projection-based model reduction methods may help to break computational limits and speed up simulations profoundly \cite{art:morParamHamHesthaven,art:comred-ecsw,art:hyperreduction-fem-cubature}. In the linear case, a low dimension of the projection directly translates into an efficient reduced model. The supplementation by complexity reduction becomes necessary in the presence of nonlinearities.
A prevailing approach to obtain reduced models is to train them towards given snapshot-data. The proper orthogonal decomposition \cite{art:pod-kunisch,art:increm-POD} and reduced basis method \cite{book:morHesRS16} are widely used examples for such model order reduction methods. As for the complexity reduction, we refer to the empirical interpolation method and its variants \cite{art:deim-state-space-err,art:empint-maday04} as well as to the quadrature-type approaches \cite{art:efficient-integration-cubature,art:comred-ecsw}.

While reduction methods have been successfully used in numerous applications, it is known that the standard approaches can suffer from poor results and stability issues for problems of complex structure.
One major cause for the issues is that the reduced models may disregard fundamental structural properties of the model they approximate. These properties might include, e.g., conservation laws, dissipative relations or symplecticity. As a remedy, structure-preserving approximation methods have been developed, see, e.g., \cite{art:topics-in-strpresdisc,art:avfPDE,art:morParamHamHesthaven} for selected overviews. Problem-adapted space discretization approaches range from finite volume methods \cite{book:leveque-finite-volume-methods} to mixed finite element methods \cite{book:BreeziMixedFE08} or mimetic finite differences \cite{art:FISHER13}. As for the time discretization, e.g., the symplectic and geometric integrators \cite{book:hairerGeomInt,art:pH-symplectictime-koty} are famous structure-preserving methods.
The development of structure-preserving reduction methods is a more recent topic. Frequency-based model order reduction methods for linear port-Hamiltonian models are discussed in \cite{inbook:interpolation-based-port-Hamiltonian-Systems,art:WOLF2010401}. The works \cite{art:hyperreduction-preserving-lagrangian-structure, art:carlberg-PetrovGal11} treat nonlinear Lagrangian dynamics by snapshot-based approaches, particularly also using empirical interpolation-type complexity reduction. The papers \cite{art:morParamHamHesthaven,art:moramHamDissHesthaven,art:symplHamMor} investigate model order reduction methods that preserve canonical symplectic structure under certain compatibility conditions and thus yield reduced systems in Hamiltonian form. Additionally, complexity reduction by the discrete empirical interpolation method is considered, but this step is not strictly structure-preserving, i.e., it does not guarantee a Hamiltonian structure. A modification of the discrete empirical interpolation that allows for a Hamiltonian representation can be found in \cite{art-mor-structure-preserving-nonlinear-pH}, but it inherits higher local errors due to an enforced symmetrization step it includes. In other contexts  \cite{art:hyperreduction-fem-cubature,art:comred-ecsw}, quadrature-type complexity reduction methods are used, which show to be structure-preserving in a more natural way. 

In this paper, we construct structure-preserving and online-efficient reduced models for a nonlinear flow problem on networks. The problem is governed by the barotropic Euler equations with a friction term and inherits a port-Hamiltonian structure. It finds application in the context of gas transport networks. The efficient simulation, control and optimization of gas networks is an important topic as the wide range of publications in this direction shows. In \cite{art:models_revisted} a model hierarchy of gradually simplified models is derived, whereby the isothermal (i.e., barotropic) Euler equations are the starting point. Strongly related are the contributions \cite{art:Domschke2015,art:Domschke2018}, which apply adaptive model-switching in a gas simulation tool. Other approaches to speed up calculations include problem-adapted preconditioners \cite{art:Qiu2020} or concepts from discrete optimization \cite{book:gasKochHiller15}. Conventional snapshot-based model reduction methods have also been applied for gas networks \cite{art:mor-gas-grundel-jansen,art:himpe-gasmor21}. Our model reduction approach focuses on structure-preserving properties, i.e., local mass conservation, an energy bound and a port-Hamiltonian structure underlying the model problem. We  establish these properties using compatibility conditions related to the mixed variational formulation \cite{art:lilsailer-nlfow}, see also \cite{art:egger-mfem-compressEuler,inproc:bls-hyp2018}. While beneficial for the robustness and performance of the reduced models, the compatibility conditions pose a challenge in the training phase, as one faces constraints in the training problems related to the snapshot data. We aim for reduced models that fulfill an optimality condition under these compatibility constraints. The conventional model reduction approaches, cf.~\cite{art:morKunV01,book:dimred2003,art:empint-maday04}, are designed towards an optimality condition related to the snapshot data, but compatibility conditions are not taken into account, i.e., structural properties are generally not guaranteed. For some applications, more problem adapted approaches have been derived. In \cite{art:Rozza13,art:stab-redbasis19,book:morHesRS16} inf-sup stable reduced models for the Stokes and Navier-Stokes equations are derived. This is done using a splitted procedure that treats the training towards the snapshot data and the inf-sup compatibility separately from each other. Consequently, no optimality condition related to the snapshot data can be guaranteed in this approach. Structure-preserving model order reduction methods that aim for optimal reduced models have been considered for canonical Hamiltonian systems, see  \cite{art:morParamHamHesthaven,art:symplHamMor}. Conceptually, our approach shares most similarities with the latter approaches, which are also known as symplectic model reduction.
However, the model problem we treat is quite different. Particularly, we cannot reformulate our problem as a canonical Hamiltonian system, we have to deal with network aspects, and our approximations rely on a mixed variational formulation. Regarding the model order reduction, we show that our training problem with its compatibility constraints can be attributed to an unconstrained principal component analysis, given appropriate norms are chosen. Based on this, we derive an efficient algorithmic implementation. Similarly, a training problem involving the snapshot data and constraints reflecting another compatibility condition is considered for the complexity reduction. It is a semi-definite program with combinatorial aspects, which we approach algorithmically by a greedy procedure.
The performance of our proposed model reduction method is investigated numerically at the example of gas networks and compared to conventional, non-structure-preserving reduction methods. Moreover, we show a well-posedness result for our reduced models.

The structure of the paper is as follows. In Section~\ref{sec-netproblem} we state the model problem and our notation concerning the network aspects. The approximation ansatz from \cite{art:lilsailer-nlfow} is shortly reviewed in Section~\ref{sec:approx-ansatz}. It is employed to set up the underlying full order model (space discretization) and also serves as a basis for the reduction steps. 
Our structure-preserving model order- and the complexity-reduction methods are derived and analyzed in Section~\ref{sec:prop-mr} and Section~\ref{sec:cr-approx}, respectively. The numerical studies are presented in Section~\ref{sec:numres}, and a well-posedness result for our reduced models is derived in Section~\ref{sec:well-posed}.

\section{Model problem} \label{sec-netproblem}
This section introduces our nonlinear model problem, which describes, e.g., gas transportation networks \cite{art:models_revisted,art:model-gasdistrib-2009}. Moreover, the variational principle our approximations are based on is stated, and the notation concerning the network aspects is settled, cf. \cite{art:kugler-dwe-net, phd:liljegren}.

Let a (gas pipe) network be described by a directed graph $({\mathcal{N}},\mathcal{E})$ with sets of nodes ${\mathcal{N}} = \{\blsVertex_1, \ldots, \blsVertex_\ell \}$ and edges $\mathcal{E} = \{ \Onepipe_1, \ldots, \Onepipe_k \} \subset {\mathcal{N}} \times {\mathcal{N}}$. Each edge $\Onepipe\in \mathcal{E}$ is attached with a length  $l^\Onepipe>0$ and a weight $A^\Onepipe >0$. 
The set of all edges adjacent to the node $\blsVertex$ is denoted by $\mathcal{E}(\blsVertex) =  \{\Onepipe \in \mathcal{E}: \, \Onepipe= ( \blsVertex, \bar{\blsVertex}), \text{ or } \Onepipe=( \bar{\blsVertex}, \blsVertex) \}$, and a weighted incidence mapping is defined by
\begin{align*}
	 n^\Onepipe[\blsVertex] =
	 \begin{cases}
	 	\,\,\,\, A^\Onepipe & \text{ for } \Onepipe = ( \blsVertex, \bar{\blsVertex}) \text{ for some }  \bar{\blsVertex} \in \mathcal{N} \\
	 	 -A^\Onepipe  & \text{ for } \Onepipe = ( \bar{\blsVertex}, \blsVertex ) \text{ for some }  \bar{\blsVertex} \in \mathcal{N} .
	 \end{cases}
\end{align*}
The nodes are grouped into interior nodes ${\mathcal{N}}_0 \subset {\mathcal{N}}$ and boundary nodes \mbox{${\mathcal{N}}_\partial = {\mathcal{N}} \text{\textbackslash} {\mathcal{N}}_0$}.
Function spaces on the network are constructed by compositions of standard Sobolev spaces for every edge. The spatial domain of the union of edges is  $\Omega_{} = \{x: \, x\in \Onepipe, \text{ for } \Onepipe \in \mathcal{E}\}$. Every edge $\Onepipe$ can be identified with an interval $(0,l^\Onepipe)$ which is tacitly employed in the upcoming integral expressions. The space of square-integrable functions on $\mathcal{E}$ is given by $\funSpace{L}^2(\mathcal{E}) =  \left\{b: \Omega_{}\rightarrow \mathbb{R} \text{ with } b_{|\Onepipe} \in \funSpace{L}^2(\Onepipe) \text{  for all } \Onepipe\in \mathcal{E} \right\}$, where the subscript $._{|\Onepipe}$ indicates the restriction of a function to the edge $\Onepipe$. The respective scalar product and norm read $\Lscal b, \tilde{b} \Rscal =  \sum_{\Onepipe \in \mathcal{E}} A^\Onepipe \int_\Onepipe b \, \tilde{b} \, dx$ and $|| b || = \sqrt{\Lscal b, b \Rscal}$ for $b$, $\tilde{b} \in \funSpace{L}^2(\mathcal{E})$. 
The weak (broken) derivative operator for functions on the network is defined by $(\partial_x b)_{|\Onepipe} =  \partial_x b_{|\Onepipe}$ for $\Onepipe \in \mathcal{E}$. The space of functions with square-integrable weak broken derivative is given as $\funSpace{H}_{pw}^1(\mathcal{E}) = \left\{b\in \funSpace{L}^2(\mathcal{E}) :  \partial_x b \in \funSpace{L}^2(\mathcal{E}) \right\}$. Accordingly, $\funSpace{C}_{pw}^k(\mathcal{E}) = \left\{b: \Omega_{} \rightarrow \mathbb{R} \text{ with } b_{|\Onepipe} \in \funSpace{C}^k(\Onepipe) \text{  for all } \Onepipe\in \mathcal{E}  \right\}$ denotes the space of piecewise $k$-times continuously differentiable functions, $k\geq 0$.
For $b \in \funSpace{H}_{pw}^1(\mathcal{E})$ we indicate node evaluations with squared brackets, i.e., $b_{|\Onepipe}[\blsVertex] \in \mathbb{R}$ for $\blsVertex \in \mathcal{N}$.
The Sobolev space
$\funSpace{H}_{div}^1(\mathcal{E}) = \{b\in \funSpace{H}^1_{pw}(\mathcal{E}) : \sum_{\Onepipe \in \mathcal{E}(\blsVertex) } n^{\Onepipe}[\blsVertex] b_{|\Onepipe}[\blsVertex] = 0, \text{ for }  \blsVertex \in {\mathcal{N}}_0 \}$ incorporates certain coupling conditions at inner nodes.
The boundary nodes $\mathcal{N}_\partial = \{ \blsVertex_1, \ldots, \blsVertex_p \}$ are assumed to be connected to exactly one edge each. Thus, a boundary operator $\TraceOp^{ }: \funSpace{H}_{pw}^1(\mathcal{E}) \rightarrow \mathbb{R}^{p}$ can be defined by $\left[ \TraceOp^{ } b \right]_i =
	 n^\Onepipe[\blsVertex_i] b_{|\Onepipe}[\blsVertex_i] $  for  $\Onepipe \in \mathcal{E}(\blsVertex_i)$, $ i =1,\ldots, p $ and  $b \in H^1_{pw}(\mathcal{E})$.

We consider the barotropic Euler equations with friction governing density and velocity $\rho,v: [0,T]\times \Omega \rightarrow \mathbb{R} $ on the edges of the network,
\begin{subequations} \label{bls-eq:abstr}
\begin{align} \label{bls-eq:abstr-a}
\partial_t \rho + \partial_x \left( \rho v \right) &= 0 , \qquad \partial_t v + \partial_x \frac{v^2}{2} + \partial_x P'(\rho) = -\tilde{r}(\rho,v) \rho v. 
\end{align}
Here, $P'$ is the derivative of the pressure potential $P:\mathbb{R}^+ \rightarrow \mathbb{R}$, which we assume to be strictly convex and twice continuously differentiable for all considered $\rho$. Moreover, the friction term $\tilde{r}$ is assumed to be pointwise positive. The solution components are interconnected by the coupling conditions at $\nu \in \mathcal{N}_0$
\begin{align} \label{bls-eq:abstr-coup}
	\sum_{\Onepipe \in \mathcal{E}(\blsVertex) } n^{\Onepipe}[\blsVertex] (\rho v)_{|\Onepipe}[\nu]= 0  , \qquad
	\frac{(v_{|\Onepipe}[\blsVertex])^2}{2} + P'(\rho_{|\Onepipe}[\blsVertex]) = \frac{(v_{|\tilde\Onepipe}[\blsVertex])^2}{2} + P'(\rho_{|\tilde\Onepipe}[\blsVertex]) \quad \text{for } \Onepipe,\tilde{\Onepipe} \in \mathcal{E}(\blsVertex).
\end{align}
The model is closed by appropriate boundary and initial conditions, e.g., 
\begin{align} \label{bls-eq:abstr-bc}
\frac{(v[\nu_i])^2}{2} + P'(\rho[\nu_i]) = u_{\nu_i}, \qquad \text{ for given } u_{\nu_i} \in \funSpace{C}([0,T];\mathbb{R}), \quad \nu_i \in \mathcal{N}_\partial,
\end{align}
and $\rho(0,x) = \rho_0(x)$, $\rho v(0,x) = m_0(x)$ with suitable $\rho_0,m_0 \in \funSpace{C}^1_{pw}(\mathcal{E})$.
\end{subequations}
Instead of the pressure potential, the flow can also be described by means of the pressure $p$ given by $p(\rho) = \rho P'(\rho) - P(\rho)$, cf., \cite{art:gieselmann-rel-energy,art:antonelli-quantumHydrodyn}. 
The physical energy of the system is given by the Hamiltonian
\begin{align*}
	\tilde{\HamPDE}(\rho,v) = \int_\Omega \rho \frac{v^2}{2} + P(\rho) dx .
\end{align*}
Conservation of mass as well as of the Hamiltonian at  inner nodes are ensured by the coupling conditions.
Furthermore, for sufficiently smooth solutions, it can be shown that the Hamiltonian dissipates over time, up to the exchange with the boundary, i.e.,
\begin{align} \label{eq:energ-diss0}
	\frac{d}{dt}\tilde{\HamPDE}(\rho,v) = -\int_\Omega \tilde{r}(\rho,v) (\rho v)^2 dx + \bv u \cdot \TraceOp^{ }(\rho v) \leq \bv u \cdot \TraceOp^{ }(\rho v).
\end{align}
The energy dissipation (in-)equality can be derived from the variational principle stated in \eqref{eq:var}.

For a strong solution $(\rho,v)\in \funSpace{C}^1([0,T];\funSpace{C}_{pw}^1(\mathcal{E}) \times \funSpace{C}_{pw}^1(\mathcal{E}))$ of \eqref{bls-eq:abstr} and the mass flux $m=\rho v$, the variational principle
\begin{subequations}\label{eq:var}
\begin{align} 
 \Lscal \partial_t \rho, \tfa \Rscal &= - \Lscal \partial_x m, \tfa \Rscal \\
  \Lscal \partial_t \frac{m}{\rho}, \tfb \Rscal &= \Lscal P'(\rho)+ \frac{m^2}{2 \rho^2}, \partial_x \tfb \Rscal  - \Lscal r(\rho,m) m, \tfb \Rscal + \bv u \cdot \TraceOp \tfb
\end{align}
\end{subequations}
holds for all $\tfa \in \funSpace{L}^2(\mathcal{E})$, $\tfb \in \funSpace{H}_{div}^1(\mathcal{E})$, whereby $\bv u = [u_{\nu_1}, \ldots, u_{\nu_p}]^T$ and $r$ is given by $r(\rho,m) = \tilde{r}(\rho,m/\rho)$.
The principle follows from \eqref{bls-eq:abstr} by standard arguments and the coupling conditions inherited in $\funSpace{H}_{div}^1(\mathcal{E})$. 

\begin{rmrk}\label{rem:EulerArepH}
The dissipation (in-)equality \eqref{eq:energ-diss0} is a consequence of the more profound structural property that the barotropic Euler equations can be stated in {port-Hamiltonian} form. In particular, \eqref{bls-eq:abstr-a} can be formally written as
\begin{align*}
	\partial_t \rho = -\partial_x \nabla_v \tilde{\HamPDE}(\rho,v) , \qquad \partial_t v = -\partial_x \nabla_\rho \tilde{\HamPDE}(\rho,v) -\tilde{r}(\rho,v) \nabla_v \tilde{\HamPDE}(\rho,v),
\end{align*}
with $\nabla_\rho \tilde{\HamPDE}$ and $\nabla_v \tilde{\HamPDE}$ denoting the functional derivatives of $\tilde{\HamPDE}$ with respect to (w.r.t.)~the variable $\rho$ and $v$, respectively. The approximation procedure from \cite{art:lilsailer-nlfow}, which we employ in this work, is designed to preserve the {port-Hamiltonian} structure on the discrete level. It crucially relies on the formulation of the variational principle \eqref{eq:var} in terms of the mass flux $m$ while still using the Hamiltonian formulation with its symmetries, as it is stated here.
\end{rmrk}

\section{Structure-preserving approximation framework} \label{sec:approx-ansatz}

The structure-preserving approximation concept from \cite{art:lilsailer-nlfow, phd:liljegren} forms the theoretical basis of this work. We briefly present the concept and then discuss our specific choices for constructing the full order model and setting up the data-based reduction.

\subsection{Generic approximation procedure}
The structure-preserving approximation procedure developed in \cite{art:lilsailer-nlfow, phd:liljegren} is based on the variational formulation~\eqref{eq:var}.
It is applicable to Galerkin projections with classical spatial finite element discretization and model order reduction as well as to complexity reduction of quadrature-type for the nonlinear terms. On all approximation levels local mass conservation and an energy bound for the approximations as well as the port-Hamiltonian structure indicated in Remark~\ref{rem:EulerArepH}  are kept under mild assumptions.

\begin{assumption}[Compatibility of spaces] \label{assum:compatV1V2}
Let $\funSpace{V}= \Vsa \times \Vsb \subset \funSpace{L}^2(\mathcal{E})\times \funSpace{H}_{div}^1(\mathcal{E})$ be a finite dimensional subspace fulfilling the compatibility conditions
\begin{enumerate}[{A}1)]
\item $\Vsa =   \partial_x \Vsb, \qquad  \text{ with }  \partial_x \Vsb= \left \{\xi: \text{It exists } \zeta \in \Vsb \text{ with } \partial_x \zeta = \xi \right \} $,
\item $\funSpace{K} \subset \Vsb$, \hspace{0.9cm} \text{ with }  $\funSpace{K}=\{\tfb \in \funSpace{H}^{1}_{div}(\mathcal{E}): \, \partial_x \tfb = 0 \}$.
\end{enumerate}
\end{assumption}
The space $\funSpace{K}$ consists of the edge-wise constant functions fulfilling a certain coupling condition. Thus it only depends on the underlying topology and $\mydim(\funSpace{K})=|\mathcal{E}|-|\mathcal{N}_0|$.
\begin{assumption}[Compatibility of scalar product]  \label{assum:quadrat-ansatz}
Let the bilinear form $\Lscal \cdot , \cdot \Rscal_\scps: \funSpace{L}^2(\mathcal{E}) \times \funSpace{L}^2(\mathcal{E}) \rightarrow \mathbb{R}$ be such that the following holds:
\begin{enumerate}[{A}1)]
\item For a constant $\tilde{C}\geq 1$ and $||b||_\scps = \sqrt{\Lscal b,b \Rscal_{\scps}}$, it holds
${\tilde{C}}^{-1} ||b||_\scps \leq ||b|| \leq \tilde{C} ||b||_\scps$  for all $b \in \Vsa \cup \Vsb$.
\item For any $f \in \funSpace{C}_{pw}^0(\mathcal{E})$ with $f \geq 0$ it holds $\Lscal f , 1 \Rscal_\scps \geq 0$.
\end{enumerate}
\end{assumption}

\begin{system} \label{sys:abstract}
 Assumptions~\ref{assum:compatV1V2} and \ref{assum:quadrat-ansatz} are supposed to hold. Given $(\sfa_0,\sfb_0) \in \Vsa\times \Vsb$ and $\bv{u}: [0,T]\rightarrow \mathbb{R}^p$, find $(\sfa,\sfb) \in \funSpace{C}^1([0,T];\Vsa \times \Vsb) $ with $\sfa(0) = \sfa_0$, $\sfb(0) =\sfb_0$ and
\begin{align*} 
 \Lscal \partial_t \rho, \tfa \Rscal &= - \Lscal \partial_x m, \tfa \Rscal \\
  \Lscal \partial_t \frac{m}{\rho}, \tfb \Rscal_\scps &= \Lscal P'(\rho)+ \frac{m^2}{2 \rho^2}, \partial_x \tfb \Rscal_\scps  - \Lscal r(\rho,m) m, \tfb \Rscal_\scps + \bv u \cdot \TraceOp \tfb 
\end{align*}
for all $\tfa \in \Vsa$, $\tfb \in \Vsb$.  Further, $\HamPDE_\scps(\rho,m) = \Lscal P(\rho)+ {m^2}/({2 \rho}), 1 \Rscal_\scps$ is referred to as the Hamiltonian of the system.
\end{system}

\begin{thrm}[Energy bound, \cite{art:lilsailer-nlfow}] \label{thrm:energy-diss-cr}
The energy dissipation (in-)equality is fulfilled by any solution $(\rho,m)$ of System~\ref{sys:abstract},  
\begin{align*}
	\frac{d}{dt} \HamPDE_\scps(\rho(t),m(t))  &= \bv u(t) \cdot \TraceOp^{ }m(t) - \left\Lscal r(\rho(t),m(t)) m(t)^2, 1 \right\Rscal_\scps \leq \bv u(t) \cdot \TraceOp^{ }m(t). 
\end{align*}
\end{thrm}
Note that the inequality in the energy bound is implied by the non-negativity of the map $f \mapsto \Lscal f , 1 \Rscal_\scps $. 
Our flow problem possesses the form of System~\ref{sys:abstract} on each approximation level. At first, a Galerkin projection is constructed by mixed finite elements yielding the full order model ({FOM}). Subsequently, a second Galerkin projection is applied to obtain a reduced order model ({ROM}). Although the ROM is of much lower dimension than the FOM, it is unfortunately not more efficient, since the evaluation of the nonlinearities is not independent of the dimension of the {FOM}. Thus, a final complexity reduction step is performed that replaces the $\funSpace{L}^2$-scalar product by a cheaper-to-evaluate bilinear form for all nonlinear expressions yielding the complexity-reduced model ({CROM}).  In this sense, $\Lscal \cdot , \cdot \Rscal_\scps$ in System~\ref{sys:abstract} is chosen as a complexity-reduced bilinear form (specified and denoted as $\Lscal \cdot , \cdot \Rscal_c$ below) for the CROM, otherwise for FOM and ROM it is considered as the $\funSpace{L}^2$-scalar product $\Lscal \cdot, \cdot \Rscal$. In any case, $\Lscal \cdot , \cdot \Rscal_\scps$ is a scalar product and induces a norm on the spaces $\Vsa$ and $\Vsb$, respectively, by Assumption~\ref{assum:quadrat-ansatz}. The energy dissipation (in-)equality (Theorem~\ref{thrm:energy-diss-cr}) ensures the stability of the model on each approximation level. Apart from that, System~\ref{sys:abstract} inherits further structural properties. It preserves the  {port-Hamiltonian} structure mentioned in Remark~\ref{rem:EulerArepH}, and its solutions are locally mass conservative, i.e., $\partial_t \rho = -\partial_x m$ holds in a pointwise sense. We refer to \cite{art:lilsailer-nlfow, phd:liljegren} for details and proofs.

\subsection{Specific full order model and data-based reduction setting} \label{subsec:fom-training-prob}
Concerning the {FOM} we follow \cite{inproc:bls-hyp2018,art:lilsailer-dwemor} and use a mixed finite element space that is compatible in the sense of  Assumption~\ref{assum:compatV1V2}. A partitioning of the spatial domain $\Omega$ into the cells $K_1, \ldots, K_J$ is defined such that
\begin{subequations}\label{eq:FOM-FE}
\begin{align}
	\bigcup\limits_{j} K_j = \Omega, \qquad \int_{K_j \cap K_{l}} 1 dx = 0, \quad \text{ for } j \neq l, \quad K_j \subset \Onepipe \, \text{ for a } \Onepipe \in \mathcal{E} \text{ for } j=1,\ldots, J. \label{eq:FOM-FE-K}
\end{align}
The finite element ansatz space $ \funSpace{V}_{f}= \Vsa_{f} \times \Vsb_{f}$, which we also refer to as {FOM} space, is given by
\begin{align}
	\Vsa_{f} &= \left\{ \phi: \Omega_{} \rightarrow \mathbb{R}: \hspace{0.26cm} \phi_{|K_{j}}(x) = \zeta_{j}, \hspace{1.25cm} \text{ with } \zeta_{j}\in \mathbb{R} \hspace{0.7cm}\text{ for } j=1,\ldots,J\right\} \\
	\Vsb_{f} &= \left\{ \phi\in \funSpace{H}^1_{div}(\mathcal{E}): \phi_{|K_{j}}(x) = \xi_{j} + \mu_j x, \hspace{0.25cm} \text{ with } \xi_{j},\mu_j \in \mathbb{R} \hspace{0.2cm}\text{ for } j=1,\ldots,J\right\}. \label{eq:FOM-FE-V2}
\end{align}
\end{subequations}
The functions in $\Vsa_{f}$ are piecewise constant and discontinuous, whereas $\Vsb_{f}$ consists of piecewise linear, edgewise continuous functions and inherits the coupling conditions that relate to mass conservation at inner nodes. 

The reduced models are trained towards data, the so-called snapshots, given as
\begin{align} \label{eq:snap-statespace}
	\ubv{S} = \{ \ubv{a}^1, \ldots , \ubv{a}^L \}, \hspace{1cm}
	\ubv{a}^\ell = (\rho^\ell ,
		m^\ell) \in \funSpace{V}_f = \Vsa_f \times \Vsb_f, \quad \ell =1,\ldots , L.
\end{align}
The snapshots are typically selected (time) points of one or several solution trajectories of the {FOM} for appropriate training scenarios. In the model order reduction we seek for a low-dimensional subspace $\funSpace{V}_r = \Vsa_r \times \Vsb_r$ of $\funSpace{V}_f$, $\mydim(\funSpace{V}_r)\ll \mydim(\funSpace{V}_f)$, that fulfills Assumption~\ref{assum:compatV1V2} and in which the snapshot data is approximated with high fidelity. The subsequent complexity reduction aims for constructing a bilinear form $\Lscal \cdot , \cdot \Rscal_c$ that fulfills Assumption~\ref{assum:quadrat-ansatz}, whose evaluation-cost is  independent of $\mydim(\funSpace{V}_f)$, and for which $\Lscal f(\ubv{a}^\ell), \tfb_r \Rscal_c \approx \Lscal f(\ubv{a}^\ell), \tfb_r \Rscal$ holds for $\ell = 1, \ldots, L$ and $\tfb_r \in \Vsb_r$ for the appearing nonlinear terms $f$. We use a quadrature-type approximation.

Our approach differs from most conventional data-based model reduction methods, e.g, \cite{book:dimred2003,art:pod-kunisch,art:deim-state-space-err}, in that compatibility conditions are imposed to ensure the preservation of structural properties for the reduced systems, such as an energy bound (Theorem~\ref{thrm:energy-diss-cr}) or well-posedness (Section~\ref{sec:well-posed}). In the training phase, the conditions pose an additional challenge, as they lead to constraints.

\section{Realization of model order reduction} \label{sec:prop-mr}

In this section we formulate the data-based model reduction problem with compatibility conditions. It can be considered as a principal component analysis (PCA) with constraints. As will be shown, this constrained problem can be attributed to an underlying standard PCA given suitable norms are chosen. Based on that, we derive an efficient algorithmic implementation.

\subsection{Constrained principal component analysis} \label{subsec:pod-approx-cond} 

Our model order reduction follows the basic notion of proper orthogonal decomposition. The {ROM} space $\funSpace{V}_r$ is obtained by a {PCA} of the snapshot data, but in contrast to the standard approach, constraints related to the compatibility conditions are posed. To be able to treat the constraints efficiently, the choice of a proper scalar product $\Lscal \cdot , \cdot \Rscal_\diamond$ for $\Vsb_{f}$ in the {PCA}  will be crucial. We denote the respective induced norm by  $||\cdot||_\diamond$, i.e., $||\tfb||_\diamond= \sqrt{\Lscal \tfb , \tfb \Rscal_\diamond}$. Moreover, the orthogonal projection onto a subspace $\funSpace{U}_{}$ w.r.t.\ the scalar product $\Lscal \cdot , \cdot \Rscal$ and $\Lscal \cdot , \cdot \Rscal_\diamond$ is denoted by $\Pfun_{\funSpace{U}_{}}$ and $\Pfun_{\funSpace{U}_{}}^\diamond$, respectively. The general approximation problem we consider for the determination of the {ROM} space $\funSpace{V}_r$ reads as follows. 

\begin{prblm} \label{prblm:PCAstart}
Let $\funSpace{V}_f = \Vsa_{f}\times \Vsb_{f}$ be a {FOM} space compatible in the sense of Assumption~\ref{assum:compatV1V2}. For the snapshots $(\rho^\ell, m^\ell)$ as in \eqref{eq:snap-statespace} and $n_1 \leq N_1=\mydim(\Vsa_f)$, find a reduced space $\funSpace{V}_{r}= \Vsa_{r} \times \Vsb_{r} \subset \funSpace{V}_f$ as solution to
\begin{align*}
		\min_{\funSpace{V} \subset \funSpace{V}_f} 	\,\, & \sum_{\ell=1}^L  \left|\left| \sfa^\ell -  \Pfun_{\Vsa} \sfa^\ell \right|\right|^2 + \left|\left| \sfb^\ell -  \Pfun^\diamond_{\Vsb} \sfb^\ell \right|\right|_\diamond^2 \\
		\text{s.t.} & \,\, \funSpace{V}=\Vsa\times \Vsb   \text{ fulfills Assumption~\ref{assum:compatV1V2}} \\
		& \,\, \mydim(\Vsa) \leq n_1 .
\end{align*}
\end{prblm}
Due to the compatibility conditions (Assumption~\ref{assum:compatV1V2}), Problem~\ref{prblm:PCAstart} is a constrained optimization problem which might be quite involved to solve. We aim for the derivation of an equivalent unconstrained standard PCA. To remove the constraints, the underlying pure data approximation problem has to be separated from the question of compatibility. 

\begin{lmm}
	Let  $\funSpace{K}=\{\tfb \in \funSpace{H}^{1}_{div}(\mathcal{E}): \, \partial_x \tfb = 0 \}$ be as in Assumption~\ref{assum:compatV1V2}. Then the bilinear form 
	\begin{align}
		\Lscal \cdot , \cdot \Rscal_\diamond : \funSpace{H}^1_{div}(\mathcal{E}) \times \funSpace{H}^1_{div}(\mathcal{E}) \rightarrow \mathbb{R}, \qquad 
		\Lscal \tfb, \tilde\tfb \Rscal_\diamond = \Lscal \Pfun_{\funSpace{K}} \tfb, \Pfun_{\funSpace{K}} \tilde\tfb \Rscal + \Lscal \partial_x \tfb, \partial_x \tilde\tfb \Rscal \label{eq:scalarproduct}
	\end{align}
is a scalar product on $\funSpace{H}^1_{div}(\mathcal{E})$. 
\end{lmm}

\begin{proof}
Using that $\funSpace{K}$ is the kernel of $\partial_x: \funSpace{H}^1_{div}(\mathcal{E})\rightarrow \funSpace{L}^2(\mathcal{E}) $, the assertion is straightforward to show. We refer to \cite{art:kugler-dwe-net}, where the norm $|| \cdot ||_\diamond$ is considered. 
\end{proof}

In view of Assumption~\ref{assum:compatV1V2}, we introduce the orthogonal decomposition of $\Vsb_{f}$ w.r.t.~$\Lscal \cdot , \cdot \Rscal_\diamond$ into $\funSpace{K}$ and its orthogonal complement $\funSpace{K}^\perp$,
\begin{align*}
	\Vsb_{f} = \funSpace{K}^\perp \oplus   \funSpace{K}, \quad \text{i.e.,} \quad  \tfb = \Pfun^\diamond_{\funSpace{K}^\perp} \tfb +\Pfun^\diamond_{\funSpace{K}} \tfb \quad \text{and} \quad \Lscal \Pfun^\diamond_{\funSpace{K}^\perp} \tfb ,  \Pfun^\diamond_{\funSpace{K}} \tfb \Rscal_\diamond= 0 \qquad \text{for } \tfb \in  \Vsb_{f}.
\end{align*}
Furthermore, we define $\partial_x^+:\Vsa_{f} \rightarrow \Vsb_{f}$ as the right-inverse of the spatial derivative $\partial_x: \Vsb_{f} \rightarrow \Vsa_{f}$, which additionally fulfills $\Pfun^\diamond_{\funSpace{K}} \partial_x^+  \equiv 0$.
Note that $\funSpace{K}$ is the kernel of $\partial_x$, and that  $\partial_x^+$ is well-defined and relates to a weighted Moore-Penrose inverse with weights given by the scalar product $\Lscal \cdot , \cdot \Rscal_\diamond$, cf.\ \cite{book:generalizedinverses-israel}. First, we show that the constraints can be formally removed for any choice of $\Lscal \cdot , \cdot \Rscal_\diamond$.

\begin{thrm}\label{theor:PCAgeneral}
Let the assumptions of Problem~\ref{prblm:PCAstart} hold. Then $\funSpace{V}_{r}=\Vsa_{r} \times \Vsb_{r} \subset \funSpace{V}_f$ is a solution to Problem~\ref{prblm:PCAstart}, if and only if $\Vsa_{r}$ is a solution to 
\begin{align*}
\min_{\substack{\Vsa\subset \Vsa_{f}\\ \hspace{0.48cm} \mydim(\Vsa) \leq n_1 }}
 \sum_{\ell=1}^L  \left|\left| \sfa^\ell -  \Pfun_{\Vsa} \sfa^\ell \right|\right|^2 + \left|\left| \Pfun^\diamond_{\funSpace{K}^\perp} \sfb^\ell -  \Pfun^\diamond_{\partial_x^+ \Vsa} \Pfun^\diamond_{\funSpace{K}^\perp}\sfb^\ell \right|\right|_\diamond^2 ,
\end{align*}
and it holds $\Vsb_{r} = \partial_x^+ \Vsa_{r} \oplus \funSpace{K}$.
\end{thrm}

\begin{proof}
Let $\Vsa \times \Vsb \subset \Vsa_{f} \times \Vsb_{f}$ be a compatible reduced space in the sense of Assumption~\ref{assum:compatV1V2}. Analogously as for the {FOM}, an orthogonal decomposition of $\Vsb$ w.r.t. $\Lscal \cdot , \cdot \Rscal_\diamond$ can be performed, which yields $\Vsb = (\Vsb \cap \funSpace{K}^\perp ) \oplus \funSpace{K}$. From $\partial_x \Vsb = \Vsa $ it then follows that $\Vsb = \partial_x^+ \Vsa \oplus \funSpace{K}$.

The orthogonal decomposition of the FOM implies $ || \tfb ||_\diamond^2 = || \Pfun^\diamond_{\funSpace{K}^\perp} \tfb ||_\diamond^2 + || \Pfun^\diamond_{\funSpace{K}} \tfb ||_\diamond^2$ for any $\tfb \in \Vsb_{f}$. Using this equality for the approximation error together with the compatibility $\funSpace{K} \subset \Vsb$, it follows $||\tfb-  \Pfun^\diamond_{\Vsb} \tfb||_\diamond^2 = ||\Pfun^\diamond_{\funSpace{K}^\perp} \tfb -  \Pfun^\diamond_{\Vsb\cap\funSpace{K}^\perp} \Pfun^\diamond_{\funSpace{K}^\perp} \tfb ||_\diamond^2$. All in all, the compatibility conditions of Assumption~\ref{assum:compatV1V2} can be removed, and the unconstrained optimization problem is equivalent to Problem~\ref{prblm:PCAstart} in the sense of the theorem.
\end{proof}

Among others, the theorem shows that the snapshots $\sfb^\ell$ should not be used directly in the {PCA} but should be altered first to $\Pfun^\diamond_{\funSpace{K}^\perp} \sfb^\ell$. Moreover note that the sub-problem determining $\Vsa_r$ is unconstrained. However, its cost functional inherits terms with different norms  depending on the choice of scalar product $\Lscal \cdot , \cdot \Rscal_\diamond$. For the specific choice \eqref{eq:scalarproduct}, we can attribute the following standard {PCA} to the training problem.

\begin{prblm} \label{prblm:PCAred}
Let $\funSpace{V}_f = \Vsa_{f}\times \Vsb_{f}$ be a {FOM} space compatible in the sense of Assumption~\ref{assum:compatV1V2}. For the snapshots $(\rho^\ell, m^\ell)$ as in \eqref{eq:snap-statespace} and $n_1 \leq N_1=\mydim(\Vsa_f)$, find a reduced space $\Vsa_{r}  \subset \Vsa_{f}$ as solution of
\begin{align*}
		\min_{\substack{\Vsa\subset \Vsa_{f}\\ \hspace{0.48cm} \mydim(\Vsa) \leq n_1 }}	 & \sum_{\ell=1}^L  \left|\left| \sfa^\ell -  \Pfun_{\Vsa} \sfa^\ell \right|\right|^2 + \left|\left| \partial_x \sfb^\ell -  \Pfun_{\Vsa} \partial_x \sfb^\ell \right|\right|^2 .
\end{align*}
\end{prblm}

\begin{thrm} \label{theor:easier-problem}
Let $\Lscal \cdot , \cdot \Rscal_\diamond$ be the scalar product defined in \eqref{eq:scalarproduct}. Then, $\funSpace{V}_r=\Vsa_{r} \times \Vsb_{r} \subset \funSpace{V}_f$ solves Problem~\ref{prblm:PCAstart}, if and only if $\Vsa_{r}$ solves Problem~\ref{prblm:PCAred} and $\Vsb_{r} = \partial_x^+ \Vsa_{r} \oplus \funSpace{K}$.
\end{thrm}
\begin{proof}
As auxiliary result we show that for any compatible pair of spaces $\Vsa$, $\Vsb$ it holds
\begin{align} \label{eq:prblm:PCAred}
	\partial_x \Pfun_{\partial_x^+ \Vsa}^{\diamond}   \Pfun_{\funSpace{K}^\perp}^\diamond  \tfb =    \Pfun_{\Vsa}  \partial_x \tfb \hspace{1cm} \text{ for } \tfb \in \Vsb.
\end{align}
By the definition of orthogonal projections, $w_2 = \Pfun_{\partial_x^+ \Vsa}^{\diamond}   \Pfun_{\funSpace{K}^\perp}^\diamond \tfb$ for $\tfb \in \Vsb$ is fully characterized by the two conditions $w_2 \in \partial_x^+ \Vsa$ and $\Lscal w_2 , \tilde\tfb \Rscal_\diamond = \Lscal \tfb, \tilde\tfb \Rscal_\diamond$ for all $\tilde\tfb \in \partial_x^+ \Vsa$.
From $\Vsb \cap \funSpace{K}^\perp = \partial_x^+ \Vsa$ it follows that
\begin{align*}
	\Lscal w_2 , \tilde\tfb \Rscal_\diamond = \Lscal \Pfun_{\funSpace{K}^\perp}^\diamond  \tfb, \tilde\tfb \Rscal_\diamond \quad \text{ for } \tilde\tfb \in \partial_x^+ \Vsa  \quad& \Longleftrightarrow \quad
	\Lscal w_2 , \partial_x^+ \tfa \Rscal_\diamond = \Lscal \Pfun_{\funSpace{K}^\perp}^\diamond  \tfb, \partial_x^+ \tfa \Rscal_\diamond && \text{ for } \tfa \in \Vsa  \\
	& \Longleftrightarrow \quad 
	\Lscal \partial_x w_2 , \tfa \Rscal \hspace{0.24cm} = \Lscal \partial_x  \Pfun_{\funSpace{K}^\perp}^\diamond  \tfb, \tfa \Rscal_\diamond  = \Lscal \partial_x  \tfb, \tfa \Rscal_\diamond && \text{ for } \tfa \in \Vsa.
\end{align*}
The last equivalence makes use of $\funSpace{K}$ being the kernel of $\partial_x$. The full equivalence together with the fact $\partial_x w_2 \in \Vsa$ yields that $\partial_x w_2$ is the orthogonal projection of $\partial_x \tfb $ onto $\Vsa$, i.e., the result \eqref{eq:prblm:PCAred}.

To show the equivalence of the minimization problems, we make use of the representation given in Theorem~\ref{theor:PCAgeneral}. Its cost function includes terms of the form $ || \Pfun^\diamond_{\funSpace{K}^\perp} \sfb^\ell -  \Pfun^\diamond_{\partial_x^+ \Vsa} \Pfun^\diamond_{\funSpace{K}^\perp}\sfb^\ell ||_\diamond$. For our specific choice of inner product, it holds 
\begin{align*}
\left|\left| \Pfun^\diamond_{\funSpace{K}^\perp} \sfb^\ell -  \Pfun^\diamond_{\partial_x^+ \Vsa} \Pfun^\diamond_{\funSpace{K}^\perp}\sfb^\ell \right| \right|_\diamond &=
\left|\left| \partial_x \Pfun^\diamond_{\funSpace{K}^\perp} \sfb^\ell -  \partial_x \Pfun^\diamond_{\partial_x^+ \Vsa} \Pfun^\diamond_{\funSpace{K}^\perp}\sfb^\ell \right| \right| = 
\left|\left| \partial_x \sfb^\ell -  \Pfun_{\Vsa}  \partial_x \sfb^\ell \right| \right|,
\end{align*}
whereby \eqref{eq:prblm:PCAred} has been used in the last step. From that the equivalence of the cost function of Problem~\ref{prblm:PCAred} and the one stated in Theorem~\ref{theor:PCAgeneral} can be concluded in a straightforward manner. Using Theorem~\ref{theor:PCAgeneral} finishes the proof.
\end{proof}

\subsection{Structured representations} \label{subsec:alg-real}

For the algorithmic realization we introduce coordinate representations of the full and reduced order models.

\subsubsection{Full order model}
Consider the FOM space $\funSpace{V}_f = \Vsa_{f} \times \Vsb_{f}$ with $N = \mydim(\funSpace{V}_f) = N_1+N_2$. Let $\{\tfa^1,\ldots, \tfa^{N_1}\}$ and $\{\tfb^1,\ldots, \tfb^{N_2}\}$ be bases for $\Vsa_f$ and $\Vsb_f$, respectively. The bijective mapping between the coordinate representation $\bv{a}= [\bv{a}_1^T, \bv{a}_2^T]^T \in \mathbb{R}^N$, $\bv{a}_i= [\scalco{a}_i^1, \ldots , \scalco{a}_i^{N_i}]^T \in \mathbb{R}^{N_i}$, and the function $(\sfa,\sfb) \in \funSpace{V}_f$ is given by
\begin{align*}
	\Psi : \mathbb{R}^N \rightarrow \funSpace{V}_f, \qquad  \Psi(\bv{a}) =
	\left(
	 \sum_{j=1}^{N_1} \tfa^j \scalco{a}_1^j	\, , \, 	 \sum_{j=1}^{N_2} \tfb^j \scalco{a}_2^j \right) = (\sfa,\sfb).
\end{align*} 

The mass matrices for $\Vsa_f$ and $\Vsb_f$ as well as the coordinate representation of $\partial_x: \Vsb_{f} \rightarrow \Vsa_{f}$ are defined by
\begin{align} \label{eq:coord-fom-mat}
	\Massa &= \left[ \Lscal \tfa^n,\tfa^m \Rscal \right]_{m,n=1,\ldots, N_1}, \hspace{0.5cm} \Massb = \left[ \Lscal \tfb^n, \tfb^m \Rscal  \right]_{m, n = 1,\ldots, N_2}, \hspace{0.5cm} \bt{J} = \left[ \Lscal \partial_x \tfb^n,\tfa^m \Rscal  \right]_{m=1,\ldots, N_1, n = 1,\ldots, N_2} .
\end{align}
The boundary operator $\TraceOp: \Vsb_f \rightarrow \mathbb{R}^p$ has the coordinate representation $\bt{B} = \left[\TraceOp \tfb^1,\ldots , \TraceOp \tfb^{N_2} \right]^T$. 
Moreover, given $(\sfa,\sfb) = \Psi(\bv a)\in \funSpace{V}_f$, we define vectors related to the nonlinear expressions by
\begin{align}
	 \bv f^\alpha(\bv a) = \left[ \Lscal \frac{m}{\rho} , \tfb^j \Rscal  \right]_{j=1,\ldots, N_2}\hspace{-0.1cm},  \hspace{0.4cm}\bv f^\beta(\bv a) = \left[ \Lscal P'(\rho)+\frac{m^2}{2\rho^2},\partial_x \tfb^j \Rscal  \right]_{j=1,\ldots, N_2}\hspace{-0.1cm}, \hspace{0.4cm}
	\bv f^\gamma(\bv a) = \left[ \Lscal -r(\rho,m) m, \tfb^j \Rscal  \right]_{j=1,\ldots, N_2}. \label{eq:coord-fom-nl}
\end{align}
Then the {FOM} (System~\ref{sys:abstract}) can be equivalently described by the following algebraic representation.
\begin{system}[Algebraic representation of {FOM}] \label{sys:coord-rep}
Given appropriate $\bv a_0 \in \mathbb{R}^{N}$ and $\bv{u}: [0,T]\rightarrow \mathbb{R}^p$, find $\bv{a} = [\bv{a}_1^T,\bv{a}_2^T]^T \in \funSpace{C}^1([0,T];\mathbb{R}^N)$ with $\bv a(0) = \bv a_0$ and
\begin{align*} 
\Massa \frac{d}{dt} \bv{a}_1(t) = -\bt{J} \bv{a}_2(t), \hspace{1cm} \frac{d}{dt} \bv f^\alpha(\bv a(t)) = \bv f^\beta(\bv a(t)) + \bv f^\gamma(\bv{a}(t))+ \bt B \bv u(t).
\end{align*}
\end{system}

\begin{rmrk}
Related to System~\ref{sys:coord-rep}, the Hamiltonian $H$ can be defined as $H(\bv{a}) = \HamPDE(\Psi(\bv{a}))$. Using this Hamiltonian, the system can be transformed into standard {port-Hamiltonian} form, see \cite{art:lilsailer-nlfow}.
\end{rmrk}

For $m \in \Vsb_f$,  $m = \sum_{j=1}^{N_2} \tfb^j \scalco{a}_2^j$ with $\bv a_2 = [ \scalco a_2^1 ,\ldots , \scalco a_2^{N_2} ]^T \in \mathbb{R}^{N_2}$, the coordinate representation of $\partial_x m $ in the basis of $\Vsa_f$ is given by
\begin{align} \label{eq:dx-coord}
\partial_x m = \sum_{i=1}^{N_i} \tfa^i \scalco{a}_1^i, \hspace{1cm} \left[ \scalco a_1^1 ,\ldots , \scalco a_1^{N_1} \right]^T = \Massa^{-1} \bt J \bv a_2 \in \mathbb{R}^{N_1},
\end{align}
which can be concluded from Assumption~\ref{assum:compatV1V2} ($\partial_x\Vsb_f=\Vsa_f$).

\subsubsection{Reduced order models} 
The projection of the {FOM} onto the low-dimensional {ROM} space $\funSpace{V}_r \subset \funSpace{V}_f$ can be equivalently described in the algebraic setting as projecting System~\ref{sys:coord-rep} with a reduction matrix $\bt{V} \in \mathbb{R}^{N,n}$ and $n=\mydim(\funSpace{V}_r) \ll N$.
As we consider projections orthogonal w.r.t.~the $\funSpace{L}^2$-inner product, weighted orthogonal projections occur naturally in the coordinate representation. Accordingly the following definitions are helpful, cf., \cite{book:generalizedinverses-israel}. 
Given a symmetric positive definite matrix $\bt{M} \in \mathbb{R}^{N,N}$, we introduce the scalar product $\langle \cdot, \cdot \rangle_{\bt{M}}$ and the respective norm $\|\cdot \|_{\bt M}$, i.e., $\langle \bv x, \bv y \rangle_{\bt M}=\bv x^T \bt M \bv y$ and $\|\bv x \|_{\bt M}=\sqrt{\langle \bv x, \bv x \rangle_{\bt M}}$ for any $\bv x, \bv y \in \mathbb{R}^N$.
For any matrix $\bt{V} \in \mathbb{R}^{N,n}$ of full column rank, we refer to
\begin{align}\label{eq:weighted inverse}
	 \Pvec_{\bt{V}}^{\bt{M}}  =\bt{V} \bt{V}^{+,\bt{M}} = \bt{V} ( \bt{V}^T \bt{M} \bt{V} )^{-1} \bt{V}^T \bt{M}
\end{align}
as the weighted orthogonal projection w.r.t.\ $\langle \cdot, \cdot \rangle_{\bt{M}}$, where $\bt{V}^{+,\bt{M}}$ particularly denotes the weighted left-inverse of $\bt{V}$. Note that $\Pvec_{\bt{V}}^{\bt{M}} \Pvec_{\bt{V}}^{\bt{M}} = \Pvec_{\bt{V}}^{\bt{M}} $ and $\myim{\Pvec_\bt{V}^{\bt{M}}} = \myim{\bt{V}}$ as for any projection onto $\myim{\bt{V}}$ (image of $\bt{V}$). Moreover, the weighted orthogonality condition $\bt{V}^T \bt{M} (\bt{I}- \Pvec_{\bt{V}}^{\bt{M}}) = \bt{0}$ holds.
In case of $\bt{M}=\bt I$ identity matrix (implying the Euclidean scalar product and norm), we suppress the index $\bt{M}$. 

\begin{lmm}\label{lem:compatMatrix}
Let $\funSpace{V}_f$ fulfill Assumption~\ref{assum:compatV1V2} and $\Massa$, $\bt J$ be as in \eqref{eq:coord-fom-mat}. For a block-structured reduction matrix $\bt V$ 
\begin{align*}
	\bt V  = \begin{bmatrix}
	\bt V_1 & \\
			& \bt V_2
	\end{bmatrix}, \qquad \bt V_i \in \mathbb{R}^{N_i,n_i}, \quad n_1 = \mydim(\Vsa_{r}), \quad n_2 = \mydim(\Vsb_{r}),
\end{align*}
the space $\funSpace{V}_r =\Vsa_{r}\times \Vsb_{r}:= \myim{ \bt V_1} \times \myim{ \bt V_2}$ fulfills  Assumption~\ref{assum:compatV1V2}, if and only if  $\myim{\Massa \bt V_1} =   \myim{\bt J \bt V_2}$ and $\mathrm{ker}(\bt J)\subset \myim{\bt V_2}$ hold.
\end{lmm}
The compatibility condition of Assumption~\ref{assum:compatV1V2} on $\funSpace{V}_r$ can be recast as a compatibility condition on the reduction matrix $\bt V$ according to Lemma~\ref{lem:compatMatrix}, as can be derived with the help of \eqref{eq:dx-coord}, cf., \cite{art:lilsailer-dwemor,phd:liljegren}.
 The {ROM} and a reduced coordinate representation $\bv{a}_r$, $\bt V \bv{a}_r \approx \bv a$, are then characterized by a projected version of System~\ref{sys:coord-rep}. Its initial conditions $\bv{a}_r(0)\in \mathbb{R}^{n_1+n_2}$ are chosen as the $\funSpace{L}^2$-projections of $\bv{a}_0=[\bv{a}_{0,1}^T,\bv{a}_{0,2}^T]^T\in \mathbf{R}^{N_1+N_2}$, accordingly.
\begin{system}[Algebraic representation of {ROM}] \label{sys:rom-coord}
Given System~\ref{sys:coord-rep} and a reduction matrix $\bt V$ as in Lemma~\ref{lem:compatMatrix}, find $\bv{a}_r = [\bv{a}_{r,1}^T,\bv{a}_{r,2}^T]^T \in \funSpace{C}^1([0,T];\mathbb{R}^n)$, with $\bv a_r(0) = [(\bt V_1^{+,\Massa} \bv a_{0,1})^T,(\bt V_2^{+,\Massb} \bv a_{0,2})^T]^T $ and
\begin{align*} 
\Massa_{r} \frac{d}{dt} \bv{a}_{r,1}(t) = -\bt{J}_r \bv{a}_{r,2}(t), \hspace{1cm} \frac{d}{dt} \bv f^\alpha_{r}(\bv a_r(t)) = \bv f^\beta_{r}(\bv a_r(t)) + \bv f^\gamma_r(\bv{a}_r(t))+ \bt B_r \bv u(t),
\end{align*}
where the reduced state matrices are $\Massa_{r} = \bt V_1^T \Massa \bt V_1$, $\bt{J}_r = \bt V_1^T \bt{J} \bt V_2$ and $\bt B_{r} = \bt V_2^T \bt B$. The reduced nonlinearities are defined by $\bv f_r(\bv a_r) = \bt V_2^T \bv f(\bt V \bv a_r)$ for $\bv f \in \{ \bv f^\alpha,  \bv f^\beta,  \bv f^\gamma \}$.
\end{system}

\subsection{Computation of structure-preserving reduction basis} \label{subsec:pod-approx-sol}
To reformulate a {{PCA}} problem in a weighted norm in terms of the Euclidean norm, we make use of the following result, cf. \cite{art:increm-POD}.
Let  $\bt{M} \in \mathbb{R}^{N,N}$ be a symmetric positive definite matrix and $\bt{M} = \bt{L}^T\bt{L}$. Then it holds 
\begin{align}\label{eq:rescalePCA}
		\|  \bv{x} -  \Pvec_{\bt{V}}^{\bt{M}}   \bv{x}  \|_\bt{M} = 
		\|  \bt{L}  \bv{x} -  \Pvec_{\bt{LV}} \bt{L}  \bv{x}  \|, \qquad \text{for } \bt{V} \in \mathbb{R}^{N,n}, \quad \bv{x} \in \mathbb{R}^N.
\end{align}
The relation follows from $\|   \bv{x} -  \Pvec_{\bt{V}}^{\bt{M}}   \bv{x}  \|_\bt{M} = 	\|  \bt{L} ( \bv{x} -  \Pvec_{\bt{V}}^{\bt{M}}   \bv{x} ) \| =	\|  \bt{L}  \bv{x} -  \Pvec_{\bt{L}\bt{V}}  \bt{L}  \bv{x}  \|$, as $\Pvec_{\bt{V}}^{\bt{M}} = \bt{L}^{-1} \Pvec_{\bt{L}\bt{V}} \bt{L}$ according to the definition of the weighted orthogonal projection \eqref{eq:weighted inverse}.

\begin{thrm}\label{theor:V1}
Let $\Massa$, $\bt J$ and $\Psi$ be given as in Section~\ref{subsec:alg-real}. Let $\bt{L } \in \mathbb{R}^{N_1,N_1}$ be such that $\Massa = \bt{L }^T \bt{L }$. Further, let the coordinate representations of the snapshots $(\rho^\ell , m^\ell) \in \funSpace{V}_f$ in \eqref{eq:snap-statespace} be given by
\begin{align*}
	\begin{bmatrix}
		\bv a_1^\ell \\
		\bv a_2^\ell
	\end{bmatrix} = \Psi^{-1}(\rho^\ell , m^\ell), \quad \ell =1,\ldots , L, \qquad \text{and} \qquad 
\bt S_i  = [\bv a_i^1, \ldots, \bv a_i^L] \in \mathbb{R}^{N_i,L}, \qquad   i =1,2
\end{align*}
with  $n_1 \leq \mathrm{rank}(\bt{K})$, $\bt{K} = \bt{L } [ \bt{S}_1 , \Massa^{-1} \bt{J} \bt{S}_2]$.
Then Problem~\ref{prblm:PCAred} can be equivalently described by
\begin{align*}
		\min_{\hspace{0.4cm} \bt{V}_1 \in\, \mathbb{R}^{N_1,n_1}}  \, & \sum_{\ell=1}^L  \| \bv{a}_1^\ell -  \Pvec_{\bt{V}_1}^{\Massa} \bv{a}_1^\ell \|_{\Massa}^2 
		+ \| \Massa^{-1} \bt{J} \bv{a}_2^\ell -  \Pvec_{\bt{V}_1}^{\Massa} \Massa^{-1}\bt{J} \bv{a}_2^\ell \|_{\Massa}^2 .
\end{align*}
Moreover, solutions $\bt{V}_1^*$ to this problem fulfill $\myim{\bt{V}_1^*} = \myim{ \bt{L}^{-1}[\tilde{\bv{v}}_1,\ldots,\tilde{\bv{v}}_{n_1} ]}$, where the vectors $\tilde{\bv{v}}_i$ denote the first left-singular vectors of the matrix $\bt{K}$.
\end{thrm}
\begin{proof}
The stated minimization problem relates to a coordinate representation of Problem~\ref{prblm:PCAred} in the basis $\{\tfa^1, \ldots, \tfa^{N_1}\}$ of $\Vsa_f$. This can be seen from the equalities
\begin{align*}
	\left| \left| \sfa^\ell -  \Pfun_{\Vsa} \sfa^\ell \right|\right| &=  \| \bv{a}_1^\ell -  \Pvec_{\bt{V}_1}^{\Massa} \bv{a}_1^\ell \|_{\Massa} \\
	\left|\left| \partial_x \sfb^\ell -  \Pfun_{\Vsa} \partial_x \sfb^\ell \right|\right| &= \| \Massa^{-1} \bt{J} \bv{a}_2^\ell -  \Pvec_{\bt{V}_1}^{\Massa} \Massa^{-1}\bt{J} \bv{a}_2^\ell \|_{\Massa}, \qquad \text{for } \ell \in 1, \ldots , L,
\end{align*}
using \eqref{eq:dx-coord}. Let 
$\bv{s}^\ell \in \mathbb{R}^{N_1}$ for $\ell = 1,\ldots, 2L$ be defined by $[\bv{s}^1,\ldots, \bv{s}^{2L}]= [\bt{S}_1 , \Massa^{-1}\bt{J} \bt{S}_2]$. Then it holds
\begin{align*}
\sum_{\ell=1}^L  \| \bv{a}_1^\ell -  \Pvec_{\bt{V}_1}^{\Massa} \bv{a}_1^\ell \|_{\Massa}^2 +
 \| \Massa^{-1} \bt{J} \bv{a}_2^\ell -  \Pvec_{\bt{V}_1}^{\Massa} \Massa^{-1}\bt{J} \bv{a}_2^\ell \|_{\Massa}^2 &=
 \sum_{\ell=1}^{2L}  \| \bv{s}^\ell -  \Pvec_{\bt{V}_1}^{\Massa} \bv{s}^\ell \|_{\Massa}^2 =
\sum_{\ell=1}^{2L}  \| (\bt{L } \bv{s}^\ell) -  \Pvec_{\bt{L }\bt{V}_1} (\bt{L } \bv{s}^\ell ) \|^2,
\end{align*}
according to \eqref{eq:rescalePCA}. Hence, $\bt{V}_1^*$ solves the problem of the theorem, if and only if $\tilde{\bt{V}}_1^* = \bt{L } \bt{V}_1^*$ solves
\begin{align*}
		\min_{\hspace{0.4cm} \tilde{\bt{V}}_1 \in\, \mathbb{R}^{N_1,n_1}}  \sum_{\ell=1}^{2L}  \| \tilde{\bv{s}}^\ell -  \Pvec_{\tilde{\bt{V}}_1} \tilde{\bv{s}}^\ell \|^2, \hspace{0.6cm}
		\text{ with }\tilde{\bv{s}}^\ell = \bt{L } \bv{s}^\ell =
		\begin{cases}
		 \bt{L } \bv a_1^\ell & 1 \leq \ell \leq L \\
		 \bt{L } \Massa^{-1} \bt J \bv a_2^{\ell-L}  & L+1 \leq \ell \leq 2 L.
		\end{cases}
\end{align*}
The latter is a standard {PCA} problem and can be solved by the method of snapshots, cf. \cite{art:pod-kunisch,art:morKunV01}. It has a solution $\tilde{\bt{V}}_1^* = [\tilde{\bv{v}}_1,\ldots,\tilde{\bv{v}}_{n_1}]$ with $\tilde{\bv{v}}_i$ as in the theorem. Consequently,  $\bt{V}_1^*$ fulfills $\myim{\bt{V}_1^*} =\myim{ \bt L^{-1}[\tilde{\bv{v}}_1,\ldots,\tilde{\bv{v}}_{n_1} ]}$.
\end{proof}

In the following algorithm we summarize our procedure for constructing a reduction matrix $\bt V$, which yields an optimal reduction space in the sense of Problem~\ref{prblm:PCAstart}. 
As the characterizations of  $\funSpace{V}_r=\Vsa_r\times \Vsb_r$ in Lemma~\ref{lem:compatMatrix} and Theorem~\ref{theor:PCAgeneral} suggest, the procedure separates into two parts. First, the space $\Vsa_r= \myim{\bt{V}_1}$ is determined according to Theorem~\ref{theor:V1}. Then, $\Vsb_r =\myim{\bt{V}_2}$ is derived using the compatibility conditions.

\begin{samepage}
\begin{lgrthm}[Optimal compatible reduction basis]\text{$ $} \label{alg:CPCA}\\
INPUT:
\begin{itemize}
\item Snapshot matrices for density and mass flux: $\bt{S}_1 \in \mathbb{R}^{N_1,L} $, $\bt{S}_2\in \mathbb{R}^{N_2,L}$ 
\item Reduced dimension for density: $n_1$ 
\item System matrices according to \eqref{eq:coord-fom-mat}: $\Massa$, $\Massb$, $\bt{J}$
\end{itemize}
OUTPUT: Reduction matrix: $\bt V \in \mathbb{R}^{N_1+N_2,n_1+n_2}$
\begin{enumerate}
\item \label{eq:alg:CPCA-s1} Calculate Cholesky factor $\bt L \in \mathbb{R}^{N_1,N_1}$ of $\Massa$ with $\Massa = \bt L^T \bt L$
\item \label{eq:alg:CPCA-s2} Determine $\tilde{\bt{V}}_1 \in \mathbb{R}^{N_1,n_1}$ consisting of the first $n_1$ left-singular vectors of $\bt K=\bt{L } [ \bt{S}_1 , \Massa^{-1} \bt{J} \bt{S}_2]$
\item \label{eq:alg:CPCA-s3} Determine matrices $\bt V_i \in \mathbb{R}^{N_i,n_i}$ for $i=1,2$\\
(whose column vectors are orthonormal in the scalar products induced by $\Massa$ and $\Massb$, respectively) with \\[0.1em]
	$\myim{\bt{V}_1} = \myim{\bt L^{-1} \tilde{\bt{V}}_1}$\\
	  $\myim{\bt{V}_2} = \myker{\bt{J}} \oplus \myim{\bt J^{\dag} (\bt Q \bt V_1)}$, \qquad	 where \, $\bt{J}^{\dag} = \Massb^{-1} \bt J^T (\bt J \Massb^{-1} \bt J^T)^{-1}$
\item Set $\bt V = \begin{bmatrix}
	\bt V_1 & \\
		  & \bt V_2
\end{bmatrix}$
\end{enumerate}
end
\end{lgrthm}
\end{samepage}
In Step \eqref{eq:alg:CPCA-s3} of Algorithm~\ref{alg:CPCA}, weighted scalar products are taken into account for better numerical stability.  Particularly, $\bt{J}^{\dag}$ is the right-inverse of $\bt{J}$ weighted w.r.t.\ the scalar product induced by  $\Massb$. Let us stress that the equality $\myim{\bt{V}_2} = \myker{\bt{J}} \oplus \myim{\bt J^{\dag} (\bt Q \bt V_1)}$ as well as $\myim{\bt{V}_1}$ and $\myim{\bt{V}_2}$ are independent of the choice of scalar product in the absence of rounding errors.

\begin{rmrk}
The computational complexity of Algorithm~\ref{alg:CPCA} is comparable to the standard {POD} approach when implemented appropriately. Since $\Massa$ is diagonal and $\Massb$ is tridiagonal, their Cholesky factorizations are of the same form. Step~\eqref{eq:alg:CPCA-s2} scales only linearly in the dimension of the {FOM} when the method of snapshots \cite{art:pod-kunisch,art:morKunV01} is employed. In Step~\eqref{eq:alg:CPCA-s3}, the kernel $\myker{\bt{J}}$ can be efficiently determined by a sparse LU-factorization \cite{art:gotsman-sparsenull}. Moreover, the right-inverse $\bt{J}^{\dag}$ does not have to be computed explicitly, but a few sparse linear equations can be solved instead, as can be seen using the Cholesky factorization of $\Massb$.
\end{rmrk}

\section{Quadrature-type complexity reduction} \label{sec:cr-approx} 

Although the {ROM}s are typically of much lower dimension than the {FOM}, i.e., $\mydim(\funSpace{V}_r) \ll N=\mydim(\funSpace{V}_f)$, they are in general not more efficient due to the lifting bottleneck related to the nonlinearities. The evaluation of the nonlinear integrals in System~\ref{sys:abstract} scales with the number of finite elements $K_{j}$, $j=1,\ldots ,J$, and hence is not independent of the FOM dimension, since $N=2J+\mydim(\funSpace{K})$, see, e.g., the friction term
\begin{align*}
	\Lscal r(\rho,m) m), \tfb \Rscal = \sum_{\substack{j =1,\ldots, J,\\
	\Onepipe \text{ given s.t.\ } K_{j} \subset \Onepipe  }} A^\Onepipe \int_{K_j} r(\rho,m) m \, \tfb \,  dx, &&(\rho,m)\in \funSpace{V}_r, \quad w\in \Vsb_r.
\end{align*}
We propose a complexity reduction by a quadrature-type approximation that restricts to integral evaluations at only few finite elements. It is trained towards the given snapshot data and, in contrast to conventional methods, regards compatibility conditions (Assumption~\ref{assum:quadrat-ansatz}). The resulting training problem is a semi-definite program with combinatorial aspects. We approach it algorithmically by a greedy procedure, as is generally done for complexity reduction.

\begin{rmrk}
The {CROM} distinguishes from the {ROM} only in the nonlinear terms. Its algebraic representation is similar to System~\ref{sys:rom-coord}, but $\bv{f}_r^\alpha$, $\bv{f}_r^\beta$, $\bv{f}_r^\gamma$ are substituted by complexity-reduced analogs. 
\end{rmrk}

\subsection{Training goal} \label{subsec:cr-optp}

The complexity reduction is done subsequently to the model order reduction and its training employs the {ROM} space $\funSpace{V}_r=\Vsa_r \times \Vsb_r$ explicitly. Similarly to \cite{art:comred-ecsw, art:hyperreduction-fem-cubature}, we consider a complexity reduction of quadrature-type. We aim for a complexity-reduced bilinear form 
\begin{align} \label{eq:quadrat-ansatz}
\Lscal b, \bar{b}\Rscal_c = \sum_{i\in I} \wei_i \int_{K_i} b\, \bar{b} \,dx, \qquad b, \bar{b} \in \funSpace{L}^2(\mathcal{E})
\end{align}
that approximates $\Lscal \cdot , \cdot \Rscal$, satisfies Assumption~\ref{assum:quadrat-ansatz} and realizes $\Lscal \cdot, \cdot\Rscal_\scps$ in System~\ref{sys:abstract}. In the training phase, we search for the index-set $I\subset \{1,\ldots, J\}$ and weights $\wei_i\geq0$ for $i\in I$. In particular, we aim for $|I|=n_c$ with $n_c \ll J$ small to ensure a reduced computational complexity for the {CROM}.

 To formulate the training goal, we assume an appropriate collection of snapshots of all nonlinear integrands in the {ROM} to be given, i.e.,
\begin{align} \label{eq:snapfuns}
	\ubv{S}_f = [f^1, \ldots, f^{\bar{L}} ], \hspace{1cm} \text{with } f^\ell \in \funSpace{L}^1(\mathcal{E}) \quad \text{for } \ell = 1,\ldots \bar{L}.
\end{align}
The complexity reduction is desired to be of high fidelity for the snapshot data, i.e., $\Lscal f^\ell, 1 \Rscal_c\approx \Lscal f^\ell, 1 \Rscal$ for $\ell=1,\ldots,\bar{L}$. Expressing that in a least-squares sense yields the following optimization problem for the training of the quadrature-rule.
\begin{prblm} \label{prob:cr-start}
Let nonlinear snapshot functions $\ubv{S}_f $ be given as in \eqref{eq:snapfuns}. Let $n_c \in \{1,\ldots, J \}$ and $\tilde{C} > 1$. Find an index-set $I^*\subset \{1, \ldots, J \}$ and weights $\wei^*_i \in \mathbb{R}$ for $i\in I^*$ as solution to the problem
	\begin{align*}
		\min_{\substack{I \subset \{1,\ldots, J \}, \,|I| =  n_c\\
			            \wei_i \geq 0, \,i \in I } }
			             & \quad \sum_{\ell=1}^{\bar{L}} \left( \sum_{i \in I} \wei_i \int_{K_i} f^\ell dx - \Lscal f^{\ell},1 \Rscal \right)^2 \\
		\text{s.t.} & \quad \frac{1}{\tilde{C}^2} \Lscal b,b \Rscal \leq  \sum_{i \in I} \wei_i \int_{K_i} b^2 dx \leq   {\tilde{C}^2} \Lscal b,b \Rscal  \,\, \quad \text{for } b \in \Vsa_r \cup \Vsb_r .
	\end{align*}
\end{prblm}
The constraints reflect the compatibility condition of Assumption~\ref{assum:quadrat-ansatz} that ensures structure preservation. In particular, the non-negativity of $\wei_i$ directly implies the non-negativity of the function $f \mapsto \Lscal f , 1 \Rscal_c$. 

As for the snapshots ${\ubv{S}}_f$, we consider the following mappings for $\ubv a = (\rho,m)$, related to the three nonlinear integral expressions occurring in the ROM (System~\ref{sys:abstract}),
\begin{align} \label{eq:Sf}
	\ubv a \mapsto  \left( P'(\rho)+ \frac{m^2}{2 \rho^2} \right) \partial_x \tfb^j_r , \qquad 
	\ubv a \mapsto r(\rho,m) m \, \tfb^j_r  , \qquad  
	\ubv a \mapsto 	v(\ubv a) \tfb^j_r  \qquad \text{ with  } v(\ubv a)= \frac{m}{\rho},
\end{align}
with $\tfb_r^j$ for $ j = 1 , \ldots , n_2$ denoting the basis for $\Vsb_r$. For the first two mappings, nonlinear snapshots are generated using the data set of the state snapshots $\ubv{S}=\{\ubv{a}^1,...,\ubv{a}^L\}$, $\ubv{a}^\ell\in \funSpace{V}_f$, i.e., the snapshots the ROM is trained for.
The third mapping appears in a time-derivative of the {ROM}. To treat this accordingly, we additionally assume that the state snapshots $\ubv{S}$ consist of one or several solution trajectories in time.  Particularly, given $\ubv{a}^l$ and $\ubv{a}^\ell$ are snapshots of the same trajectory at different time points $t_{k-1}$ and $t_{k}$, we collect the terms $(v(\ubv{a}^{\ell}) - v(\ubv{a}^l))/(t_{k}-t_{k-1})\tfb^j_r$ for $j = 1 , \ldots , n_2$, which relates to using a finite difference approximation on $\partial_t v(\ubv a)$ at time $t_k$. In case that ${\ubv{S}}$ represents one trajectory, we can approximate the time-derivative at $L-1$ instances, and the collection of all nonlinear snapshots ${\ubv{S}}_f$ has $\bar{L} = (3 L-1)n_2$ entries in total.

\begin{rmrk}\label{rem:cr-store-A}
By construction, the number of integrand-snapshots in ${\ubv{S}}_f$ is significantly larger than the number of state snapshots ${\ubv{S}}$ in \eqref{eq:snap-statespace}. When the problem at hand is of very large scale, a pre-processing of ${\ubv{S}}_f$ or other strategies to speed up the training phase might become necessary. We do not discuss this issue here but refer to \cite{art:efficient-integration-cubature , art:hyperreduction-fem-cubature} for some attempts in this direction.
\end{rmrk}


\subsection{Greedy implementation of training} \label{subsec:cr-solve}

For algorithmic reasons, we rewrite Problem~\ref{prob:cr-start} in algebraic form. Let $\vwei = [\wei_1, \ldots, \wei_J]^T \in \mathbb{R}^J$ denote the extended vector of weights, which is composed of the weights $\wei_i$ for $i\in I$ used in the quadrature and the entries $\wei_j=0$ for $j\notin I$. With $\bv{f}^\ell \in \mathbb{R}^J$ defined by $[\bv{f}^\ell]_i = \int_{K_i} f^\ell \,dx$ for the nonlinear integrands $f^\ell \in \funSpace{L}^1(\mathcal{E})$, $\ell =1,...,\bar L$ in \eqref{eq:snapfuns}, we can rewrite the cost functional in terms of $\vwei$ and
\begin{align}
\bt{A} = 
\begin{bmatrix}
	(\bv{f}^1)^T \\
	\vdots \\[0.2em]
	(\bv{f}^{\bar{L}})^T 
\end{bmatrix} \in \mathbb{R}^{\bar{L},J}, 
\qquad \qquad \bv{b} = \begin{bmatrix}
\Lscal f^1,1\Rscal \\
\vdots \\[0.2em]
\Lscal f^{\bar{L}},1\Rscal
\end{bmatrix}
 \in \mathbb{R}^{\bar{L}}.  \label{eq:prob:cr-start-matrix}
\end{align}
The constraint in Problem~\ref{prob:cr-start} enforces a norm equivalence between the $\funSpace{L}^2$-norm and its complexity-reduced counterpart $|| \cdot ||_c$. And as the {ROM} bases $\{\tfa_r^1,\ldots, \tfa_r^{n_1}\}$ and $\{\tfb_r^1,\ldots, \tfb_r^{n_2}\}$ are orthonormal in the $\funSpace{L}^2$-norm according to Algorithm~\ref{alg:CPCA}, the related mass matrix becomes the unit matrix. Its complexity-reduced counterpart reads
\begin{align*}
	\MassBoth_c(\vwei) &= \begin{bmatrix}
		\Massa_{c}(\vwei)  & \\
				& \Massb_c(\vwei) 
	\end{bmatrix}\\
	 \Massa_c(\vwei) &=  \sum_{i\in I} \wei_i \left[ \int_{K_i} \tfa_r^k \tfa_r^\ell  dx \right]_{\ell,k=1,\ldots, n_1}, \hspace{0.5cm}	\Massb_{c}(\vwei) =  \sum_{i\in I} \wei_i \left[ \int_{K_i} \tfb_r^k \tfb_r^\ell  dx \right]_{\ell,k=1,\ldots, n_2},
\end{align*}
from which it follows that $|| \cdot ||_c$ is a norm on the {ROM} space, if and only if $\MassBoth_c$ is positive definite. The more specific constraint in Problem~\ref{prob:cr-start} can be expressed as a bound on the spectrum $\sigma(\MassBoth_c)$ of $\MassBoth_c$, as can be shown using the estimate $\|\MassBoth_c^{-1} \|_2^{1/2} \| \bv a_r\| \leq  ( \bv a_r^T \MassBoth_c \bv a_r)^{1/2} \leq \|\MassBoth_c \|_2^{1/2} \| \bv a_r\|$, with $\|\MassBoth_c \|_2$ denoting the spectral norm, cf. \cite{art:morKunV01}.
In sum, Problem~\ref{prob:cr-start} can be equivalently written as the constrained least squares problem
 \begin{subequations} \label{eq:prob:cr-start}
	\begin{align}
		\min_{ \substack{I\subset \{1,\ldots J \}, \, |I| = n_c \\  \vwei=[\wei_1,\ldots,\wei_J]^T } } & \quad \| \bt{A} \vwei - \bv{b} \|^2  \label{eq:prob:cr-start-a} \\
		\text{s.t.} & \quad \wei_i \geq 0 \, \text{ for } i \in I, \quad \wei_j  = 0 \text{ for } j \notin I\\ 
					& \quad \sigma\left(\MassBoth_{c}(\vwei)\right) \subset \left[ \tilde{C}^{-2}, \tilde{C}^2 \right]. \label{eq:prob:cr-start-c}
	\end{align}
\end{subequations}

The constraint~\eqref{eq:prob:cr-start-c} makes this a semi-definite program \cite{book:Jarre2000}. And inherently in the complexity reduction, there is a combinatorial aspect, in our case the selection of the index set $I$ of active finite elements. Particularly the latter makes this problem too hard to be solved to global optimality. Thus, we rely on a greedy approximation, as most complexity-reduction methods do, cf., \cite{art:deim-state-space-err,art:empint-maday04,art:moramHamDissHesthaven}, see also \cite{art:comred-ecsw, art:efficient-integration-cubature} for quadrature-type approaches similar to ours.
The basic idea is to alternate between enlarging the index set $I$ by a greedy search and constructing an optimal weight for that fixed index set. Algorithm~\ref{alg:p1b:empQuad} summarizes the procedure. 
While the non-negativity constraint for the weights $\wei_i$ is included explicitly, the eigenvalue constraint \eqref{eq:prob:cr-start-c} is only included as a safeguard at the end of the algorithm. As will be shown in the numerical results in Section~\ref{sec:numres}, the latter condition is obtained by our greedy search for reasonable choices of $n_c$ without further doing for our problem at hand. 

\begin{samepage}
\begin{lgrthm}[Greedy empirical quadrature weights]\text{$ $} \label{alg:p1b:empQuad}\\
INPUT:
\begin{itemize}
\item Data matrix and vector as in \eqref{eq:prob:cr-start-matrix}: $\bt{A}$, $\bv{b}$
\item Number of active finite elements: $n_c$
\item Constant: $\tilde{C} > 1$
\item Complexity-reduced mass matrix as function of weights: $\MassBoth_c(\cdot)$
\end{itemize}
OUTPUT: Vector of quadrature-weights: $\vwei^{nc}$
\begin{enumerate}
\item Define function $F(\vwei) = || \bt{A} \vwei - \bv{b} ||^2$ with gradient $\nabla F(\vwei) = 2\bt{A}^T(\bt{A} \vwei- \bv b)$
\item Initialize $I_0 = \{ \}$ and $\vwei^0 = \bv{0} \in \mathbb{R}^{J}$ 
	\item for $k = 1:n_c$
		\begin{enumerate}[a)]
				\item \label{eq:alg:p1b:empQuad-S0} Define set of candidates $I_c = \{1, \ldots, J\} \setminus I_{k-1}$ 
				\item Find $j_{max} = \text{argmax}_{j\in I_c} -\left[ \nabla_{{}} F(\vwei^{k-1}) \right]_j $
				\item Set $I_k = I_{k-1} \cup \{ j_{max} \}$
				\label{eq:alg:p1b:empQuad-S2} 
				\item \label{eq:alg:p1b:empQuad-S3} Find $\vwei^k$ as solution to
				\begin{align*}
							\min_{ \vwei=[\wei_1,\ldots,\wei_J]^T }  & \,\, F(\vwei) \\
		\text{s.t.} & \,\, \wei_i \geq 0 \, \text{ for } i \in I_k, \quad \wei_j  = 0 \text{ for } j \notin I_k 
				\end{align*}	
		\end{enumerate}
	\item[] endfor
		\item if $\sigma\left(\MassBoth_{c}(\vwei)\right) \not\subset\left[ \tilde{C}^{-2}, \tilde{C}^2 \right] $\\
			\text{ } \hspace{0.5em} Set $n_c \hookleftarrow n_c+1$ and go to \eqref{eq:alg:p1b:empQuad-S0}\\
			endif			
\end{enumerate}
end
\end{lgrthm}
\end{samepage}
Let us note that only $n_c$ columns or less of the data matrix $\bt A$ are used in most steps of the algorithm. The exception is Step \eqref{eq:alg:p1b:empQuad-S2}, where the gradient is evaluated once. It can be avoided to store the full matrix, if storage requirements are critical, cf.,  \cite{art:orthogonalgreedy-nquyen19,inproc:gradientpursuit}  and Remark~\ref{rem:cr-store-A}. Most importantly, the non-negative least-squares problem of Step \eqref{eq:alg:p1b:empQuad-S3} needs only to be solved over very low-dimensional linear spaces. We do this using the built-in \texttt{Matlab} routine \texttt{lsqnonneg} in our implementation, which realizes an active-set method.

\section{Numerical results} \label{sec:numres}
For the numerical studies of this section we consider our model problem in parameter regimes relevant for gas distribution networks. The model equations are given by the isothermal Euler equations in a friction-dominated regime. The network parameters are extracted from the \texttt{gaslib} \cite{art:gaslib-2017} with minor adaptions. Networks of pipes are considered, where each edge $\Onepipe\in \mathcal{E}$ represents a pipe of a cross-sectional area $A^{\Onepipe}= \pi/4\, (D^\Onepipe)^2$. The $A^{\Onepipe}$ act as edge weights and are prescribed by the diameter $D^\Onepipe$. The reference values for density and mass flow are taken as $\rho_\star=1~[\mathrm{kg\,m^{-3}}]$ and $(Am)_\star= 1~[\mathrm{kg\,s^{-1}}]$. An isothermal pressure law is used for $p$ $[\mathrm{Pa}]$, $p(\rho) = RT {\rho}/(1-RT\alpha \rho)$, with $T= 283~[\mathrm{K}]$,  $R = 518~[\mathrm{J} (\mathrm{kg \, K})^{-1}]$ and $\alpha = -3 \cdot 10^{-8}~[\mathrm{Pa}^{-1}]$,
which implies the pressure potential $P(\rho) = RT \rho \log{\left( {(1-RT\alpha \rho)\rho_\star}/{\rho} \right)}$. 

General boundary conditions are treated according to \cite{art:lilsailer-nlfow}. Specifically, for each mass flow boundary condition a Lagrange multiplier is added and also kept in the reduced models. The {FOM} is realized according to Section~\ref{subsec:fom-training-prob} using a uniform mesh on each pipe with the maximal spatial step size $\Delta_x \leq 200~[\mathrm{m}]$. Time discretization is carried out by an implicit Euler method in the primitive variables $\rho$ and $v$, cf. \cite{art:lilsailer-nlfow}, with a constant time step $\Delta_t = 1~[s]$, i.e., 3600 steps per hour of simulation time. The resulting nonlinear systems in each time step are solved by the Newton's method. If this fixed-point iteration diverges for any time step, we consider this as a simulation-breakdown. Given the {FOM} solution $\ubv{a}$, the error of an approximation $\tilde{\ubv{a}}$ (obtained by a reduced model or orthogonal projection) is measured by the $\funSpace{L}^2$-norm in space and the supremum-norm in time, yielding the relative error $E_T = {( \max_{t\in [0,T]}  ||\ubv{a}(t)- \tilde{\ubv{a}}(t)||)} / {\max_{t\in [0,T]}  ||\ubv{a}(t)||}$.
All numerical results have been generated using \texttt{MATLAB} Version 9.1.0 (R2016b) on an Intel Core~i5-7500 CPU with 16.0GB RAM.

In Section~\ref{subsec:numDiamNet} some qualitative properties of our complexity reduction approach are showcased at a smaller academic network scenario (diamond network).  More quantitative studies with comparisons to other well-established, but non-structure-preserving, model reduction methods are presented in Section~\ref{subsec:numLargeNet} using a network of realistic size.

\subsection{Qualitative study for complexity reduction}\label{subsec:numDiamNet}

\begin{figure}[tb]
\begin{center}
  \begin{tikzpicture}[scale=2.3]
  \node[circle,draw,inner sep=2pt] (v1) at (-1.87,0) {$\blsVertex_1$};
  \node[circle,draw,inner sep=2pt] (v2) at (-0.87,0) {$\blsVertex_3$};
  \node[circle,draw,inner sep=2pt] (v3) at (0,0.5) {$\blsVertex_4$};
  \node[circle,draw,inner sep=2pt] (v4) at (0,-0.5) {$\blsVertex_5$};
  \node[circle,draw,inner sep=2pt] (v5) at (0.87,0) {$\blsVertex_6$};
  \node[circle,draw,inner sep=2pt] (v6) at (1.87,0) {$\blsVertex_2$};
  \draw[->,thick,magenta,line width=2.0pt] (v1) -- node[above,sloped] {$\Onepipe_1$} ++ (v2);
  \draw[->,thick,magenta,line width=1.5pt] (v2) -- node[above,sloped] {$\Onepipe_2$} ++ (v3);
  \draw[->,thick,line width=2.0pt] (v2) -- node[above,sloped] {$\Onepipe_5$} ++ (v4);
  \draw[->,thick,line width=1.5pt] (v3) -- node[above,sloped] {$\Onepipe_7$} ++ (v4);
  \draw[->,thick,magenta,line width=1.5pt] (v3) -- node[above,sloped] {$\Onepipe_3$} ++ (v5);
  \draw[->,thick,line width=2.0pt] (v4) -- node[above,sloped] {$\Onepipe_6$} ++ (v5);
  \draw[->,thick,magenta,line width=1.5pt] (v5) -- node[above,sloped] {$\Onepipe_4$} ++ (v6);
  \end{tikzpicture}
\end{center}  
\caption{Diamond network with lengths [km] $\{l^{\Onepipe_i} :\Onepipe_i\in\mathcal{E} \} =  \{40  ,   38    ,18    ,15    ,28    ,27    ,25\}$ and diameters [m] $\{D^{\Onepipe_i}:  \Onepipe_i\in \mathcal{E} \} = \{1.3  ,  1, 1 , 1 , 1.3, 1.3 ,1\}$. The domain $ \{ \Onepipe_i \in \mathcal{E}: \, i =1, \ldots, 4  \}$ is marked in magenta.\label{fig:diamtop}}
\end{figure}
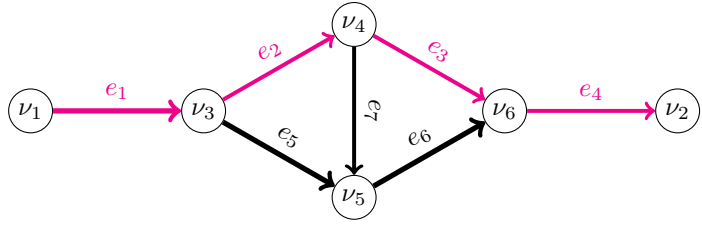

The topology and parameters for the academic network example (diamond network) are given in Fig.~\ref{fig:diamtop}. For the friction factor, we consider two different choices $\lambda= 0.01$ (high friction) and $\lambda= 0.002$ (lower friction). At the two boundary nodes $\nu_1$, $\nu_2$, we prescribe the boundary conditions according to
\begin{align*}
	\rho(t,\nu_1) &= (60 + { u }(t))\rho_*, \hspace{1.2cm}  A  m(t,\nu_2) = 200 (Am)_*, \hspace{2cm} t \in [0,2\,t_\star]
\end{align*}
with reference time $t_* = 1 [h]$ and the time-varying input profile\\
\begin{tabular}{lrl}
{\hspace{0.6cm}
\begin{minipage}{0.43\textwidth}
 $	{ u }(t) = \begin{cases} 5( \,{4t}/{t_\star}), & \hspace*{0.5cm} 0 \le t<t_\star/4 \\ 5\, (2-{4t}/t_\star), & \hspace*{0.4cm} t_\star/4 \le t<3/8 t_\star \\
 2.5, & \hspace*{0.5cm}  3/8 t_\star \le t. \end{cases} $
 \end{minipage}
} 
\vspace{1.5cm}
&
\hspace*{1.5cm}
\begin{minipage}{0.022\textwidth}
{\hspace{0.6cm}
{\normalsize \rotatebox{90}{Profile ${ u }$}\\ \vspace{0.0cm}
}}
\end{minipage}
&
{\hspace{-1.5cm}
\begin{minipage}{0.43\textwidth}
\center
\includegraphics[width=0.7\textwidth]{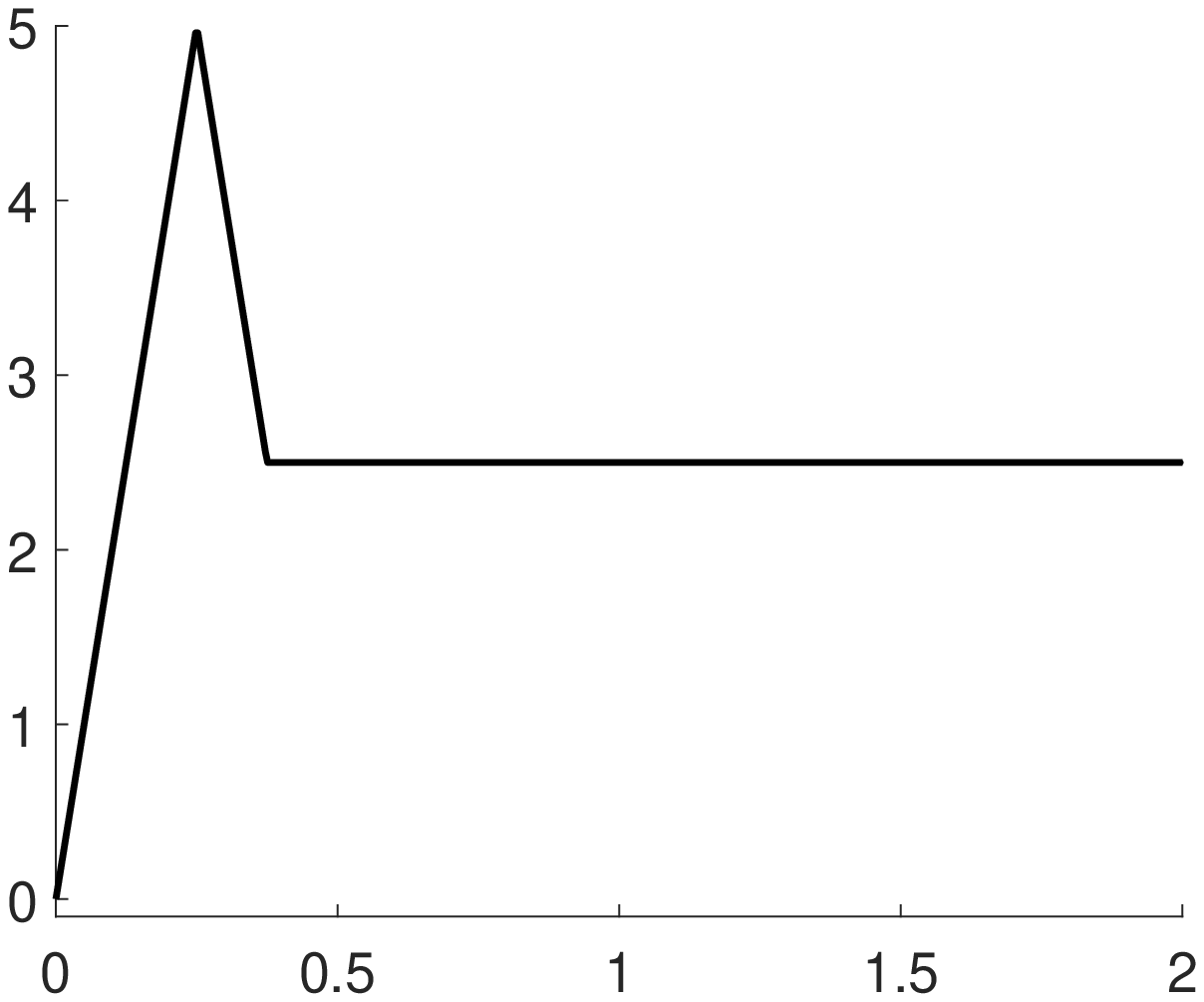}\\
 \hspace{0.5cm} {\normalsize time $t$ [h]}
\end{minipage}
}
\end{tabular}\\[-3.5em]
The initial values are chosen as the stationary solution related to the boundary conditions at $t=0$. The {FOM} yields a system of dimension $N=1922$, and the reduction methods (model order- and complexity-reduction) are trained towards 1000 equally distributed time snapshots of the solution trajectory of the {FOM}. The dimension of the {ROM} space $\funSpace{V}_r$ is fixed to $n=14$. As regards the complexity reduction, we vary the parameter $n_c$, as the qualitative behavior of this step is is the focus in the upcoming.

\begin{figure}[tb]
\begin{tabular}{rl|l}
\begin{minipage}{0.022\textwidth}
{\hspace{-0.3cm}
{\normalsize\rotatebox{90}{Density $\rho$}  \vspace{1cm}  \, 
\rotatebox{90}{Mass flow $Am$} 
}}
\end{minipage}
&
{\hspace{-0.4cm}
\begin{minipage}{0.45\textwidth}
\center
\hspace{0.5cm} {\large{\underline{Case: $\lambda = 0.01$ }}} \\
{\large{\underline{Solution}}} \hspace{1.6cm}  {\large{\underline{Difference to {FOM}}}}\\
\includegraphics[height = 0.4\textwidth, width = 0.44\textwidth]{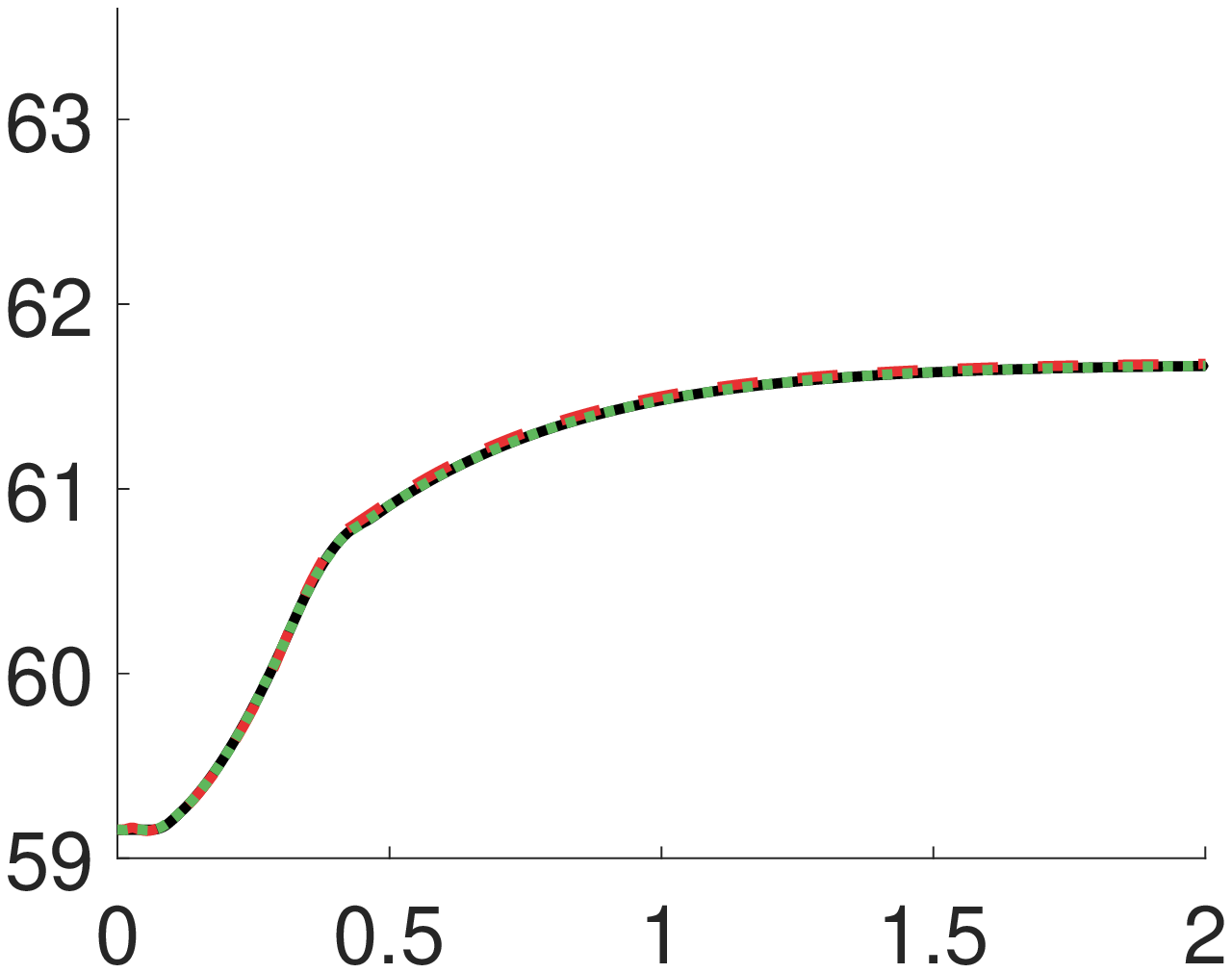}
 \includegraphics[height = 0.4\textwidth, width = 0.44\textwidth]{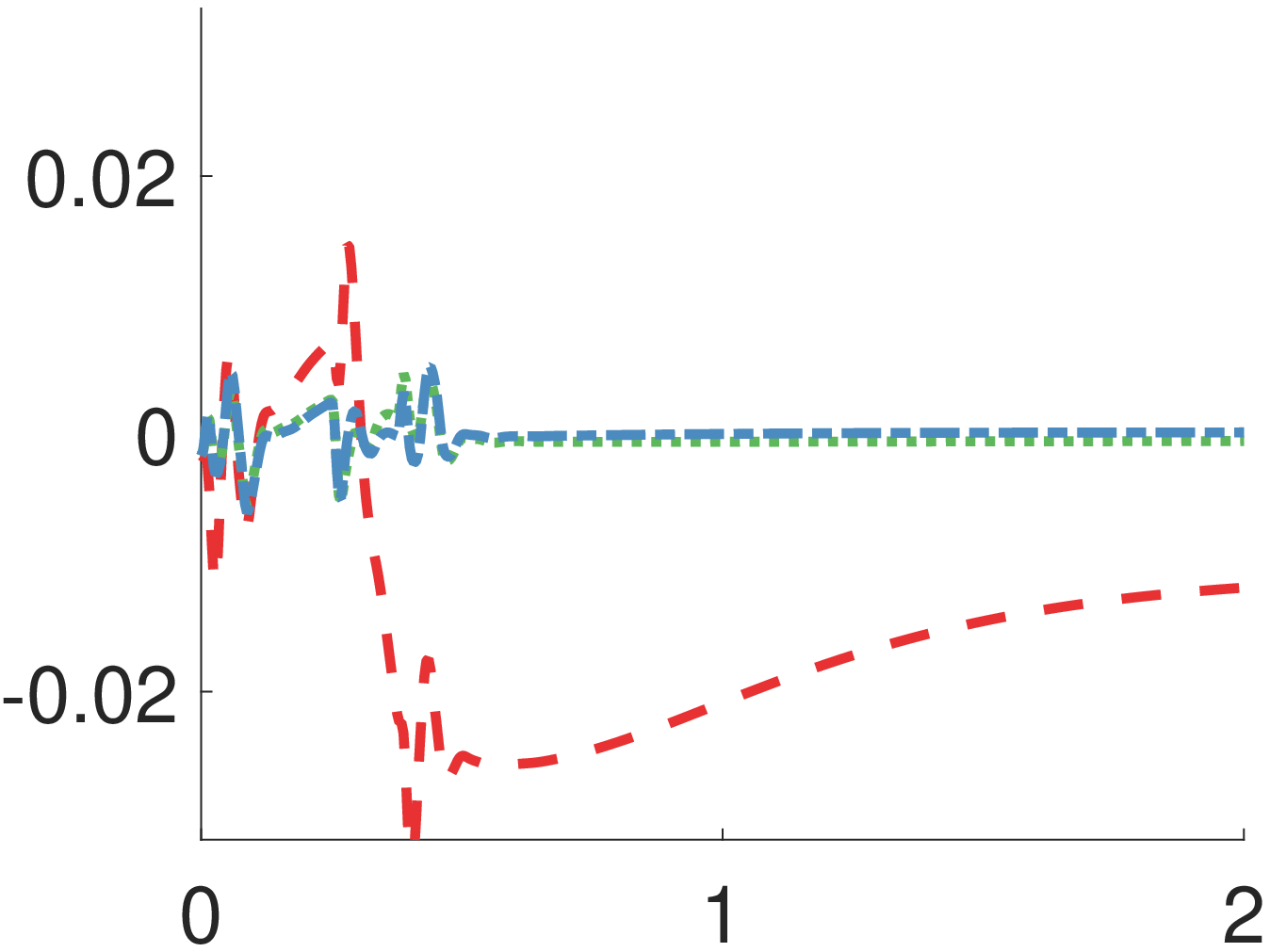}\\
\hspace{-0.1cm} 
\includegraphics[height = 0.4\textwidth, width = 0.44\textwidth]{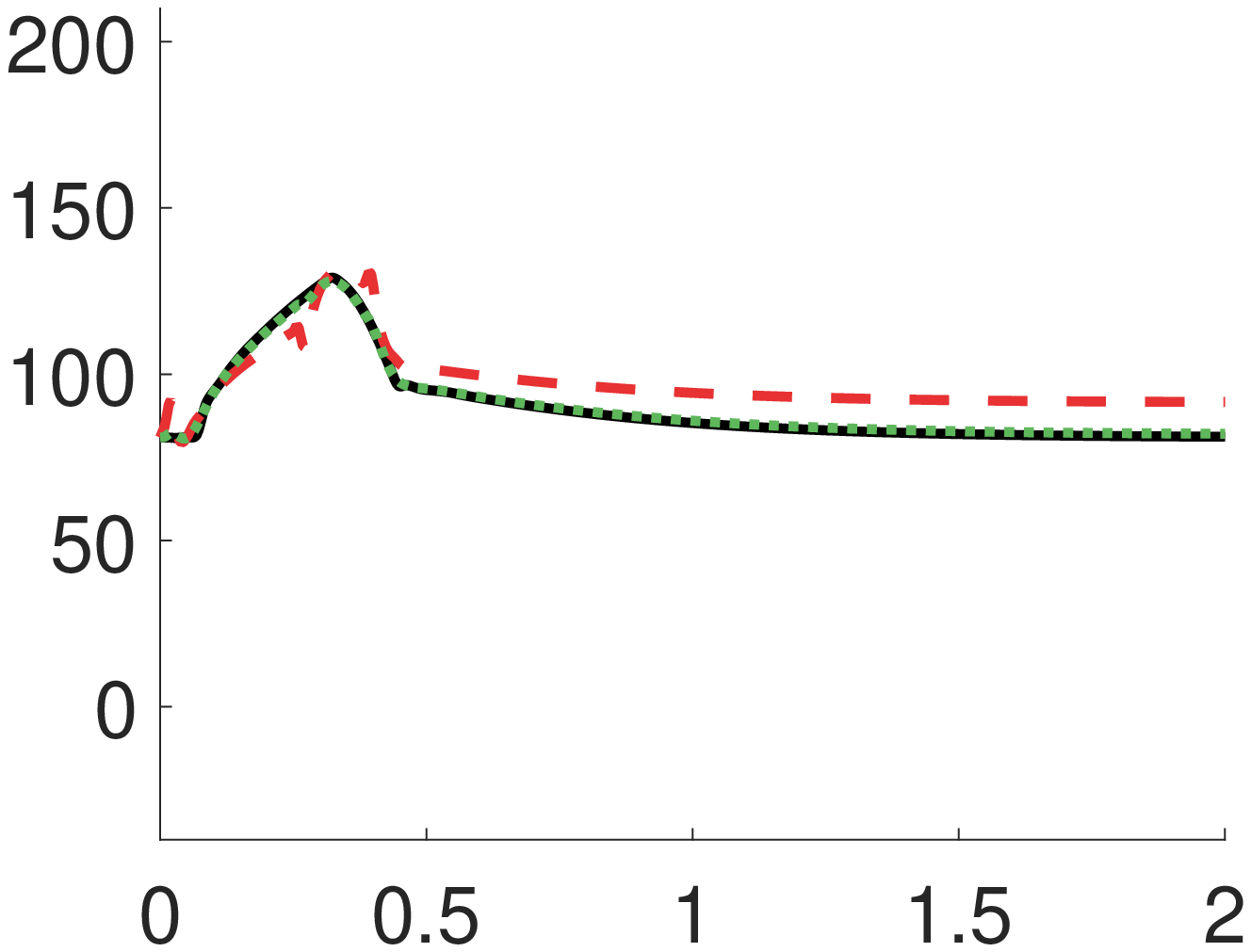} 
\includegraphics[height = 0.4\textwidth, width = 0.44\textwidth]{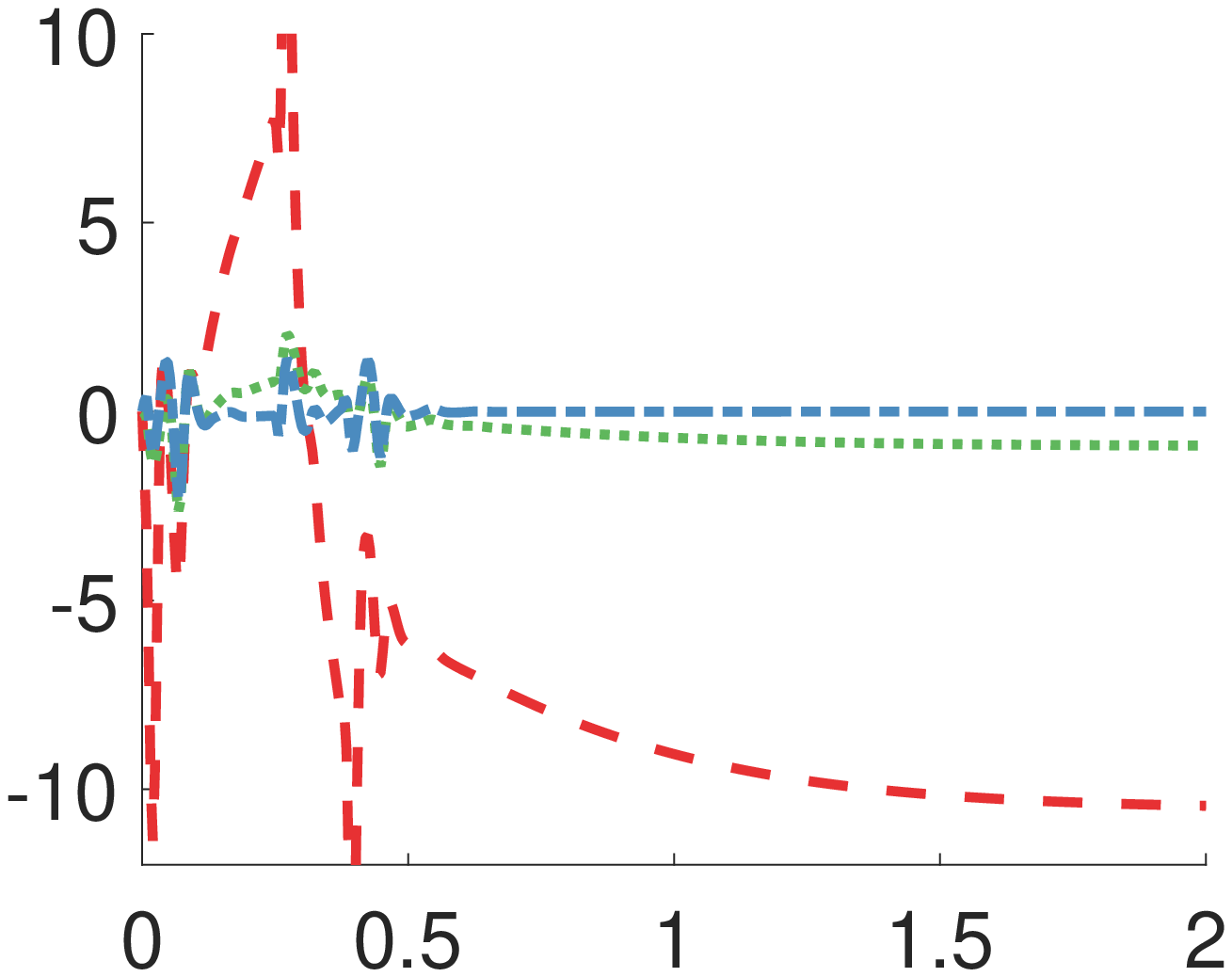}\\
 {\normalsize Time $t [h]$}  \hspace{1.0cm}  {\normalsize Time $t [h]$} 
\end{minipage}
}
&
\begin{minipage}{0.45\textwidth}
\center
\quad {\large{\underline{Case $\lambda = 0.002$}}}\\
{\large{\underline{Solution}}} \hspace{1.6cm}  {\large{\underline{Difference to {FOM}}}}\\
\includegraphics[height = 0.4\textwidth, width = 0.44\textwidth]{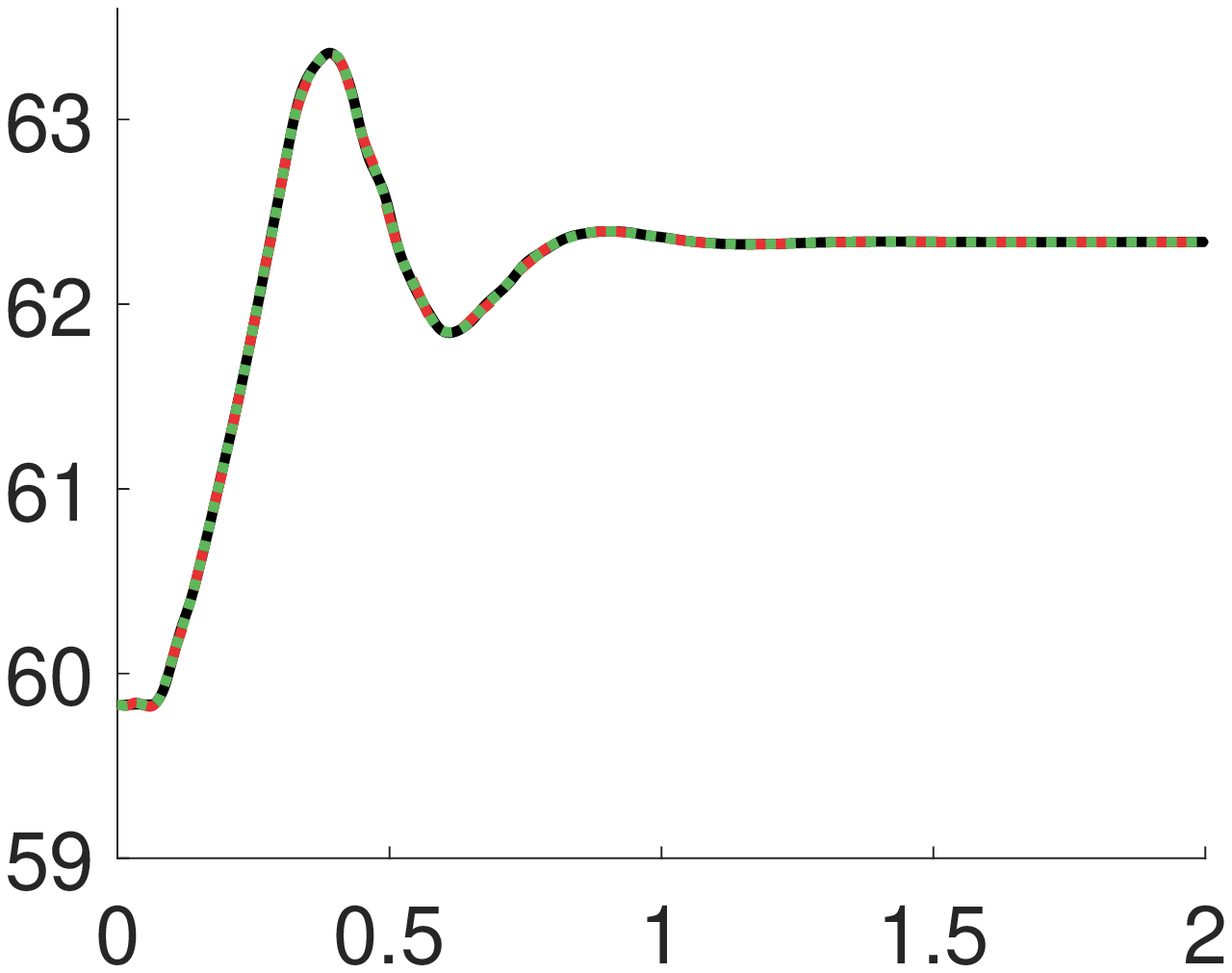}
 \includegraphics[height = 0.4\textwidth, width = 0.44\textwidth]{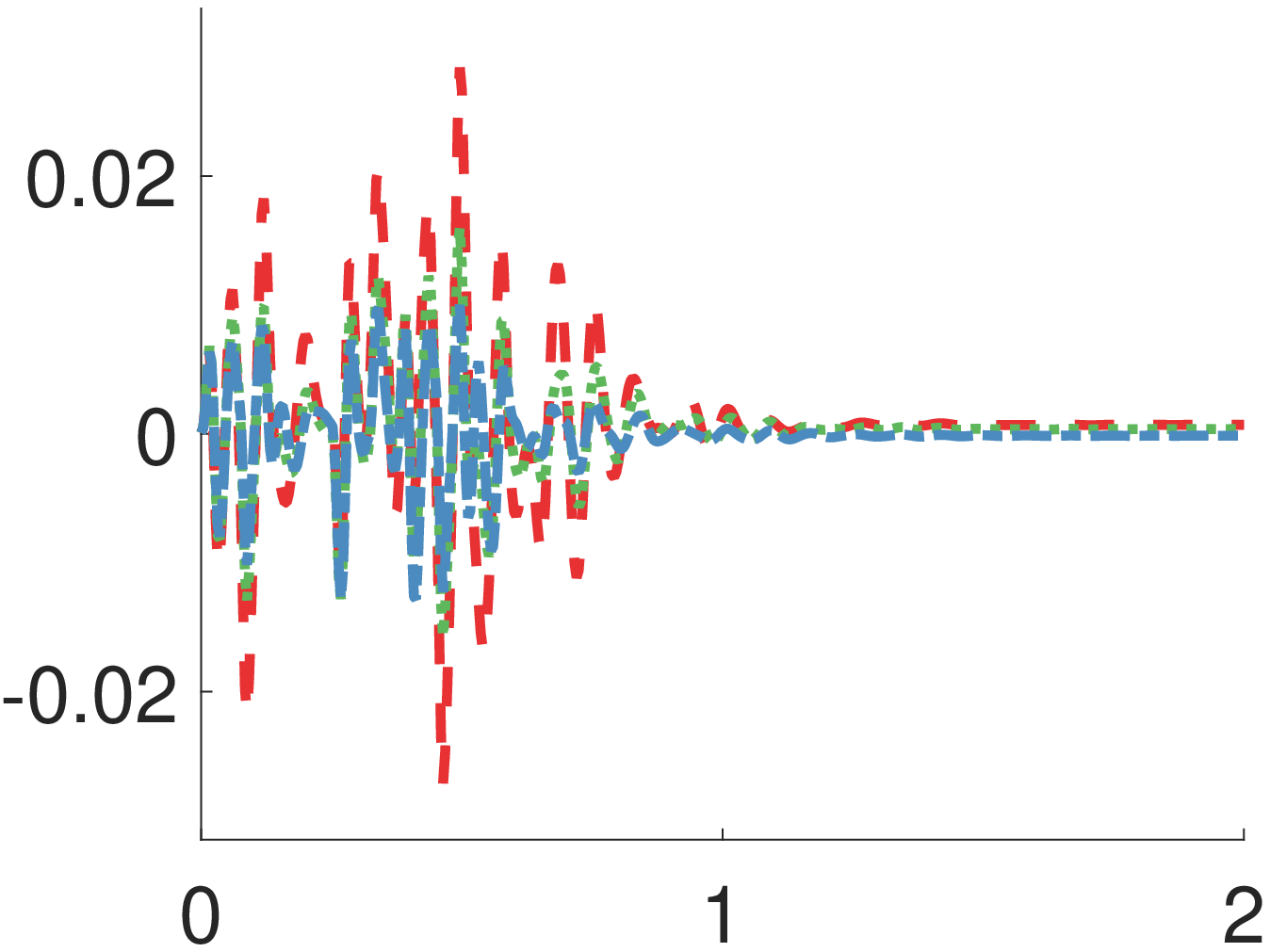}\\
 \hspace{0.1cm} 
\includegraphics[height = 0.4\textwidth, width = 0.44\textwidth]{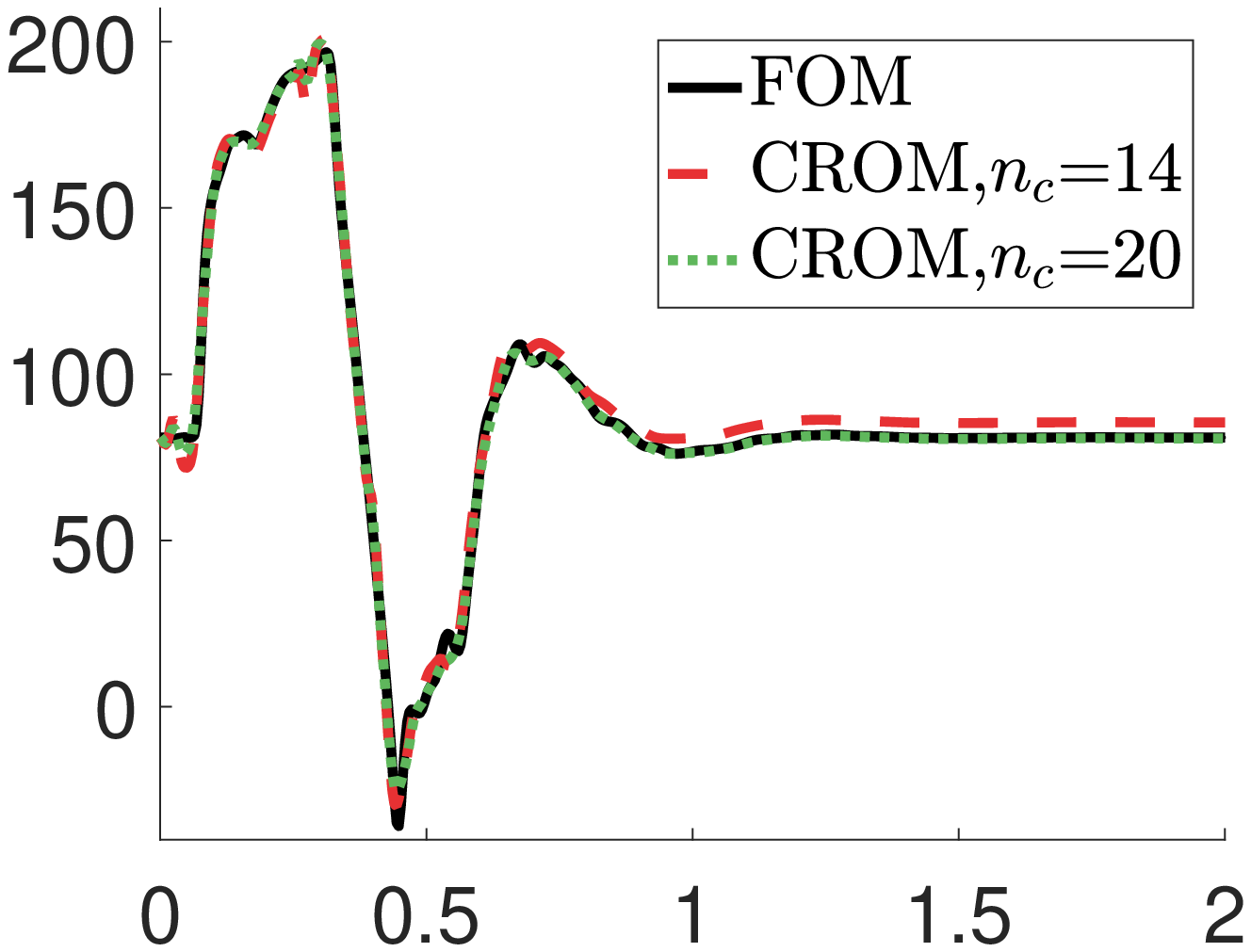} 
\hspace{0.1cm} 
\includegraphics[height = 0.4\textwidth, width = 0.44\textwidth]{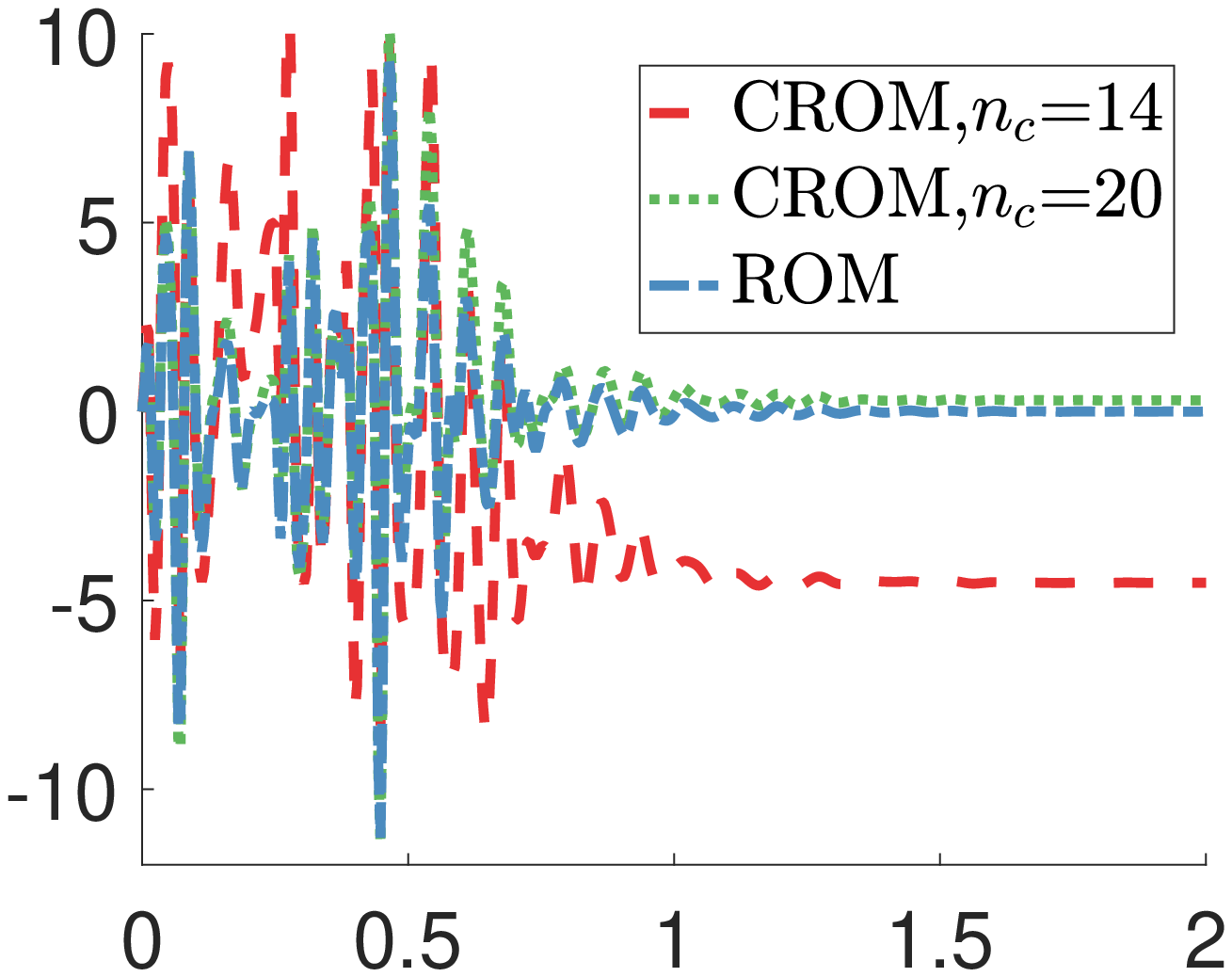}\\
 {\normalsize Time $t [h]$}  \hspace{1.0cm}  {\normalsize Time $t [h]$} 
\end{minipage}
\end{tabular}
\caption{Diamond network, solutions and differences to the {FOM} at end of pipe $\Onepipe_2$, plotted over time.\label{fig:diam-time-ev}}
\end{figure}

\begin{figure}[h]
\begin{tabular}{rll|lrl}

\begin{minipage}{0.00\textwidth}
\end{minipage}
&
\begin{minipage}{0.022\textwidth}
{
{\normalsize\rotatebox{90}{Quadrature weights}\\ 
}}
\end{minipage}
&
\begin{minipage}{0.43\textwidth}
\center
\hspace{0.1cm} \includegraphics[height = 0.4\textwidth, width = 0.9\textwidth]{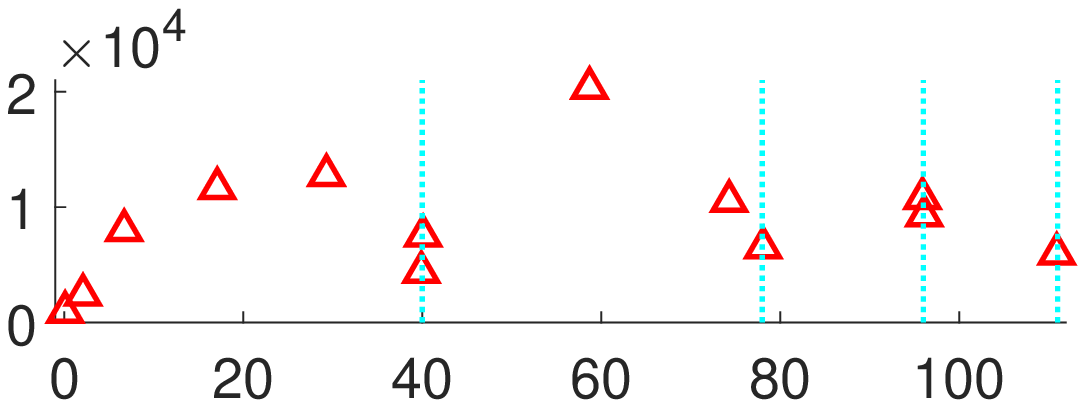}\\
\hspace{0.1cm} \includegraphics[height = 0.4\textwidth, width = 0.9\textwidth]{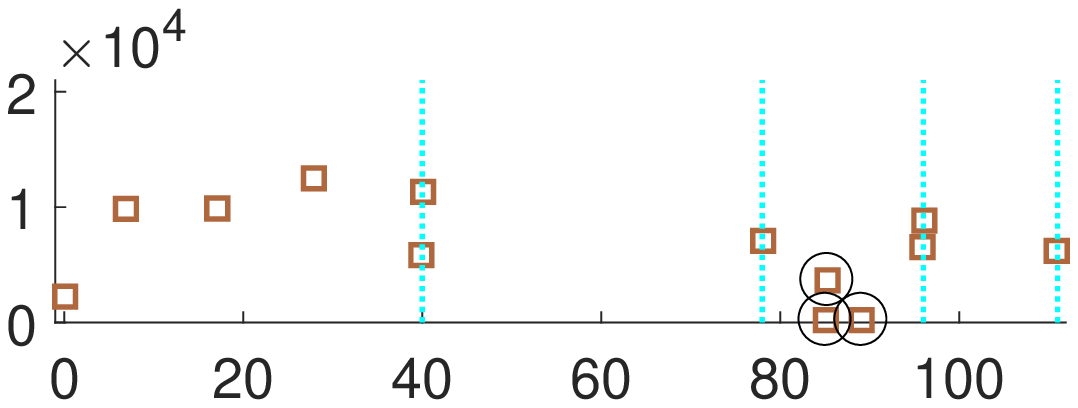}\\
 \hspace{0.5cm} {\normalsize Space $x [km]$}
\end{minipage}
&
&
\begin{minipage}{0.01\textwidth}
{
{\normalsize\rotatebox{90}{Relative error $E_T$}\\ 
}}
\end{minipage}
&
{\hspace{-0.6cm}
\begin{minipage}{0.43\textwidth}
\center
\hspace{0.1cm} \includegraphics[height = 0.77\textwidth, width = 0.8\textwidth]{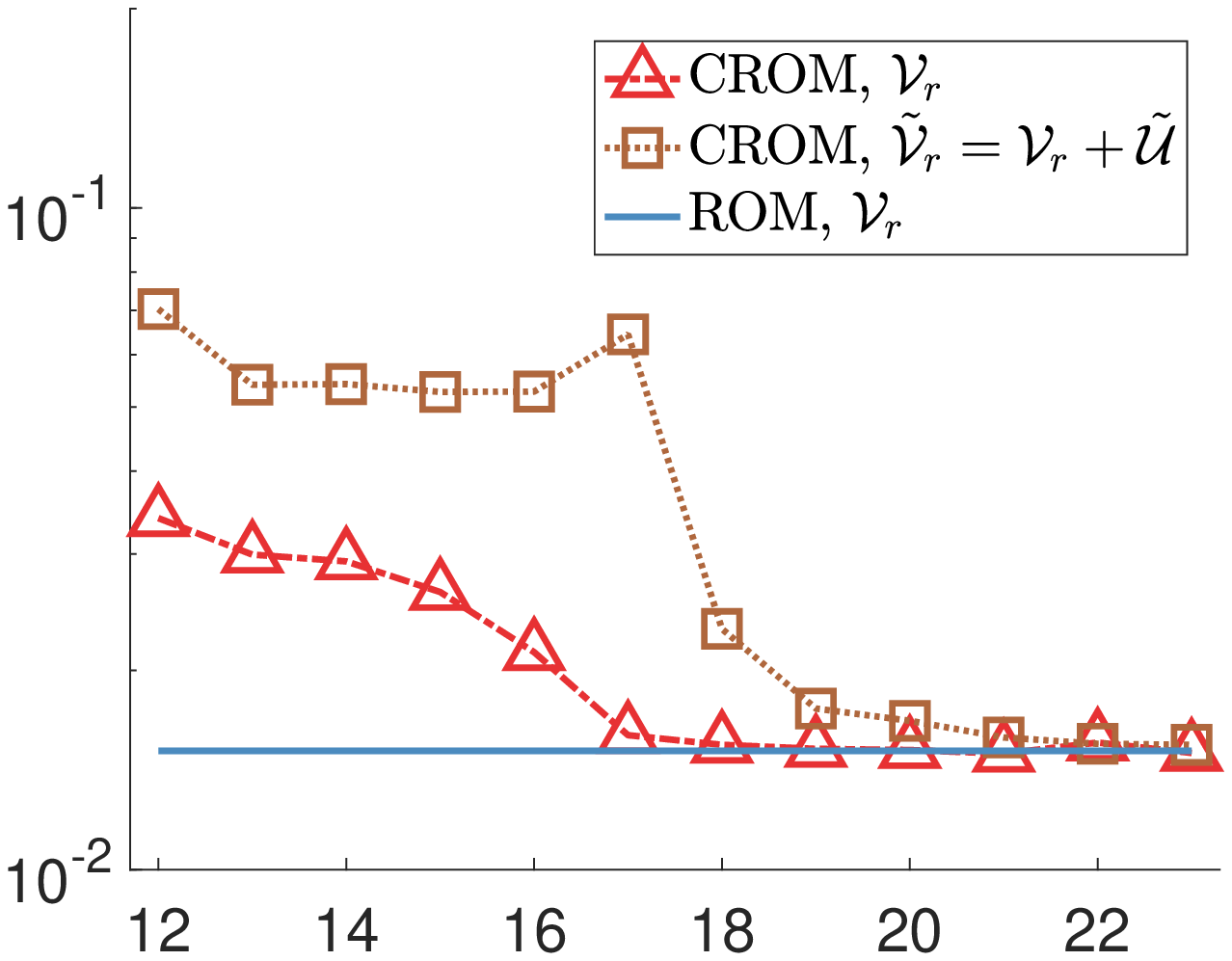}\\
 \hspace{0.5cm} {\normalsize Parameter $n_c$}
\end{minipage}
}
\end{tabular}
\caption{Diamond network, influence of {ROM} space on {CROM}. The space $\tilde{\funSpace{V}}_r = \funSpace{V}_r + \tilde{\funSpace{U}}$ includes $\mydim(\tilde{\funSpace{U}})$=4 artificial additional basis functions on edge $\Onepipe_3$. \textit{Left:} Locations and values of quadrature-weights along edges $\Onepipe_1$ to $\Onepipe_4$ for $n_c=20$, and $\funSpace{V}_r$ \textit{(top)} or  $\tilde{\funSpace{V}}_r$ \textit{(bottom)}. The encircled weights are the ones on the support of  $\tilde{\funSpace{U}}$, and the cyan vertical lines indicate junctions. \textit{Right:} Relative error $E_T$ for different choices of $n_c$.
\label{fig:diam-art}}
\end{figure}

For the two different cases of friction factor, the time evolution of the solution is visualized in Fig.~\ref{fig:diam-time-ev}. Hereby, the end of the edge $\Onepipe_2$ (at node $\nu_4$) is used as fixed spatial position. As to be expected, stronger damping effects are observed for larger $\lambda$, especially for the density. Generally speaking, models with higher damping effects are mostly better suited for model reduction, but the effect seems rather minor in the parameter range we encounter in the context of gas transportation.
For both choices of $\lambda$, the {CROM}s with $n_c=20$ produce very accurate time responses, see Fig.~\ref{fig:diam-time-ev}. But also when the complexity reduction is not well-resolved, as is the case for $n_c=14$, stable simulations are produced, which is in accordance to our energy bound (Theorem~\ref{thrm:energy-diss-cr}). The relative errors over the full simulation with $\lambda=0.002$ are shown Fig.~\ref{fig:diam-art}-\textit{right}. We observe that the error the complexity reduction step adds becomes negligible for $n_c\geq 18$. Moreover, the compatibility condition on the spectrum of the complexity reduced mass matrix, $\sigma(\MassBoth_{c}) \subset [\tilde{C}^{-2},\tilde{C}^2]$, is fulfilled for $n_c \geq 12$ and moderately large $\tilde{C}$, as is indicated by the condition number of $\MassBoth_{c}$, see Table~\ref{table:diam-cond}. 



%

\begin{table}[tb]
\renewcommand{\arraystretch}{1.2}
\small 
\begin{center}
\begin{tabular}{l || c|c|c|c|c|c|c|c|c|c |c  |c}
Parameter $n_c$ & 12 & 13   &   14   & 15  &  16  &  17  &  18  &  19  &  20 & 21  &  22  &  23 \\ 
\hline
\hline
CROM, $\funSpace{V}_r$  & 2.53    &    2.03   &  1.96  &   1.98  &  1.82  &  1.43  &  1.19  &  1.22  &  1.19  &  1.18  &  1.17		& 1.16 \\ 
\hline
CROM, $\tilde{\funSpace{V}}_r = \funSpace{V}_r + \tilde{\funSpace{U}}$ & 69.24  &  66.74  &  33.88  &  13.77  &  11.99  &  10.20  &   2.50  &   1.75  &   1.50  &   1.53  &   1.32    &   1.30\\
\end{tabular}
\medskip
\caption{
Condition number \texttt{cond}($\MassBoth_c$) w.r.t.~the spectral norm for the two {CROM} models for different choices of $n_c$. The measure is strongly related to the compatibility condition $\sigma(\MassBoth_{c}) \subset [\tilde{C}^{-2},\tilde{C}^2]$ from Section~\ref{subsec:cr-solve}.
\label{table:diam-cond}
}
\end{center}
\end{table}

An important qualitative property concerning the training of our complexity reduction method, which we want to showcase, is its adaption to the {ROM} space. To demonstrate that, we consider reduced models with an artificially enlarged {ROM} space given as $\tilde{\funSpace{V}}_r = \funSpace{V}_r + \tilde{\funSpace{U}}$, where $\tilde{\funSpace{U}}$ consists of functions with very local support. That is, two of the basis functions of the {FOM} for the density, which lie symmetrically around the midpoint of edge $\Onepipe_3$, are added to $\Vsa_r$, and the space $\Vsb_r$ is supplemented according to Assumption~\ref{assum:compatV1V2}-(A1). Thus $\tilde{\funSpace{V}}_r$ is compatible, and it holds $18 = \mydim(\tilde{\funSpace{V}}_r)> \mydim({\funSpace{V}}_r)=14$, but almost no improvement in terms of fidelity is observed for the scenario. In fact, the extension $\tilde{\funSpace{U}}$ causes the complexity reduction to produce less accurate models with larger \texttt{cond}($\MassBoth_c$), at least for small $n_c$, see Fig.~\ref{fig:diam-art}\textit{-right} and Table~\ref{table:diam-cond}. The reason becomes evident when looking at the quadrature points obtained for the {CROM} without and with $\tilde{\funSpace{U}}$. In Fig.~\ref{fig:diam-art}\textit{-left} the locations and values of the quadrature weights are visualized along the path of the edges $\Onepipe_1$ to $\Onepipe_4$ (marked in magenta in Fig.~\ref{fig:diamtop}). When the space $\tilde{\funSpace{U}}$ is added, we find three weights located on its support, see the encircled markers in the figure. Let us also mention that the weights do not cluster, the pairs that are observed around the junctions (the cyan vertical lines) lie on different pipes. Summarizing, the enrichment of the reduced space by $\tilde{\funSpace{U}}$ hampers the complexity reduction, but our method still produces stable simulations. This illustrates that our training procedure promotes Assumption~\ref{assum:quadrat-ansatz} and can cope with an enrichment of the reduction basis with elements chosen independently of the snapshot data. In more relevant scenarios, the space enrichment might be necessary due to compatibility conditions such as Assumption~\ref{assum:compatV1V2}, and the effects on the complexity reduction can be expected to be less exaggerated. Let us stress that other complexity reduction methods, such as the empirical interpolation methods \cite{art:deim-state-space-err,art:empint-maday04},
are trained independently of the choice of the {ROM} space and thus cannot adapt to it. Possibly related to that, we observe that the discrete interpolation method (DEIM) inherits stability issues when combined with our model order reduction approach, see the end of Section~\ref{subsec:numLargeNet}.

\subsection{Comparison to non-structure-preserving reduction methods}\label{subsec:numLargeNet}
\begin{figure}[tb]
\center
\includegraphics[width= 0.45\textwidth]{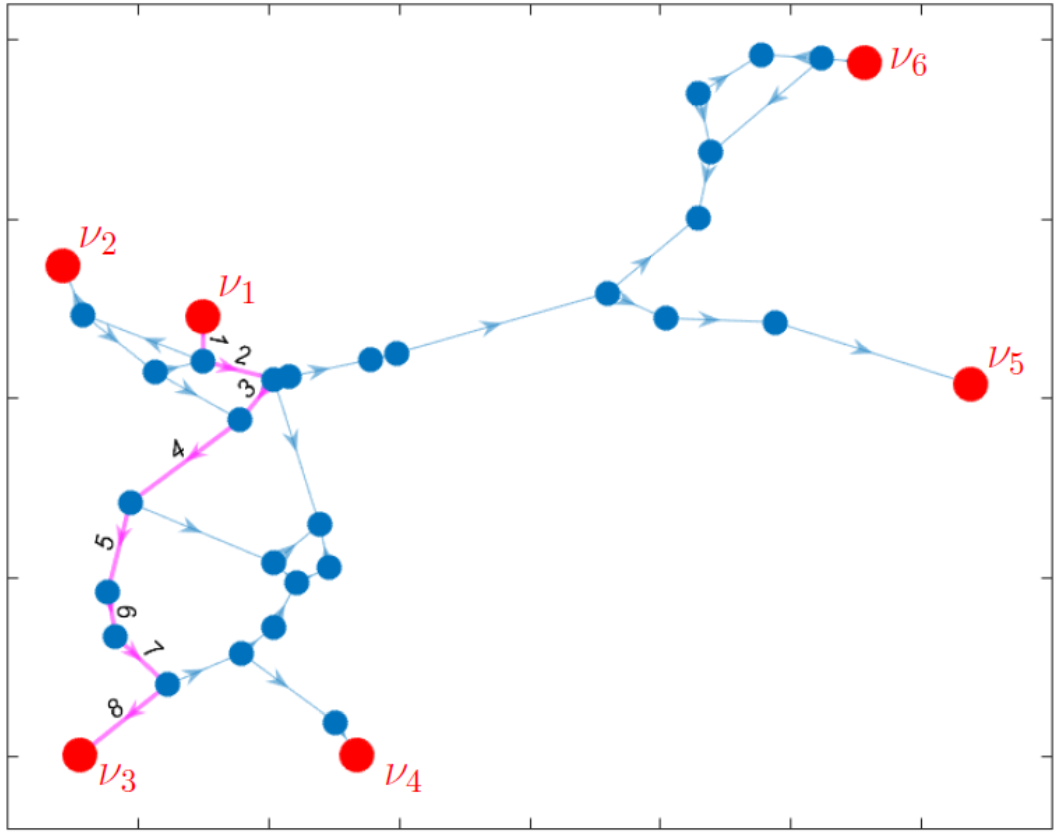}
{ 
\caption{Topology of large network with boundary nodes $\nu_i$, $i=1,...,6$ (red circles). The spatial domain of the pipes $\omega_j$, $j=1,...,8$, a path from $\nu_1$ to $\nu_3$, is colored in magenta.}\label{fig:gl38top}}
\end{figure}

\begin{figure}[tb]
\begin{tabular}{lrll | lrl}
\hspace{0.1cm}
&
\begin{minipage}{0.022\textwidth}
{\hspace{-0.3cm}
{\normalsize\rotatebox{90}{Profile $u_A$}  
}}
\end{minipage}
&

&
\hspace{-0.4cm}
\begin{minipage}{0.422\textwidth}
\hspace{1cm} {\large{\underline{Case A}}}\\[0.3em]
 \includegraphics[width=0.51\textwidth]{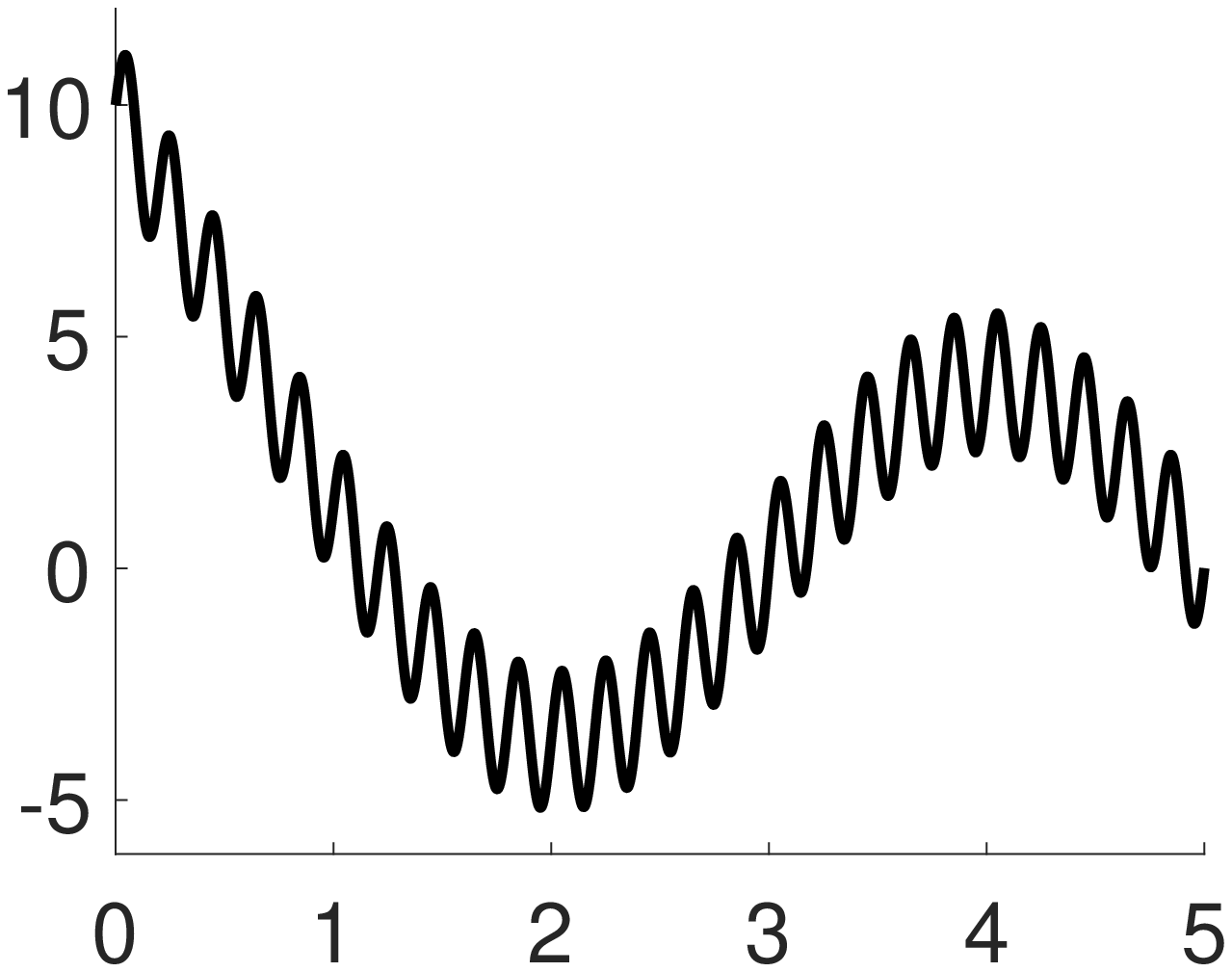} \, \\ \,
 $\text{\hspace{1cm}}$ { \normalsize Time $t$ [h]} 
 \end{minipage}
  &
 &
\begin{minipage}{0.022\textwidth}
{\hspace{-0.3cm}
{\normalsize\rotatebox{90}{Profile $u_B$}  
}}
\end{minipage}
& 
\begin{minipage}{0.45\textwidth}
\hspace{1cm}{\large\underline{Case B}} \\[0.3em] 
\includegraphics[width=0.51\textwidth]{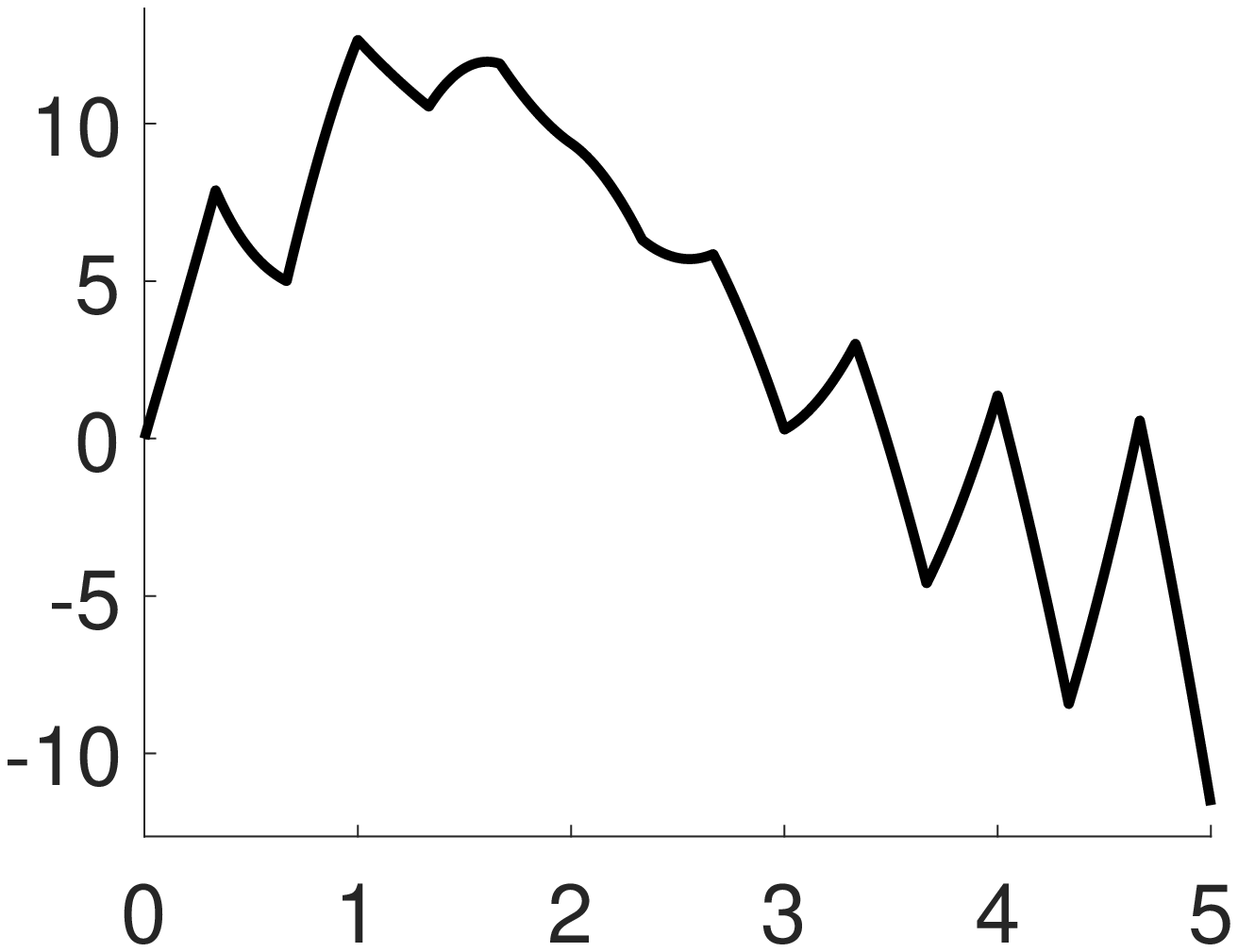}\\
 $\text{\hspace{1cm}}${\normalsize Time $t$ [h]}
\end{minipage}
\end{tabular}
\caption{Input profiles used for the boundary conditions of the large network.\label{fig:gl38-input}}
\end{figure}

The benefits of our structure-preserving approach are demonstrated for a benchmark with a realistic gas network, illustrated in Fig.~\ref{fig:gl38top}. The network is a slight modification of \cite[GasLib-40]{art:gaslib-2017}, cf. \cite{art:lilsailer-nlfow}. Let us note that a few compressors and valves would be needed in an actual application. We do not consider these active elements for ease of presentation, as they are prescribed by low-order algebraic models, which would not be regarded during the reduction process anyway. The network consists of 38 pipes with diameters $D^\Onepipe$ between 0.4 and 1~$[\mathrm{m}]$ and lengths $l^\Onepipe$ between 5 to 74~$[\mathrm{km}]$. The total pipe length is 1008~$[\mathrm{km}]$. Over the whole network, we set the friction factor to $\lambda = 0.008$.
At the six boundary nodes $\nu_i$, $i=1,...,6$, we prescribe the following boundary conditions for $t \in [0,5 t_*]$ and reference time $t_* = 1 [h]$,
\begin{align*}
	\rho(t,\nu_1) &= (65 + { u }(t/t_*))\,\rho_\star, \qquad \rho(t,\nu_2) = (50 + { u }(t/t_*))\,\rho_\star, \qquad   \rho(t,\nu_4) = (60 - { u }(t/t_*))\,\rho_\star, \\
	\rho(t,\nu_5) &= 60\,\rho_\star, \hspace*{2.65cm} \rho(t,\nu_6) = 45\, \rho_\star, \hspace*{2.25cm}  A  m (t,\nu_3) = - 100\, (Am)_\star,
\end{align*}
with the input profile ${ u } \in \{ u_A, u_B \}$ varied over two cases according to
\begin{align*}
\underline{\text{Case A:}} & \qquad u_A(t) = 6\exp\left( -\frac{3}{2} \, { t }\right)+4\cos\left(\frac{\pi}{2} { t }\right)+ \frac{3}{2} \, \sin\left( 10\pi \, { t }\right) \\[0.2cm]
\underline{\text{Case B:}} & \qquad
u_B(t) = 8 { t }^3\exp\left(-{ t }\right) - 4 \left( { t }-2 \right)f( 3 { t }), \qquad \text{ with } f(t) = 1- \left| (t \text{ mod } 2)-1 \right|,
\end{align*}
see Fig.~\ref{fig:gl38-input}. As initial condition the stationary solution belonging to the boundary conditions at $t=0$ is taken. The FOM is a system of dimension $N =10156$. The training of all reduction methods is realized with the trajectory of the {FOM} solution for Case~A with 1000 equally distributed snapshots. Thus, we refer to Case~A as the (perfectly) trained case and to Case~B as the not trained one. In general, one can expect model reduction methods to perform worse in a not trained case, even more so if they do not have good stability qualities. In the upcoming, we first show temporal and spatial visualizations of the reference solution, which is given by the {FOM}. Then we compare our structure-preserving model order reduction approach against a conventional method, and subsequently do the same for the complexity reduction step.

\begin{figure}[tb]
\begin{tabular}{rl|l}
\begin{minipage}{0.022\textwidth}
{\hspace{-0.3cm}
{\normalsize\rotatebox{90}{Density $\rho$}  \vspace{1cm}  \, 
\rotatebox{90}{Mass flow $Am$} 
}}
\end{minipage}
&
{\hspace{-0.4cm}
\begin{minipage}{0.45\textwidth}
\center
\quad {\large{\underline{Case A}}} \\[0.3em]
{\large{\underline{Pipe-end 2}}} \hspace{1.4cm}  {\large{\underline{Pipe-end 5}}}  \\
\hspace{0.3cm} 
\includegraphics[height = 0.4\textwidth, width = 0.44\textwidth]{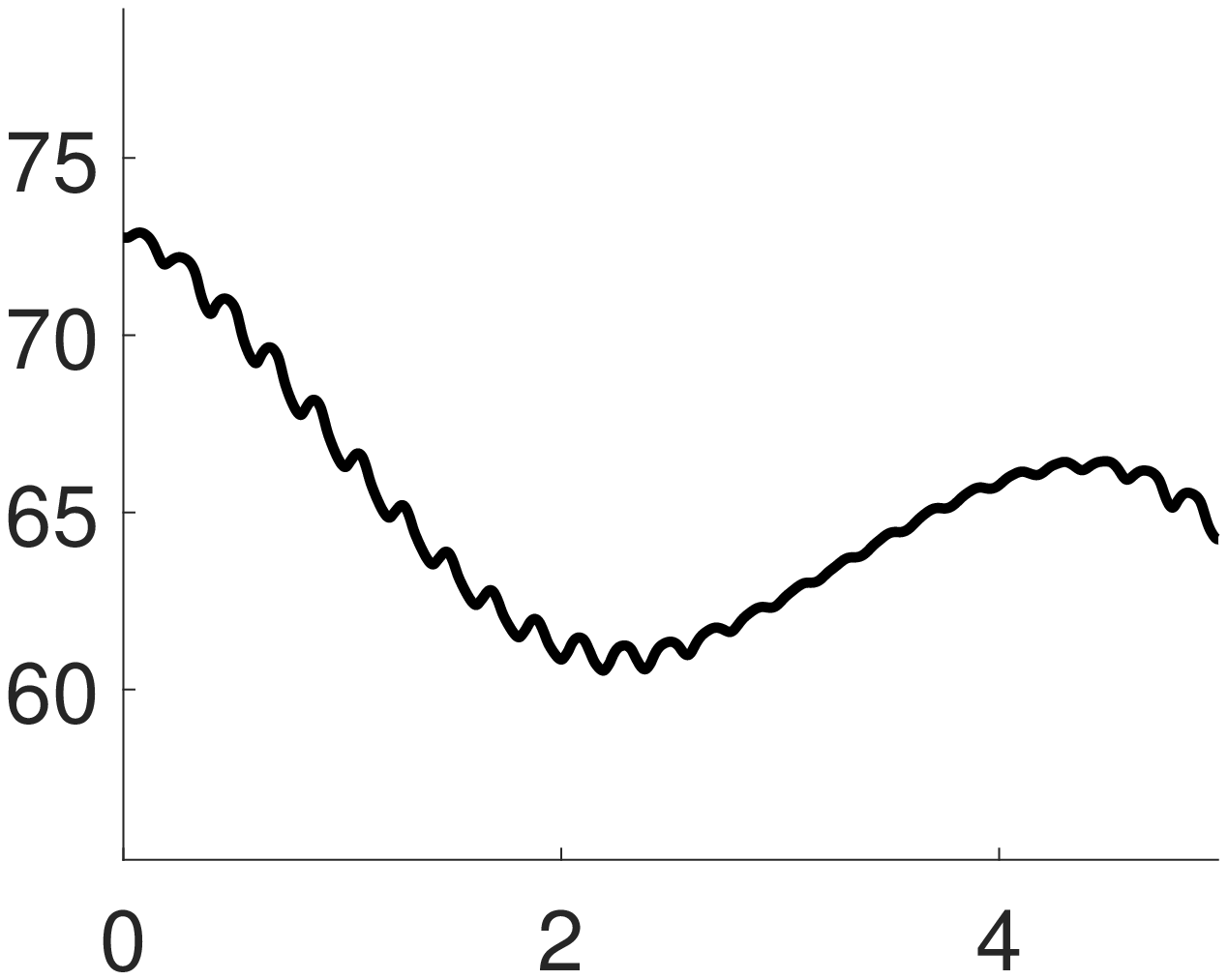}
 \includegraphics[height = 0.4\textwidth, width = 0.44\textwidth]{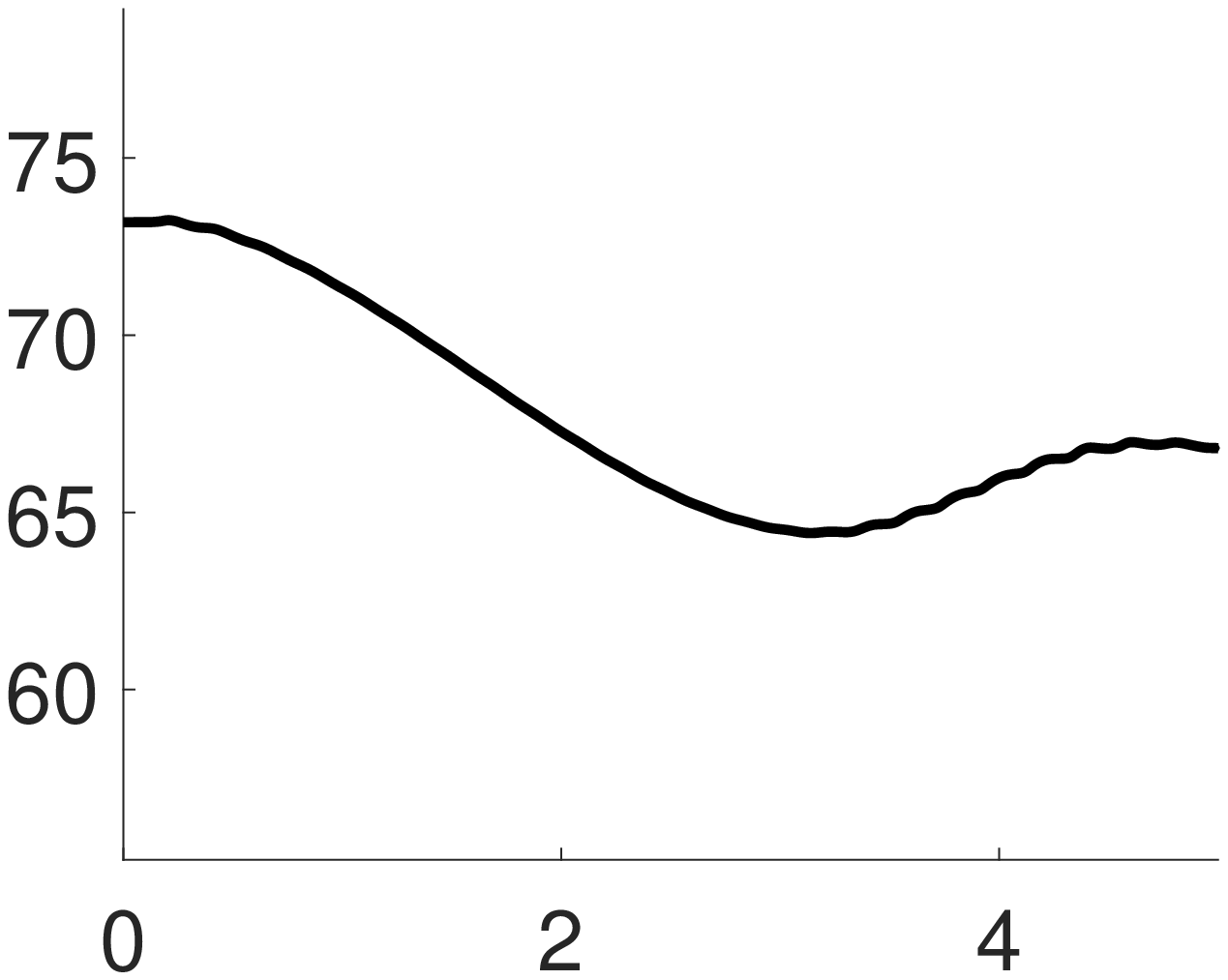}\\
 \hspace{0.1cm} 
\includegraphics[height = 0.4\textwidth, width = 0.44\textwidth]{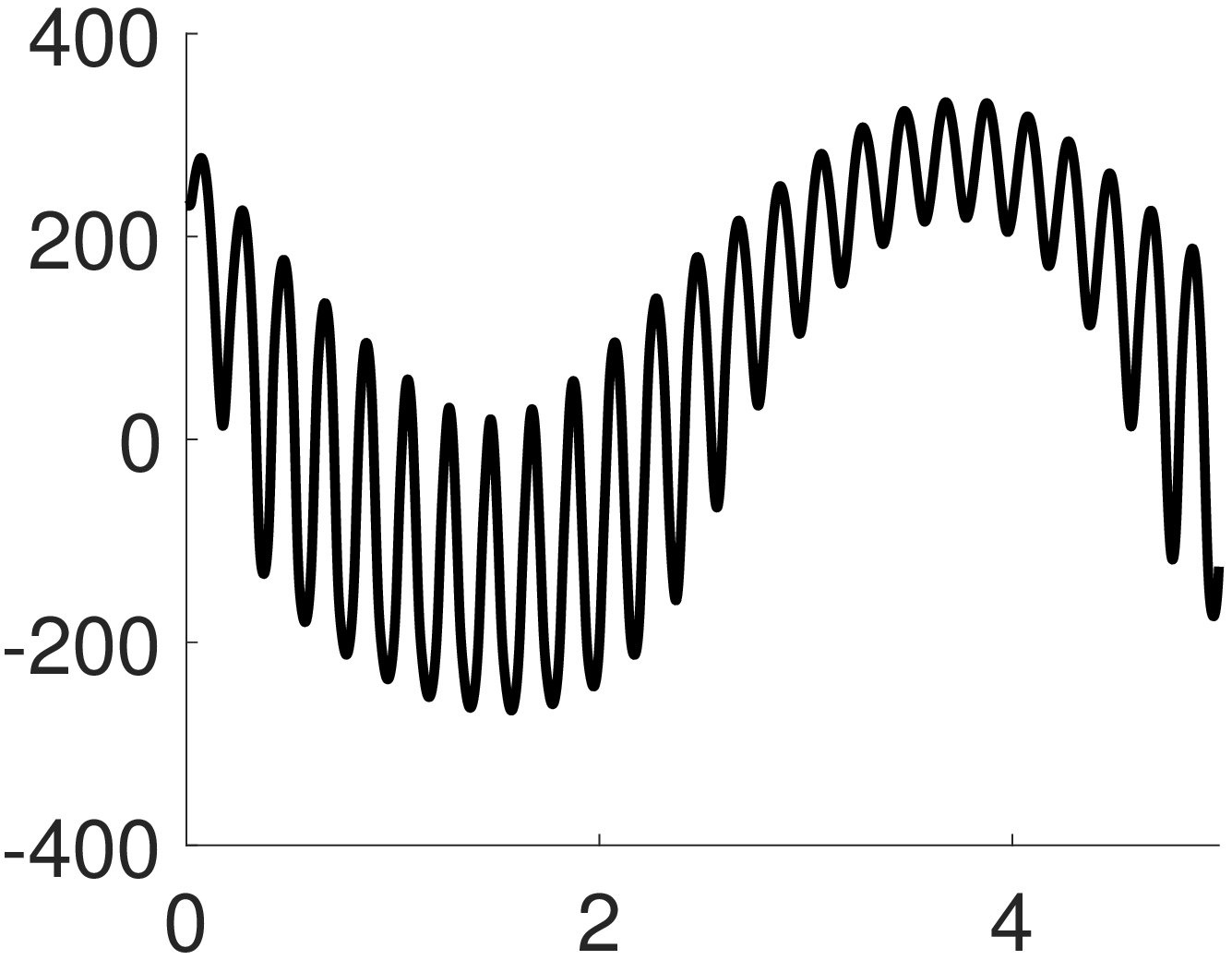} 
\includegraphics[height = 0.4\textwidth, width = 0.44\textwidth]{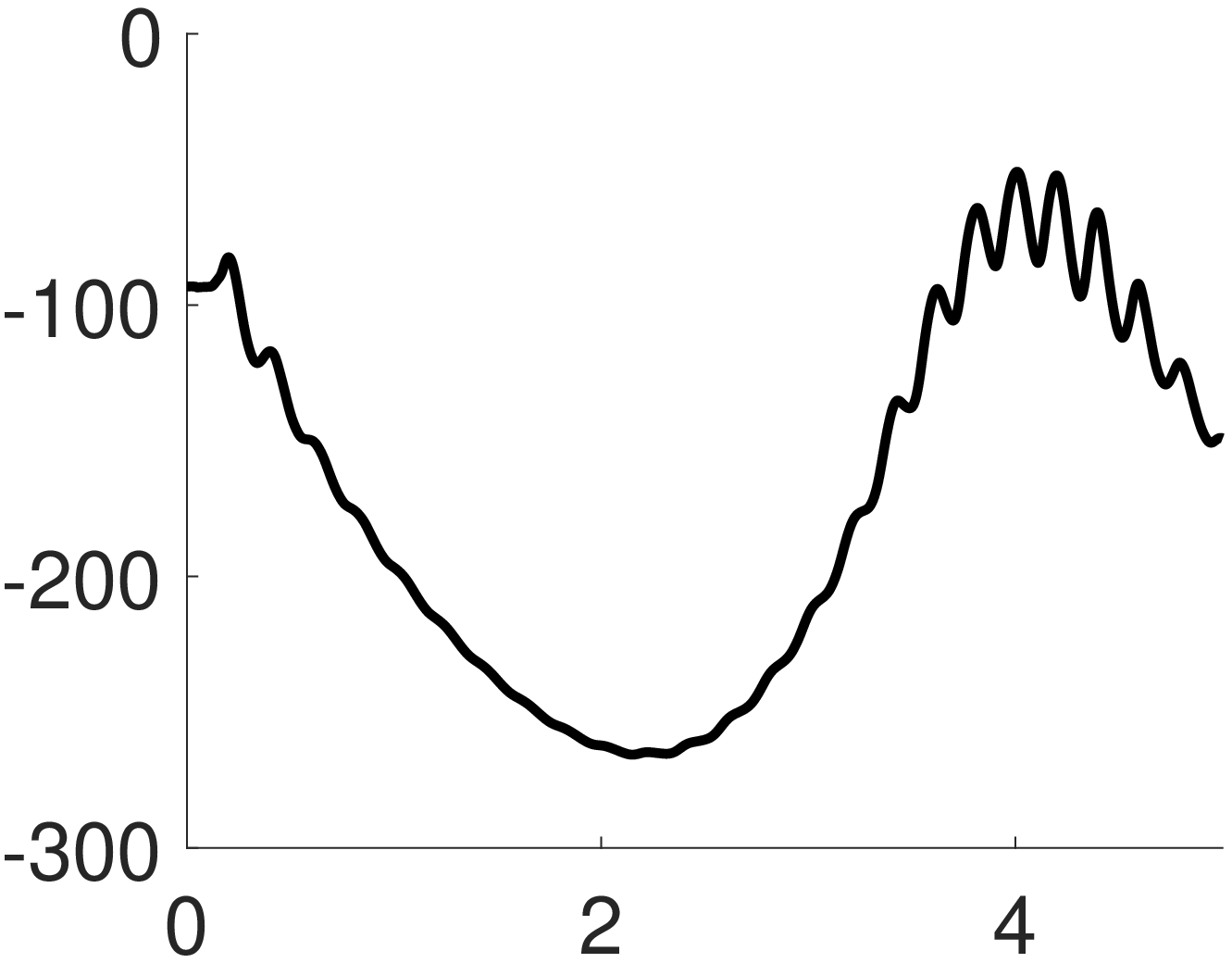}\\
 {\hspace{0.8cm} \normalsize Time $t [h]$}  \hspace{0.8cm}  {\normalsize \, $\text{\hspace{0.6cm} }$ Time $t [h]$} 
\end{minipage}
}
&
\begin{minipage}{0.45\textwidth}
\center
\quad {\large{\underline{Case B}}} \\[0.3em]
{\large{\underline{Pipe-end 2}}} \hspace{1.4cm}  {\large{\underline{Pipe-end 5}}}  \\
\hspace{0.3cm} 
\includegraphics[height = 0.4\textwidth, width = 0.44\textwidth]{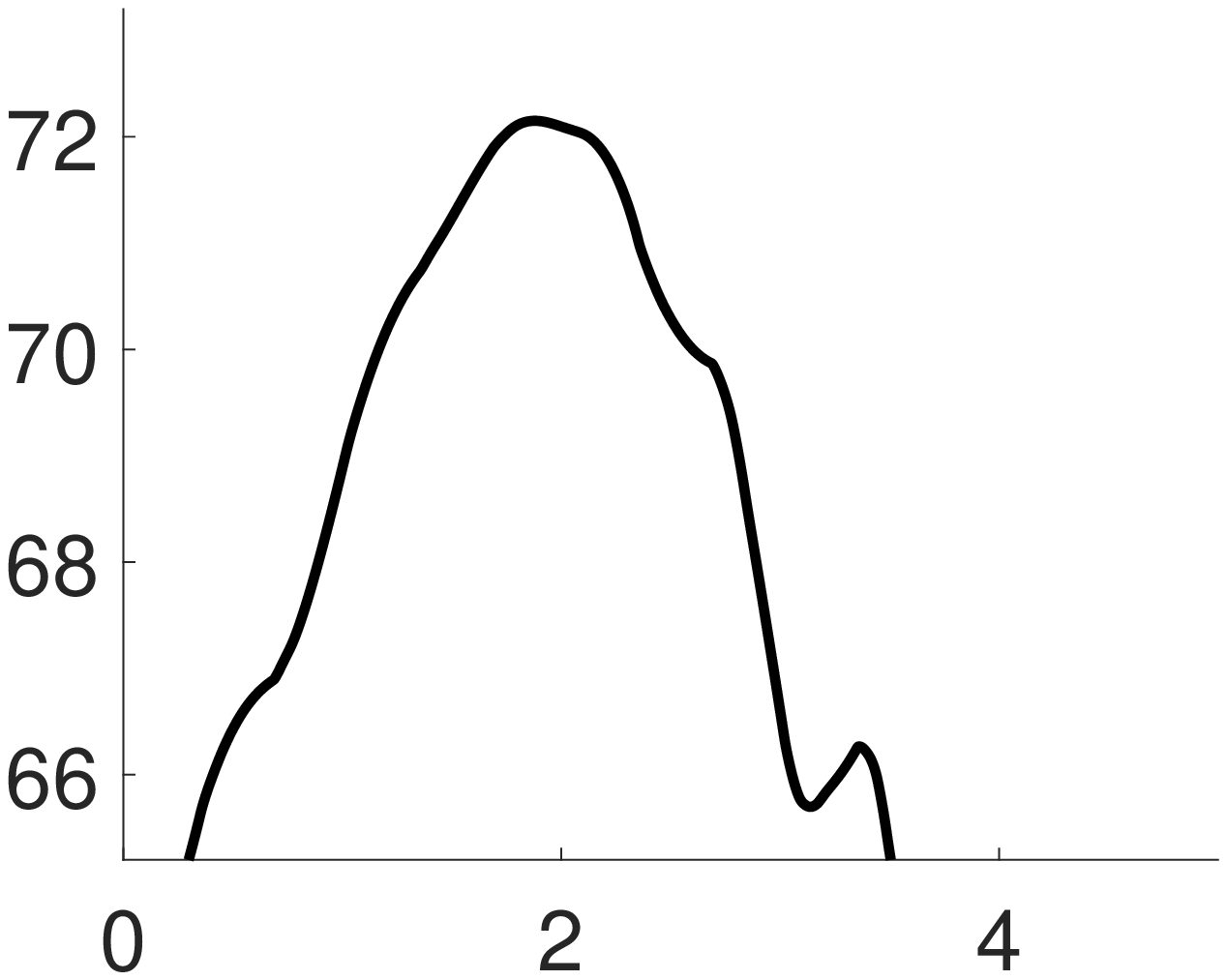} 
 \includegraphics[height = 0.4\textwidth, width = 0.44\textwidth]{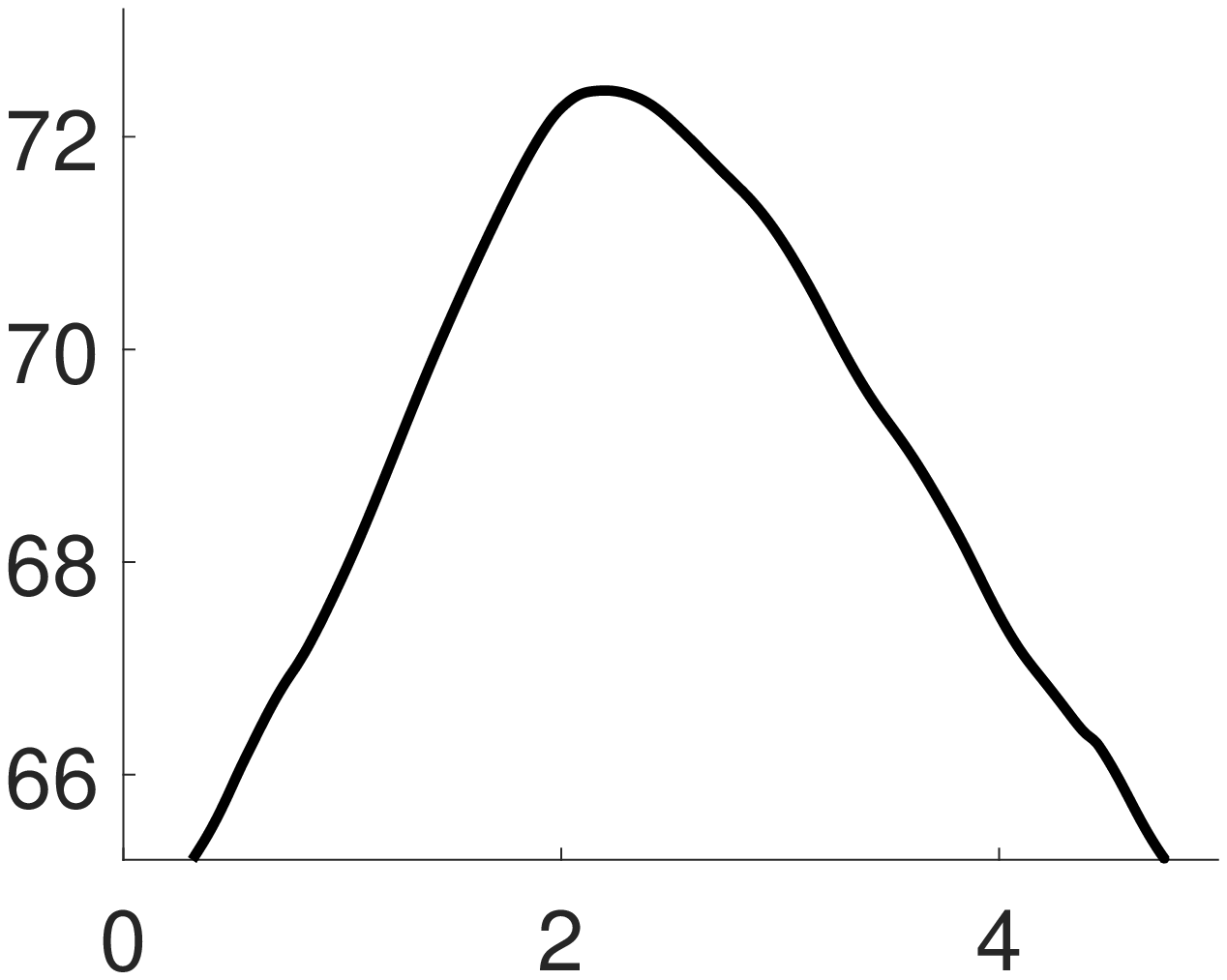}\\
 \hspace{0.1cm} 
\includegraphics[height = 0.4\textwidth, width = 0.44\textwidth]{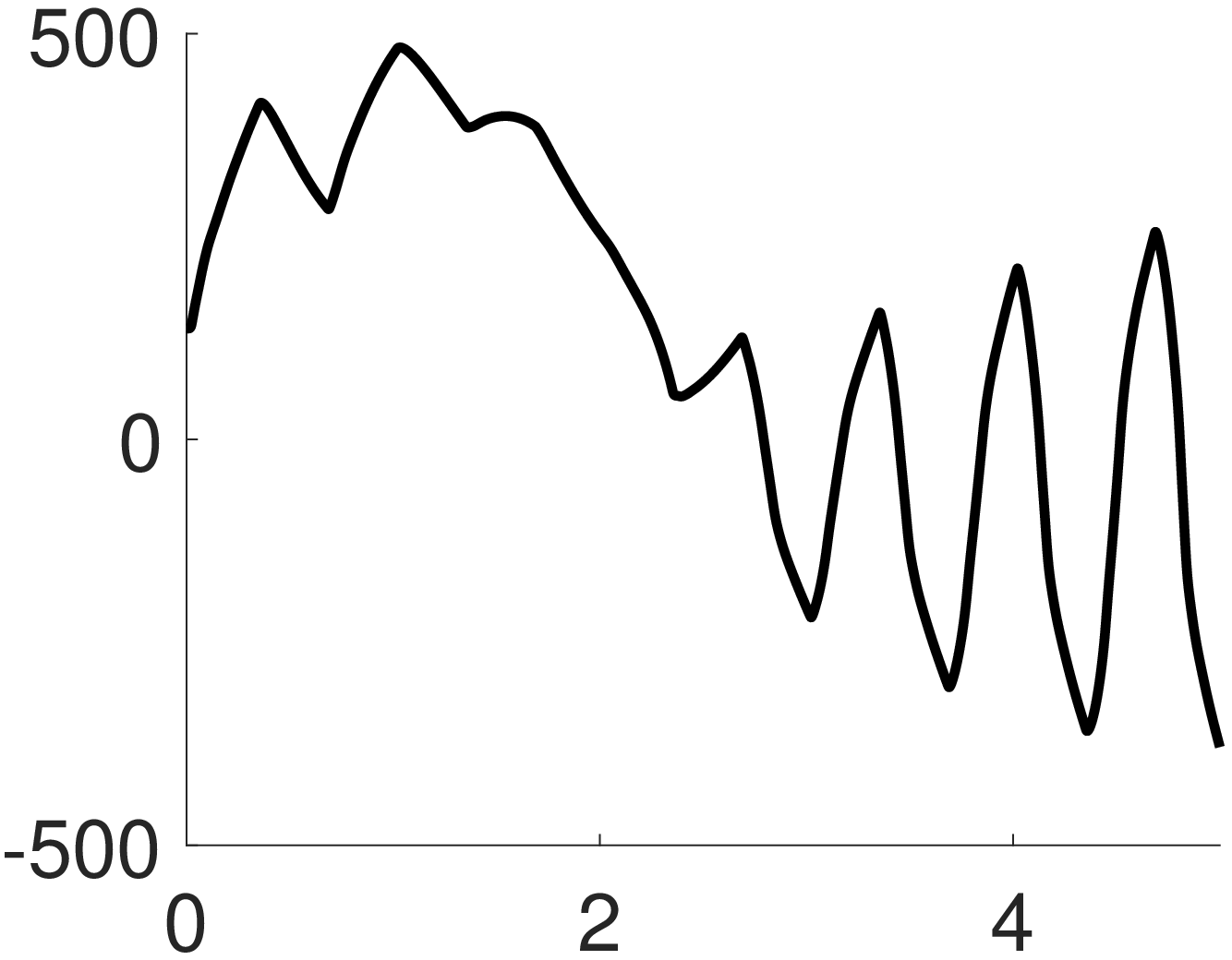} 
\includegraphics[height = 0.4\textwidth, width = 0.44\textwidth]{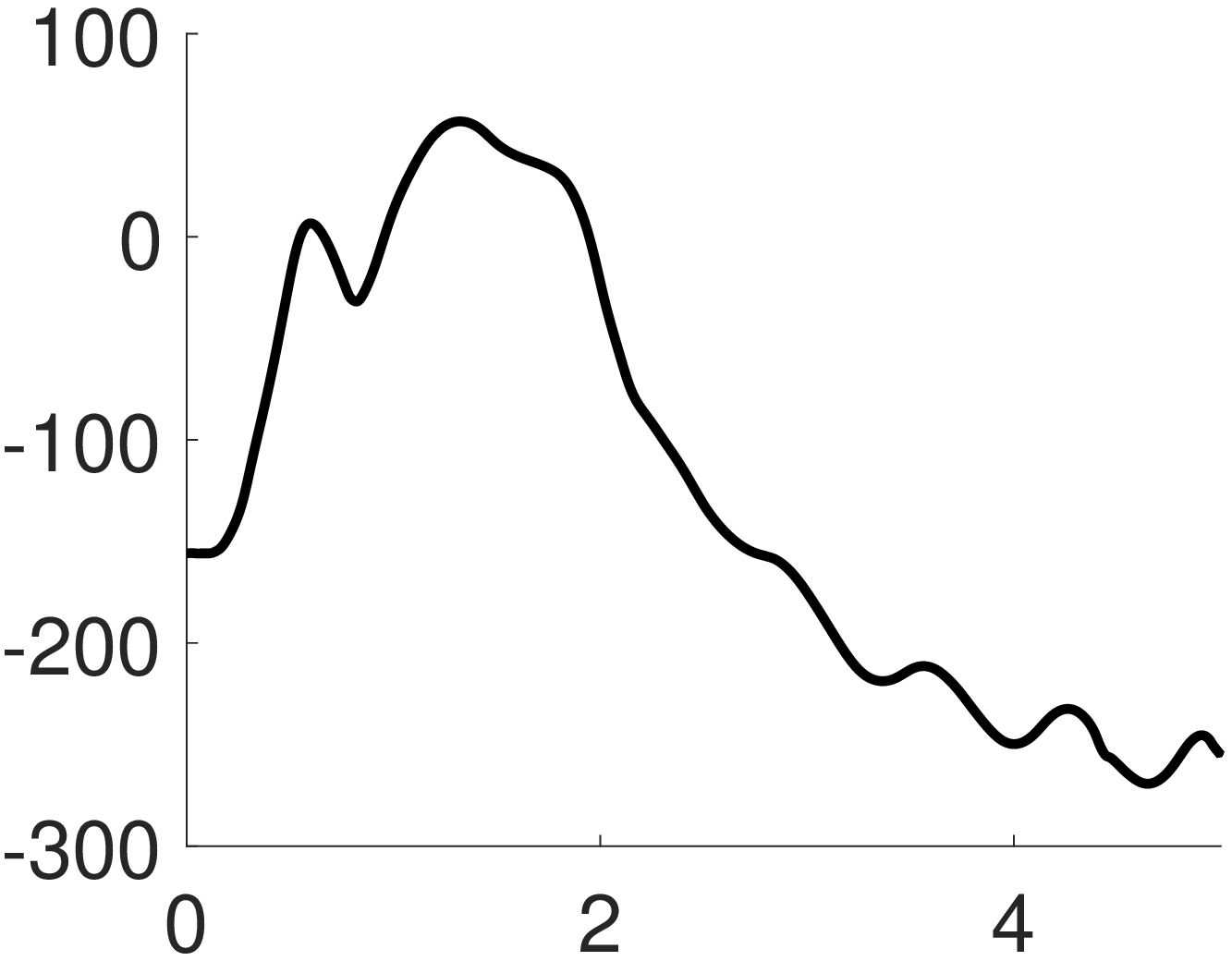}\\
 {\hspace{0.8cm} \normalsize Time $t [h]$}  \hspace{0.8cm}  {\normalsize \, $\text{\hspace{0.6cm} } $ Time $t [h]$} 
\end{minipage}
\end{tabular}
\caption{Large network, temporal representation of reference solution at the ends of pipes $\Onepipe_2$ and $\Onepipe_5$, respectively (divided by cases for input profile).\label{fig:gl38-timep}}
\end{figure}

\begin{figure}[tb]
\begin{tabular}{rll|l}
\begin{minipage}{0.022\textwidth}
{\hspace{0.3cm}
{\normalsize\rotatebox{90}{Density $\rho$}\\ \vspace{2.0cm}\\
\rotatebox{90}{Mass flow $A  m$\hspace*{1cm}} 
}}
\end{minipage}
&
{\hspace{-0.4cm}
\begin{minipage}{0.43\textwidth}
\center
\hspace{0.5cm} {\large{\underline{Case A}}} \\
\hspace{0.25cm} \includegraphics[height = 0.6\textwidth, width = 0.95\textwidth]{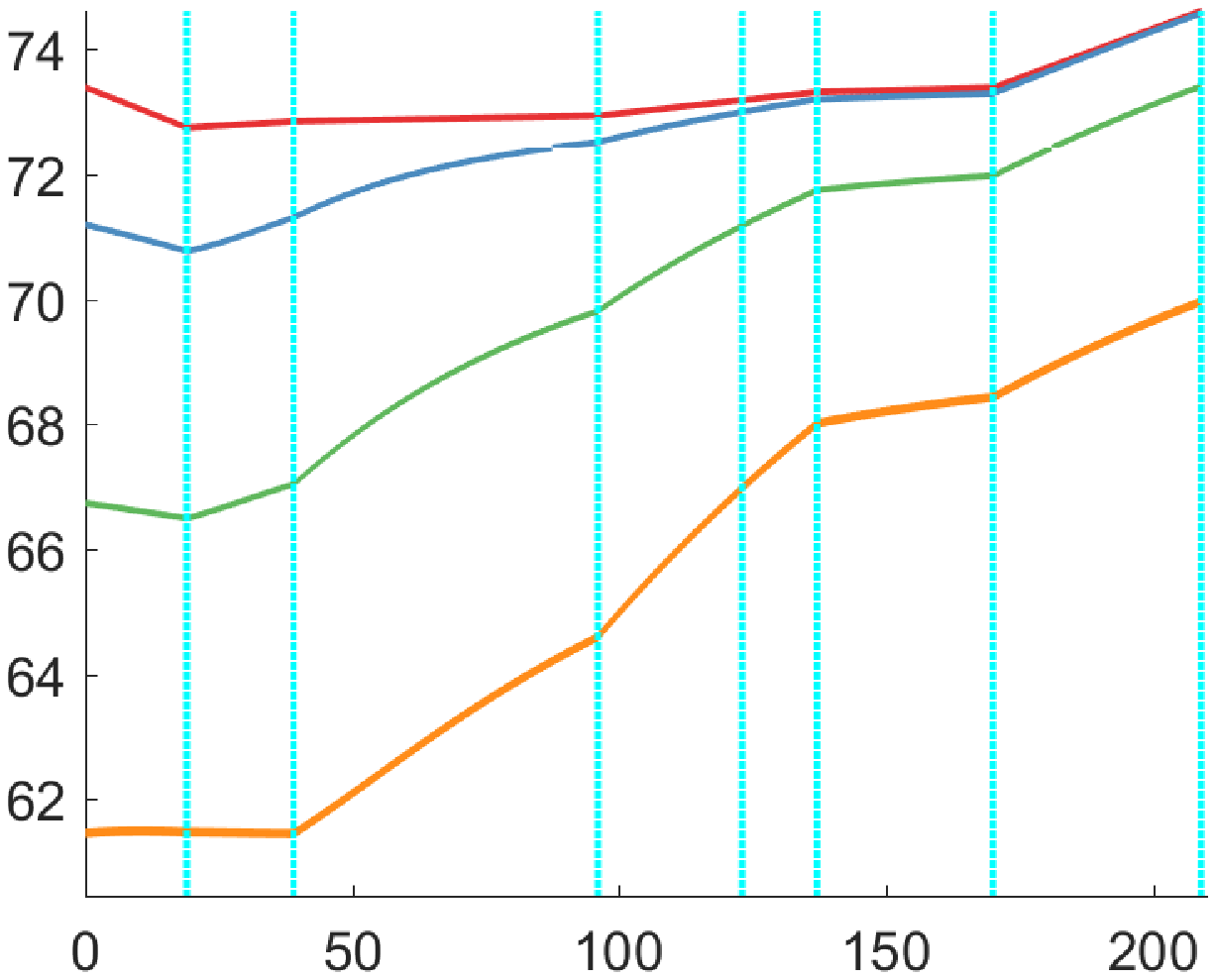} \\
\includegraphics[height = 0.6\textwidth, width = 0.83\textwidth]{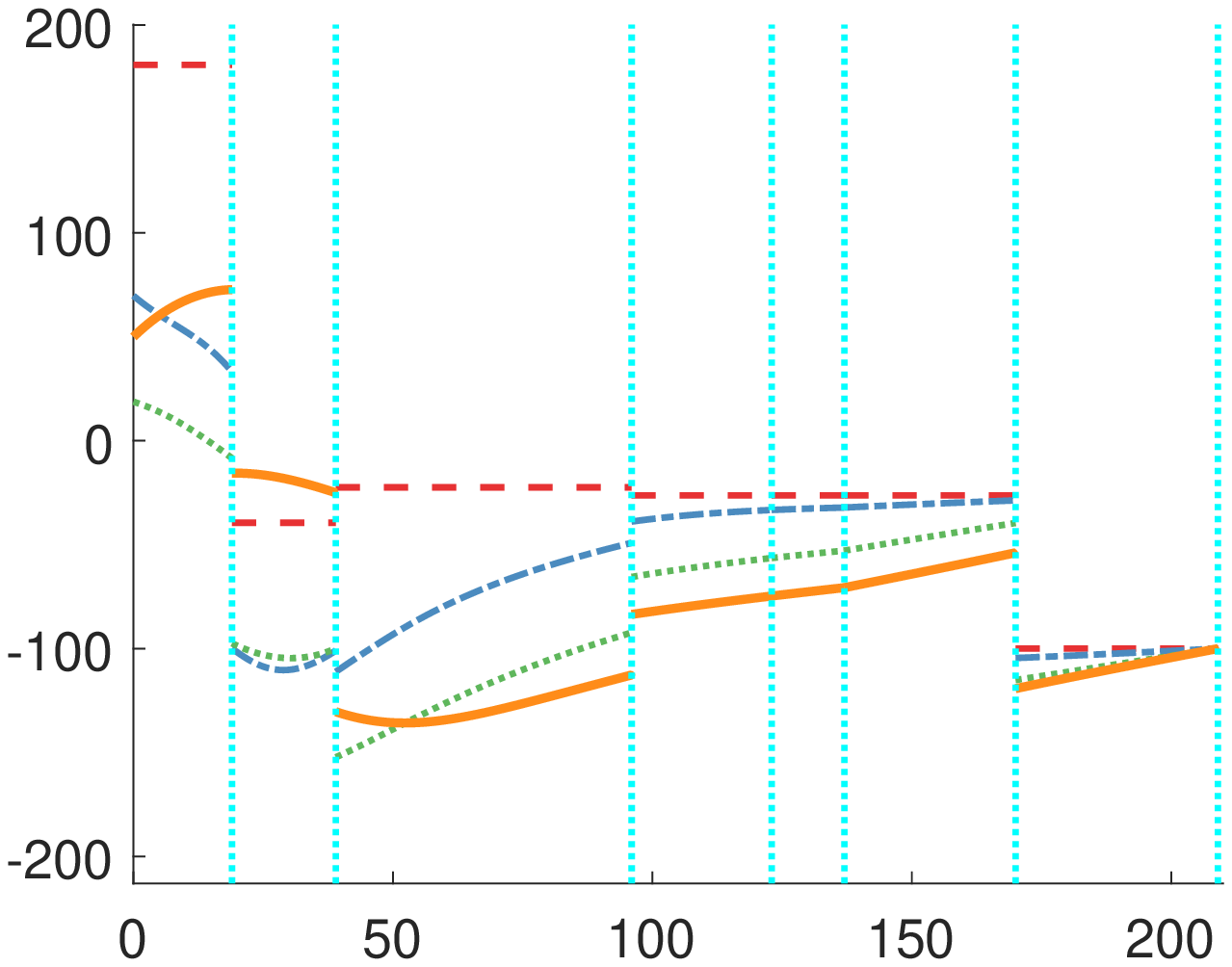} \\
 \hspace{0.5cm} {\normalsize Space $x/\overline{km}$}
\end{minipage}
}
&
\begin{minipage}{0.00\textwidth}
\end{minipage}
&
\begin{minipage}{0.43\textwidth}
\center
\quad {\large{\underline{Case B}}}\\
\hspace{0.1cm} \includegraphics[height = 0.6\textwidth, width = 0.95\textwidth]{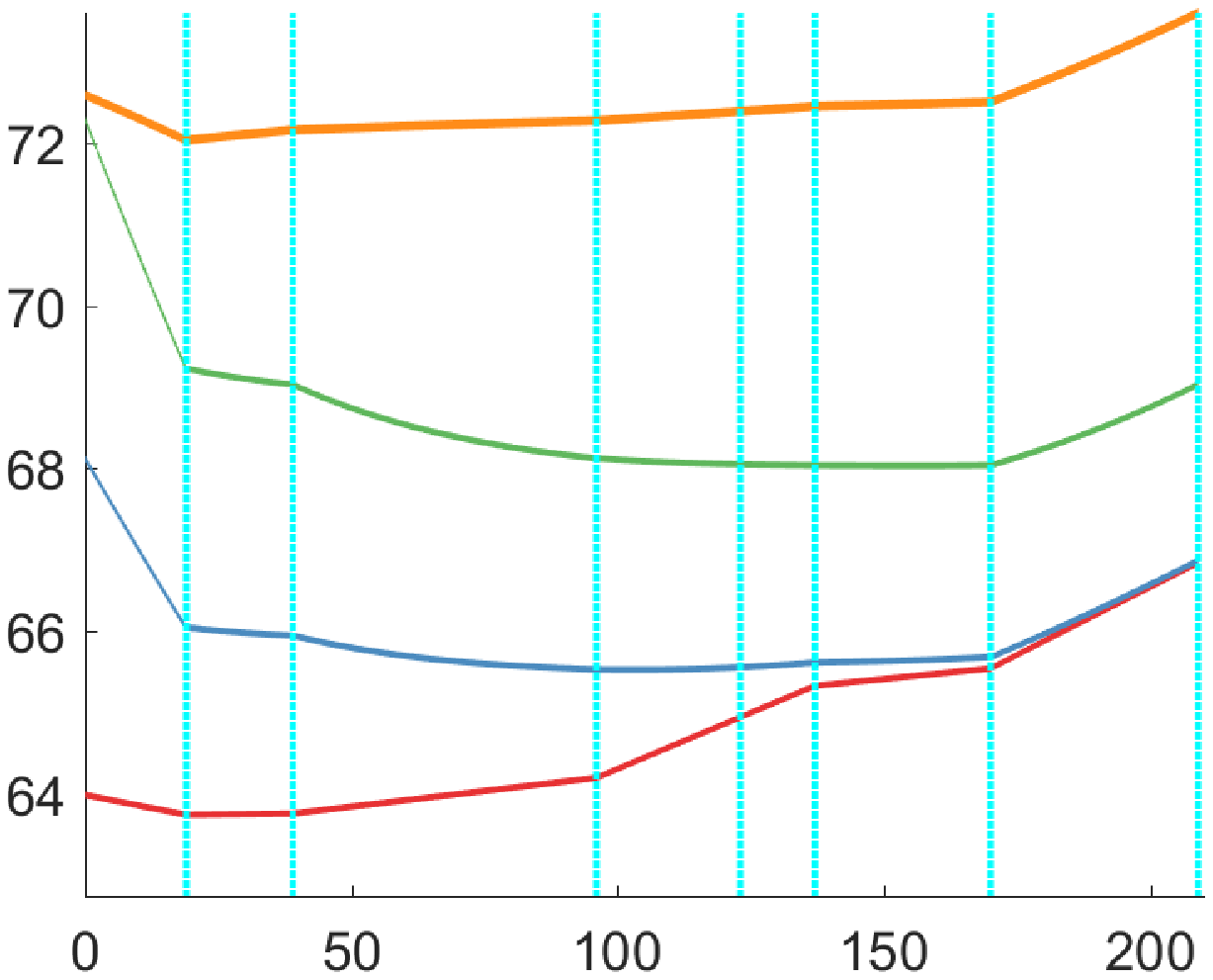} \\
\includegraphics[height = 0.6\textwidth, width = 0.83\textwidth]{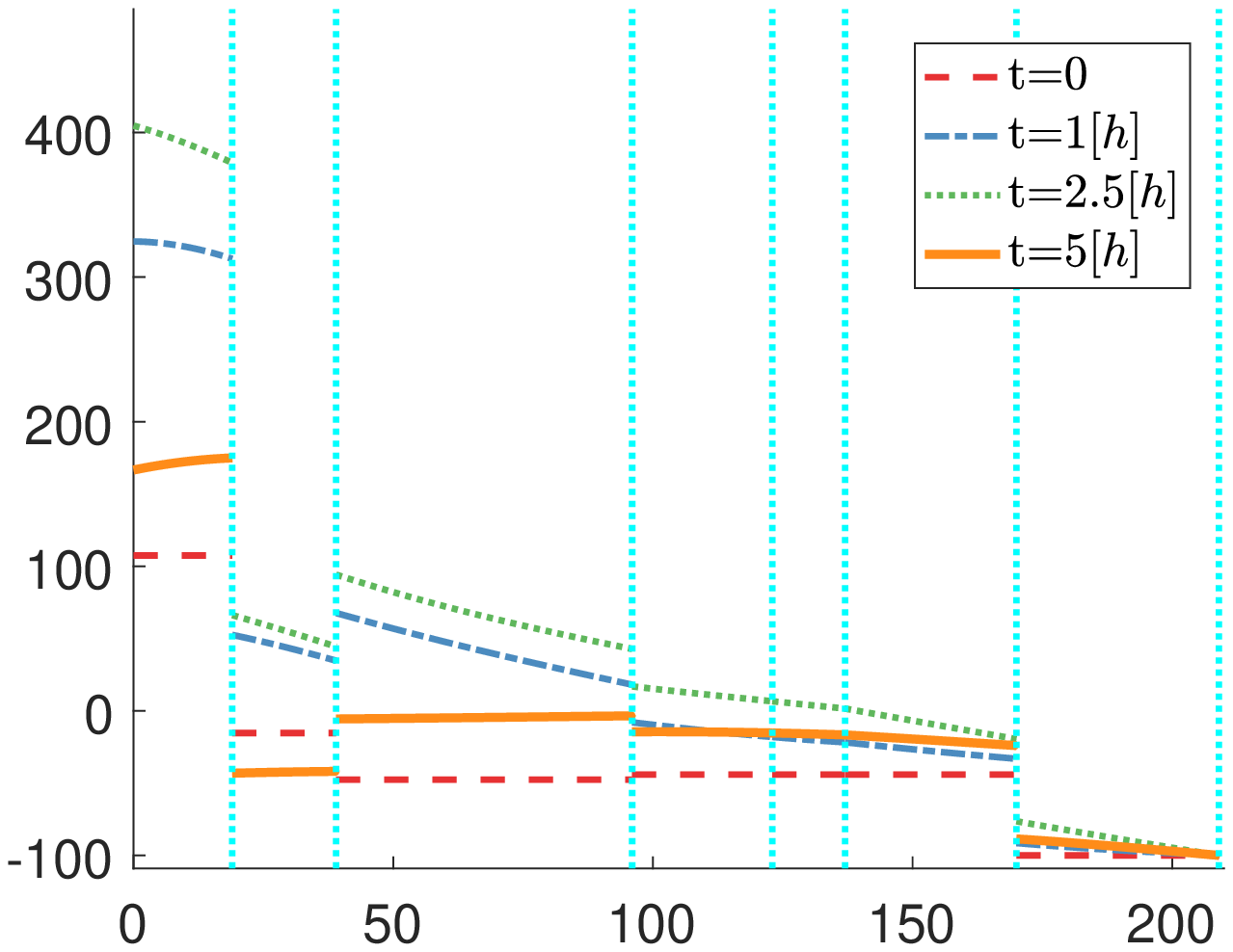} \\
 \hspace{0.5cm} {\normalsize Space $x/\overline{km}$}
\end{minipage}

\end{tabular}
\caption{Large network, spatial representation of reference solution with domain representing pipes 2-8 (marked in magenta in Fig.~\ref{fig:gl38top}). The cyan vertical lines indicate junctions.\label{fig:gl38-space-fom}}
\end{figure}

The time response of the {FOM} at the ends of pipe 2 and pipe 5 is shown in Fig.~\ref{fig:gl38-timep}. Similarly as for the academic test scenario, we observe notable damping effects. The time response is smoothed out stronger at pipe-end 5 than at pipe-end 2, because the former has a larger distance to the driving boundary nodes. In Fig.~\ref{fig:gl38-space-fom}, a spatial representation of the solution along the path from pipe 2 to pipe 8 (marked in magenta in Fig.~\ref{fig:gl38top}) is shown  for different time instances. The mass flows show to be discontinuous over junctions with more than two pipes meeting, as to be expected from the model. Apart from that, the spatial representations are rather smooth, which indicates that the solution manifold can be approximated well in low dimensions and model reduction methods could potentially be efficient.

%

\begin{figure}[h]
\begin{tabular}{rll|l}
\begin{minipage}{0.00\textwidth}
\end{minipage}
&
\begin{minipage}{0.022\textwidth}
{
{\normalsize\rotatebox{90}{Relative error $E_T$}\\ 
}}
\end{minipage}
&
{\hspace{-0.4cm}
\begin{minipage}{0.43\textwidth}
\center
\hspace{0.5cm} {\large{\underline{Case A (trained)}}} \\
\hspace{0.1cm} \includegraphics[height = 0.77\textwidth, width = 0.9\textwidth]{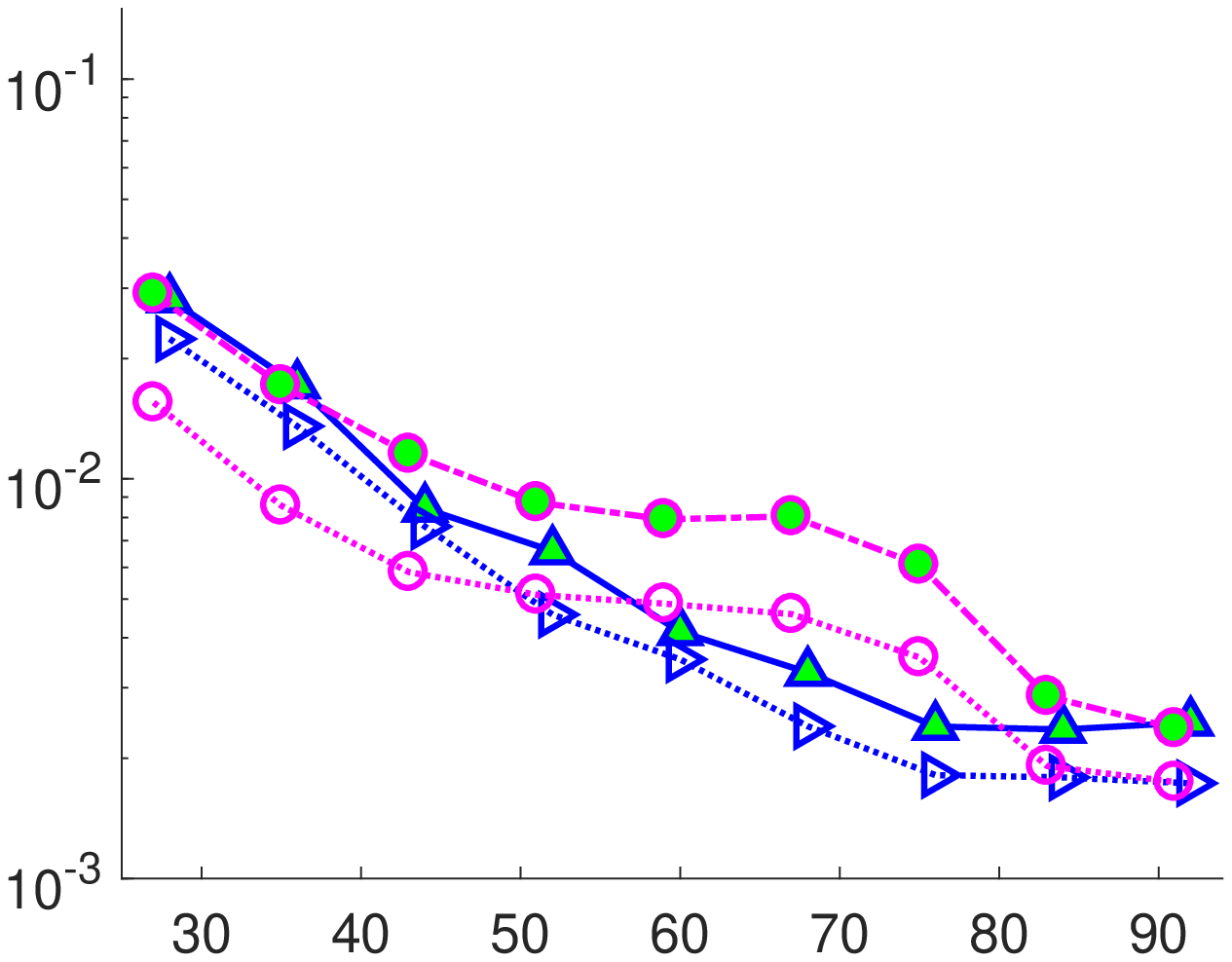}\\
 \hspace{0.5cm} {\normalsize Reduced dimension $n$}
\end{minipage}
}
&
{\hspace{-0.4cm}
\begin{minipage}{0.43\textwidth}
\center
\hspace{0.5cm} {\large{\underline{Case B (not trained)}}} 
\hspace{0.1cm} \includegraphics[height = 0.77\textwidth, width = 0.9\textwidth]{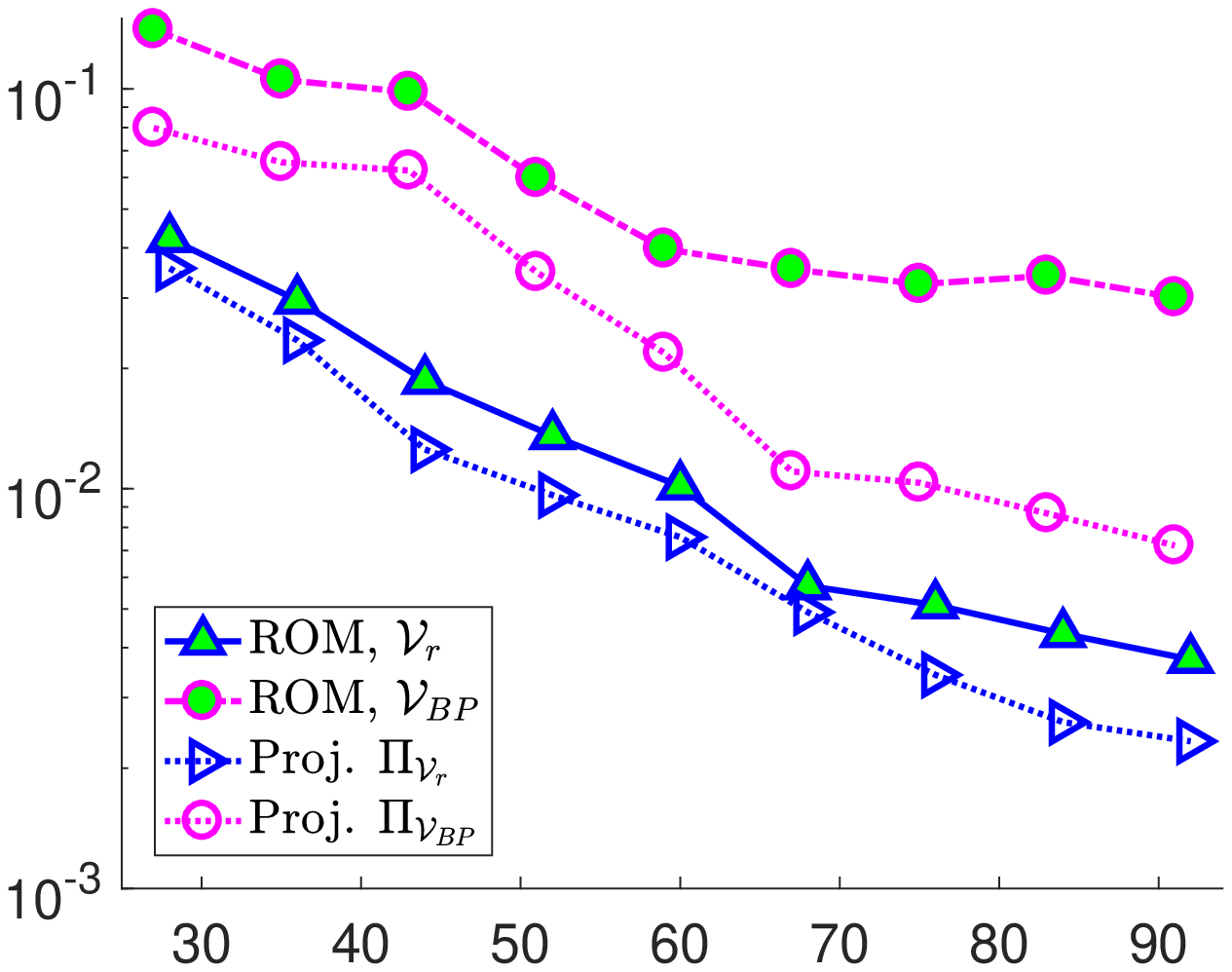}\\
 \hspace{0.5cm} {\normalsize Reduced dimension $n$}
\end{minipage}
}

\end{tabular}
\caption{Large network, errors of {ROM}s and orthogonal projection onto the respective reduction space (Proj. $\Pi_{\funSpace{V}_r}$/$\Pi_{\funSpace{V}_{BP}}$), using our structure-preserving basis $\funSpace{V}_r$ and block-structured {POD}-basis $\funSpace{V}_{BP}$. Dimension FOM: $N=10156$.\label{fig:gl38poder}}
\end{figure}

We compare our proposed model order reduction against a conventional block-structured {POD}, which is very similar to the approaches \cite{art:himpe-gasmor21,art:mor-gas-grundel-jansen} used in the gas network context. The reduction basis $\funSpace{V}_{BP}$ for the conventional method is obtained by separately applying a principal component analysis to the density- and mass flux-snapshots and extracting reduction spaces of the same dimension for them  ($n_1=n_2$). This approach is preferable to applying {POD} onto the full state, but it is not structure-preserving in the sense of Assumption~\ref{assum:compatV1V2}. In order to investigate the fidelity of the reduction spaces, independent of the stability properties of the reduced models, we additionally consider the $\funSpace{L}^2$-projections of the {FOM} solutions onto the {ROM} spaces. The respective projectors are denoted by $\Pi_{\funSpace{V}_{r}}$ and $\Pi_{\funSpace{V}_{BP}}$, and they yield the pure projection errors. 
The relative errors are shown in Fig.~\ref{fig:gl38poder}. There is one scenario, where the conventional {POD} shows slightly better results than our method. That is, better projection errors are observed for the perfectly trained case for the reduced dimension $n$ small (Fig.~\ref{fig:gl38poder}-\textit{left}, $n\leq 50$). This is not surprising, as taking the $\funSpace{L}^2$-error favors $\funSpace{V}_{BP}$. Recall that our method is derived using another norm and regards compatibility conditions. Despite that, our reduction basis $\funSpace{V}_{r}$ has slightly better projection errors in the trained case when the dimension is chosen to be larger, and shows significantly smaller projection errors for all parameter choices in the not trained case. The theoretical advantage our approach has over the other is that it fulfills an optimality condition in the full state. We suppose this is also the reason for the comparably low projection errors we observe. 
%
Practically more relevant is the comparison of the reduction errors. In this respect, our method outperforms the conventional reduction method much clearer, and we observe better errors for all dimensions $n$. The difference is most evident for the non-trained case, where the second smallest reduced order model with $\funSpace{V}_r$ is of higher fidelity than any of the reduced models for $\funSpace{V}_{BP}$, see Fig.~\ref{fig:gl38poder}-\textit{right}.
 Further, we want to highlight the small gap between reduction- and  projection-error for our proposed method compared to the conventional one, which strongly indicates the superior stability and robustness of our approach.

\begin{figure}[tb]
\begin{tabular}{rll|l}

\begin{minipage}{0.00\textwidth}
\end{minipage}
&
\begin{minipage}{0.022\textwidth}
{
{\normalsize\rotatebox{90}{Relative error $E_T$}\\ 
}}
\end{minipage}
&
{\hspace{-0.4cm}
\begin{minipage}{0.43\textwidth}
\center
\hspace{0.5cm} {\large{\underline{Case A (trained)}}} \\[0.2em]
\hspace{0.1cm} \includegraphics[height = 0.77\textwidth, width = 0.9\textwidth]{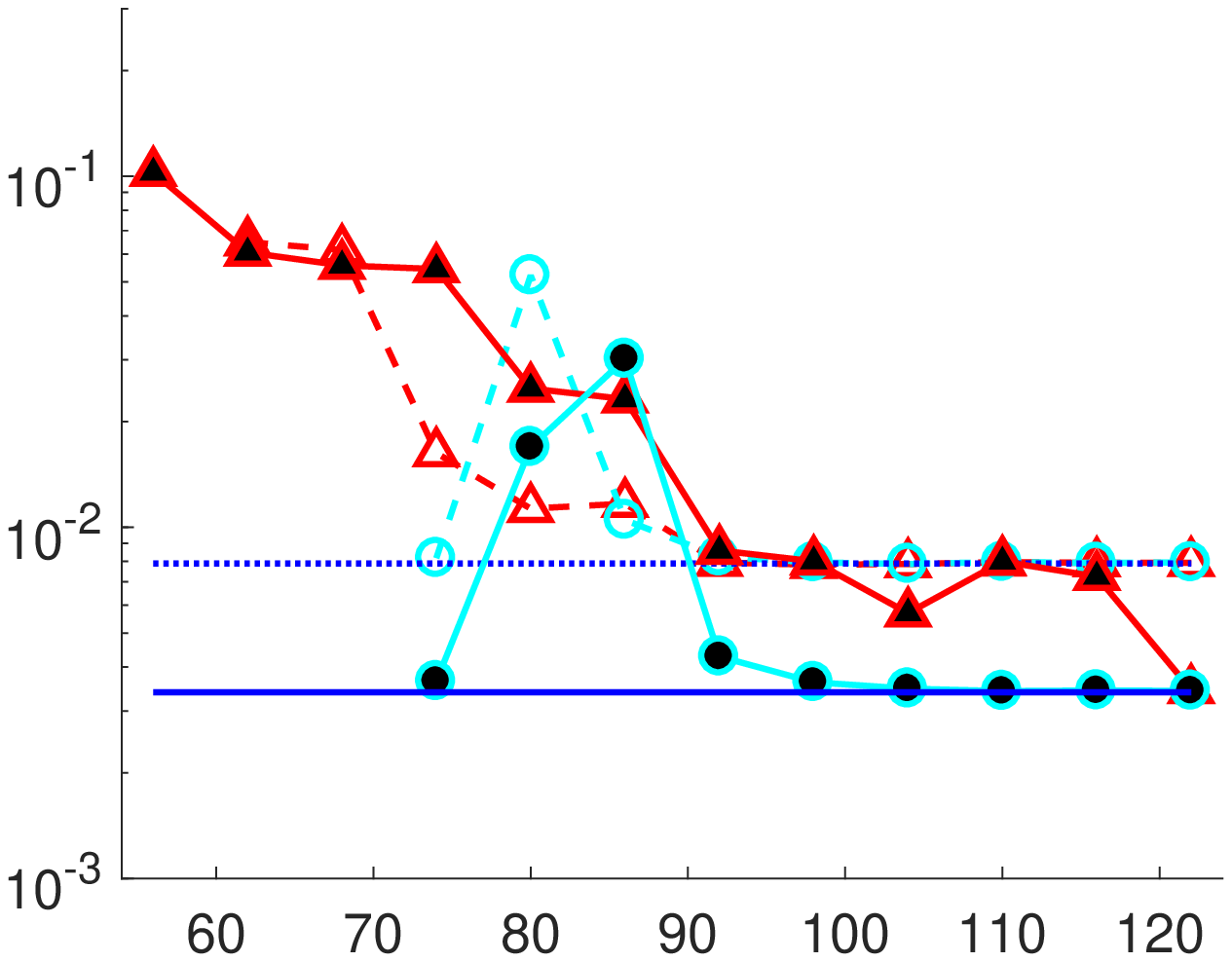}\\
 \hspace{0.5cm} {\normalsize Parameter $n_c$} 
\end{minipage}
}
&
{\hspace{-0.4cm}
\begin{minipage}{0.43\textwidth}
\center
\hspace{0.5cm} {\large{\underline{Case B (not trained)}}} \\[0.2em]
\hspace{0.1cm} \includegraphics[height = 0.77\textwidth, width = 0.9\textwidth]{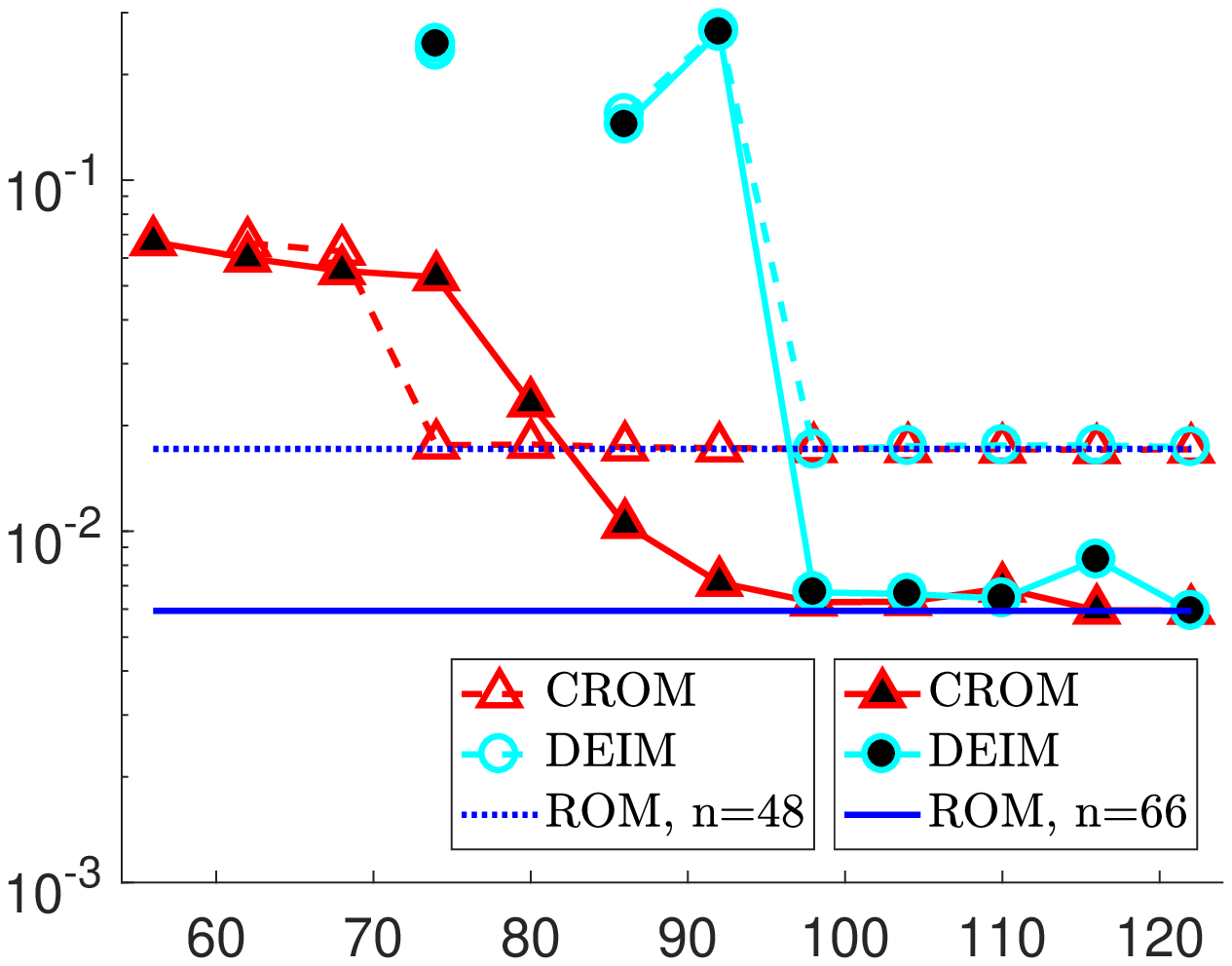}\\
 \hspace{0.5cm} {\normalsize Parameter $n_c$}
\end{minipage}
}

\end{tabular}
\caption{Large network, error of proposed complexity reduction {CROM}, non-structure-preserving {DEIM}. Underlying are {ROM}s of dimension $n = 48$ \textit{(dashed-dotted lines)} and $n = 66$ \textit{(solid lines)}.\label{fig:gl38crer}}
\end{figure}

\begin{table}[tb]
\renewcommand{\arraystretch}{1.2}
\small 
\begin{center}
\begin{tabular}{l || c|c|c|c|c|c|c|c|c|c |c  |c}
Parameter $n_c$ & 56   & 62   & 68   & 74   & 80   & 86   & 92   & 98  & 104  & 110  & 116  & 122 \\ 
\hline
\hline
 \texttt{cond}($\MassBoth_c$)   &  85.28   &      55.94   &      12.59    &      9.56     &     5.27   &      2.42     &     2.44     &     2.43     &     1.83      &    2.14     &     1.86 	&  1.20 \\ 
\end{tabular}
\medskip
\caption{Large network (case-independent) condition number \texttt{cond}($\MassBoth_c$) w.r.t.~the spectral norm using $n=66$. An almost monotone decrease can be observed. The measure is strongly related to the compatibility condition $\sigma(\MassBoth_{c}) \subset [\tilde{C}^{-2},\tilde{C}^2]$ for CROM.
\label{table:gl38-cond}
}
\end{center}
\renewcommand{\arraystretch}{1.2}
\small 
\begin{center}
\begin{tabular}{l || c|c|c|c|c|c|c|c|c|c |c  |c}
Parameter $n_c$ & 56   & 62   & 68   & 74   & 80   & 86   & 92   & 98  & 104  & 110  & 116  & 122 \\ 
\hline
\hline
CROM  &  126.7 &  125.8 &  127.5 &  132.5 &  132.6 &  135.9  & 137.8 &  140.6 & 141.8 &  144.4 &  147.1 &  152.1 \\ 
\hline
DEIM    &  - &    - & - & 432.0 &- & 470.0 &  451.2 &   454.7   &     461.2 &  466.6 & 473.4 & 482.4 \\
\end{tabular}
\medskip
\caption{Large network, Case B. Runtimes (in seconds) for the complexity reduced models with underlying {ROM} dimension $n =66$. For four choices of $n_c$ the simulation for DEIM fails. The simulation time of FOM is $2253.6$ seconds. 
\label{table:gl38-runtime}
}
\end{center}
\vspace{-2.1em}
\end{table}

As regards the comparisons for the complexity reduction step, the starting point is a ROM obtained by our proposed model order reduction method. Our quadrature-type complexity reduction is compared to a non-structure preserving alternative, which we representatively choose as the discrete empirical interpolation method (DEIM) \cite{art:deim-state-space-err,art:deim-introduction}. For convenience, we denote the DEIM-space dimension used in the  training as $n_c$. Note that the parameter in our {CROM}, which we also denote by $n_c$, has a quite different interpretation, which is the number of nonlinear integrals that need to be evaluated. The DEIM is not directly related to the integral expressions, but rather to the algebraic representation of the system, and needs several nonlinear integral evaluation for each $n_c$, which makes it less online-efficient in our setting. We apply the DEIM to each nonlinearity separately, i.e, given the algebraic representation of the {ROM} (System~\ref{sys:rom-coord}), that means the terms $\bv f_r^\alpha$, $\bv f_r^\gamma$ and $\tilde{\bv f}_r^\beta$ are separately complexity reduced. The latter term $\tilde{\bv f}_r^\beta$ is defined by $\bv f_r^\beta = \bt J^T \tilde{\bv f}_r^\beta$, and it is introduced as the resulting DEIM approach regards some of the structural properties of our model problem, cf., Remark~\ref{rem:useOfDeim} and \cite{art:lilsailer-nlfow,phd:liljegren}. But let us emphasize that this DEIM approach still cannot guarantee port-Hamiltonian structure or a provable energy bound, in contrast to our quadrature-based approach. We test the complexity reduction methods for {ROM}s of dimensions $n=48$ and $n=66$, the resulting reduction errors are shown in Fig.~\ref{fig:gl38crer}.
The first observation to be made is that DEIM yields in several parameter settings unstable results. Particularly in the untrained case, all DEIM models with $n_c\leq 92$ have either a simulation breakdown or very poor fidelity. Our CROM method shows to be much more robust and no simulation breakdowns occur. Let us also emphasize that we have a strong a priori indicator for stability for the CROM, which is the compatibility condition $\sigma(\MassBoth_{c}) \subset [\tilde{C}^{-2},\tilde{C}^2]$ for a reasonably small constant $\tilde{C}$. As Table~\ref{table:gl38-cond} shows, the condition is fulfilled in our experiments.
Another difference of the complexity reduction methods we want to stress is that our quadrature-based complexity reduction in CROM depends stronger on the underlying ROM dimension $n$. Particularly, when perfectly trained and $n$ is chosen large, the error of DEIM becomes negligible for $n_c \approx 92$, whereas this holds only from $n_c \approx 120$ on for the CROM, see Fig.~\ref{fig:gl38crer}-\textit{left} with $n=66$. However, this does not mean that DEIM is more efficient in this or any other setting. As mentioned above, more nonlinear integrals for each $n_c$ have to be evaluated in DEIM than for CROM. This is directly reflected in the simulation times, which we can be found in Table~\ref{table:gl38-runtime} for the untrained case. The DEIM models have runtimes between $432$ to $482$ seconds, as compared to the CROMs with about $126$ to $152$ seconds. Thus, the CROM is more than three times faster than DEIM for models of similar fidelity. Compared to the $2253.6$ seconds runtime of FOM, CROM shows a speedup of about $16$.

To summarize, our quadrature-based reduction yields significantly more robust and efficient results as DEIM when combined with our model order reduction approach. Moreover, the compatibility condition of Assumption~\ref{assum:quadrat-ansatz} is conveniently promoted by our greedy training procedure in all our tests.

\begin{rmrk}\label{rem:useOfDeim}
Instead of applying  DEIM to $\tilde{\bv f}^\beta$ as we did, one could apply it directly to $\bv f^\beta$. This alternative showed severe stability issues in our numerical tests. We assume this is related to a more profound loss of structural properties. Specifically, the anti-symmetry (symplectic structure) revolving around the terms $-\bt J$ and $\bt{J}^T$ is destroyed by the more naive DEIM version, cf., \cite{art:symplecticMOR-Buchfink19,art:symplHamMor} for related discussions. 
\end{rmrk}

\section{Well-posedness of approximations} \label{sec:well-posed}

In this section, we derive a well-posedness result for our model order- and complexity-reduced approximation ({CROM}). The {CROM} is a realization of System~\ref{sys:abstract}, and its ansatz space $\funSpace{V}_r$ is a subspace of a finite element space $\funSpace{V}_f$ with partitioning $K_1,\ldots,K_J$ given as in \eqref{eq:FOM-FE}. Further, the complexity-reduced bilinear form and Hamiltonian take, according to \eqref{eq:quadrat-ansatz}, the form
\begin{align*}
\Lscal b, \bar{b}\Rscal_c = \sum_{i\in I} \wei_i \int_{K_i} b(x) \bar{b}(x) dx,  \quad \text{ with } \wei_i > 0 \text{ for }i \in I, \hspace{1cm}  \text{and } \hspace{1cm}   \HamPDE_c(\rho,m)= \Lscal h(\rho,m),1 \Rscal_c.
\end{align*}
Here, we assume $\wei_i>0$ instead of $\wei_i \geq 0$ for technical reasons, which can be done w.l.o.g.~by simply restricting $I$ accordingly. 
The well-posedness result is established under the situation described in System~\ref{sys:abstract} (i.e.,\,for boundary conditions as in \eqref{bls-eq:abstr-bc}) and the following additional assumptions on the model problem.
\begin{assumption} \label{assum:wellposed} \quad
\begin{enumerate} [{A}1)] 
\item The pressure potential $P: (0,M) \rightarrow \mathbb{R}$ for $M>0$ is two times continuously differentiable. It holds $P''(y) >0$, and $y \leq \max \{P(y),1 \}$ for $y\in (0,M)$, and $\lim_{y\rightarrow M} P(y) = \infty$.
\item The initial conditions $(\rho_0,m_0) \in \funSpace{V}_r$  are chosen such that $\rho_{0}(x) \in (0,M)$ for $x\in K_i$, $i\in I$.
\item  The friction model takes the form $r(\rho,m) = \rho^{-2}$ (laminar friction model for gas pipelines).
\end{enumerate}
\end{assumption}

Assumptions $(A1)$-$(A2)$ are used to show that $\HamPDE_c$ is well-defined and bounds the norm of the solution. Moreover, it should be mentioned that $P$ is also bounded from below by $(A1)$. To establish a global existence result, we derive uniform boundedness of the solution and strict positivity of the density. The positivity of the density is needed to avoid zeros in the denominators of System~\ref{sys:abstract}. Our proof of that relies on $(A2)$-$(A3)$. Further note that the restriction of the density to a single element $K_i$, $\rho(t)_{|K_i}$, is constant.
\begin{rmrk}
Similar results as the ones derived in this section can be established using other boundary conditions or friction/dissipation terms, but the derivations might get slightly more technical \cite{phd:liljegren}. We also refer to \cite{art:egger-mfem-compressEuler, art:fD-1D-Eulerisentr} for related results for other space discretizations of the Euler equations. However, the latter references do not treat model order- and complexity-reduction.
\end{rmrk}

\begin{lmm}\label{lem:iso1-comred-rhoInv}
Under the assumptions of this section, the solutions of System~\ref{sys:abstract} fulfill $\rho(t)_{|K_i} >0$ for $t>0$ and $i\in I$. Moreover, there exists a constant $C$, independent of the discretization parameters, such that
\begin{align*}
		\frac{1}{\rho(t)_{|K_i}} &\leq \exp{\left(C \frac{\sqrt{t}}{{\Delta_{x_i}} \sqrt{\wei_i} } \sqrt{R(t)}   \right) } \frac{1}{\rho(0)_{|K_i}} \\
		& \text{ for }R(t) = \HamPDE_c ( \rho(0), m(0) )+ \max \left \{0 \, , \, -\mathrm{inf}_{y\in(0,M)} P(y)  \right\}+ \int_{0}^t  \bv{u}(s)	\cdot \TraceOp m(s) ds , \hspace{0.8cm}  \text{for $i\in I$},
\end{align*}
where ${\Delta_{x_i}}$ is the grid size of the finite element $K_i$.
\end{lmm}
\begin{proof}
	At several instances, we employ that $\rho(t)_{|K_i}$ and $\partial_x {{m}}_{|K_i}(t)$ are constant in space. As long as $1/\rho(t)_{|K_i}$ is well-defined, it therefore holds
	\begin{align*}
		\frac{d}{dt} \left( \frac{1}{\rho(t)_{|K_i}} \right) &=  \left\Lscal \partial_t \frac{1}{\rho(t)},1\right\Rscal_{K_i} = -\left\Lscal \partial_t \rho(t),\frac{1}{\rho(t)^2}\right\Rscal_{K_i} \\
		&= \left\Lscal\partial_x {{m}}(t),\frac{1}{\rho(t)^2}\right\Rscal_{K_i} 
		 \leq  \left( \frac{1}{\rho(t)_{|K_i}} \left|\left|\partial_x {{m}}(t)\right|\right|_{K_i,\infty} \right) \frac{1}{\rho(t)_{|K_i}} ,
	\end{align*}
where the subscript $K_i$ indicates a restriction of the spatial domain to one finite element. By the inverse estimate, there exists a constant $C$ with $||\partial_x m(t)||_{K_i} \leq {C}/{{\Delta_{x_i}}} ||m(t)||_{K_i}$. Together with $1/\rho(t)$ and $\partial_x m(t)$ both being constant on $K_i$, this yields
\begin{align*}
\frac{1}{\rho(t)_{|K_i}} \left|\left|\partial_x {{m}}(t)\right|\right|_{K_i,\infty}  = \frac{1}{\rho(t)_{|K_i}} \left|\left|\partial_x {{m}}(t)\right|\right|_{K_i} 
\leq  c_{i}(t) , \hspace{1cm}
 \text{ with   } c_{i}(s) = \frac{C}{{\Delta_{x_i}}} \left|\left|\frac{ m(s)}{\rho(s)}\right|\right|_{K_i} .
\end{align*}
	
Setting together the two estimates, we thus get by the Gronwall lemma
\begin{align*}
	\frac{1}{\rho(t)_{|K_i}}  \leq  \exp \left( \int_{0}^t c_{i}(s) ds \right) \ \frac{1}{\rho(0)_{|K_i}} .
\end{align*}
It remains to bound $\int_{0}^t c_{i}(s) ds $ for $i \in I$. As a preparation, we introduce the auxiliary function $\hat{\HamPDE}( \rho, m) = \HamPDE_c( \rho,m) + \max \left \{0 \, , \, -\mathrm{inf}_{y\in(0,M)} P(y)  \right\}$, which is a constant shift of the Hamiltonian $\HamPDE_c$. By construction, it yields only non-negative values for all considered states and fulfills the same energy-dissipation equality as $\HamPDE_c$, cf., Theorem~\ref{thrm:energy-diss-cr}. Assuming $\rho(t)_{|K_i} > 0$ and $i \in I$, we can follow
\begin{align*}
	\wei_i \int_{0}^T \left|\left|\frac{ m(t)}{\rho(t)}\right|\right|_{K_i}^2  dt &\leq \sum_{i \in I} \wei_i \int_{0}^T \left|\left|\frac{ m(t)}{\rho(t)}\right|\right|_{K_i}^2  dt = 
	\int_{0}^T	\left \Lscal \left(\frac{m(t)}{\rho(t)}\right)^2 , 1 \right \Rscal_c  dt \\
	&\leq \int_{0}^T	\left \Lscal \left(\frac{m(t)}{\rho(t)}\right)^2 , 1 \right \Rscal_c  dt + \hat{\HamPDE} (\rho(t),m(t)) = R(t)
\end{align*}
wit $R$ defined as in the lemma. The latter equality follows from integrating the energy dissipation equality for $\hat{\HamPDE}$ (Theorem~\ref{thrm:energy-diss-cr}) in time. Moreover, $R(t)\geq 0$ holds for $t\geq 0$, as the last equality also implies that $R$ can be bounded from below by $ \hat{\HamPDE}$, which itself is non-negative. Now first applying the Jensen-inequality \cite{book:RockWets98} and then inserting the former estimate yields
\begin{align*}
	\int_{0}^t c_i(s) ds &\leq \left( t \int_{0}^t c_i(s)^2 ds \right)^{1/2}  = C \frac{\sqrt{t}}{{\Delta_{x_i}}}  \left( \int_{0}^t \left|\left|\frac{ m(t)}{\rho(t)}\right|\right|_{K_i}^2  ds \right)^{1/2} \leq  C \frac{\sqrt{t}}{{\Delta_{x_i}} \sqrt{\wei_i}} \sqrt{R(t)}.
\end{align*}
Inserting the latter bound on $\int_{0}^t c_i(s) ds$ into the estimate obtained by the Gronwall lemma finishes the proof.
	\end{proof}
	
\begin{rmrk}
In the special case of $\Lscal \cdot, \cdot \Rscal_c$ chosen as the $\funSpace{L}^2$-scalar product, we recover the case of pure Galerkin approximation, without complexity reduction. Lemma~\ref{lem:iso1-comred-rhoInv} then shows strict positivity of $\rho(t)$ on all of the spatial domain. Such a result has been derived in \cite{art:egger-mfem-compressEuler} in a similar setting.
\end{rmrk}
Note that positivity of $\rho(t)_{|K_i}$ can only be guaranteed for $i\in I$ by Lemma~\ref{lem:iso1-comred-rhoInv}. This turns out to be sufficient, as the $\rho(t)$-terms in the denominator are only evaluated for $i\in I$. Next, we derive a boundedness result for the solution.


\begin{thrm}\label{theor:iso1-comred-bound}
	Under the assumptions of this section, there exist constants $C_1,C_2$, independent of the discretization parameters, such that for $(\rho,m) \in \funSpace{V}_r$ with $\rho_{|K_i}>0$ for $i \in I$ it holds
	\begin{align*}
		||\rho|| + ||m|| &\leq \left(\max_{i\in I}{\wei_i}^{-\frac{1}{2}}\right) \left[ C_1 \HamPDE_c(\rho,m) + C_2 \right].
	\end{align*}
\end{thrm}
\begin{proof}
Let $k \in \text{argmax}_{i\in I} \rho_{|K_i}$, and let $\bar{I} = \{i\in I: \rho_{|K_i} > 1]\}$ w.l.o.g.~be non-empty. Otherwise we trivially can bound the terms $||\rho||_c, \rho_{|K_k}$ by a constant. From Assumption~\ref{assum:quadrat-ansatz} we get
	\begin{align*}
		\frac{1}{\tilde{C}^2} ||\rho||^2 &\leq ||\rho||_c^2 = \Lscal\rho^2,1\Rscal_c =  \sum_{i\in I} \wei_i \rho_{|K_i}^2 \leq \rho_{|K_k}  \sum_{i\in I} \wei_i \rho_{|K_i} \\
		&= \rho_{|K_k}  \left(\sum_{i\in \bar{I}} \wei_i \rho_{|K_i} + \sum_{j\in I \text{\textbackslash} \bar{I}} \wei_j \rho_{|K_j}\right)
		\leq  P(\rho_{|K_k} ) \left( \sum_{i\in \bar{I}} \wei_i P(\rho_{|K_i}) + \sum_{j\in I \text{\textbackslash} \bar{I}} \wei_j \right).
	\end{align*}
	By Assumption~\ref{assum:wellposed}-$(A1)$ it follows that the function $P$ can be bounded from below (not necessarily by zero), which implies that there exists a constant $\hat{c}$ such that
	\begin{align*}
		\wei_k P(\rho_{|K_k} ) \leq \sum_{i\in \bar{I}} \wei_i P(\rho_{|K_i}) \leq \Lscal P(\rho),1 \Rscal_c + \hat{c}.
	\end{align*}
	Together, this shows that for some $\hat{c}_1,\hat{c}_2>0$ it holds
	\begin{align*}
		||\rho|| \leq \frac{1}{\sqrt{\wei_k}} \left(\hat{c}_1  \Lscal P(\rho),1 \Rscal_c +\hat{c}_2 \right) .
	\end{align*}
	Similarly, it follows
	\begin{align*}
		\frac{1}{\tilde{C}^2} ||m||^2 &\leq ||m||_c^2 = \left\Lscal \frac{m^2}{\rho} \rho,1 \right\Rscal_c \leq \rho_{|K_k} \left\Lscal \frac{m^2}{\rho} ,1 \right\Rscal_c \leq P(\rho_{|K_k}) \left\Lscal \frac{m^2}{\rho} ,1 \right\Rscal_c \\
		&\leq  \frac{1}{\wei_k}  ( \Lscal P(\rho),1\Rscal_c + \hat{c}) \left\Lscal \frac{m^2}{\rho} ,1 \right\Rscal_c,
\end{align*}
where in the last step the same estimate on $P(\rho_{|K_k} )$ as before has been used. With the help of Young's inequality, it follows
	\begin{align*}
		||m|| \leq   \frac{\tilde{C}}{\sqrt{\wei_k}} \sqrt{(\Lscal P(\rho),1\Rscal_c + \hat{c})}  \, \sqrt{\left\Lscal \frac{m^2}{\rho} ,1 \right\Rscal_c} \leq \frac{\tilde{C}}{2\sqrt{\wei_k}} \left( \Lscal P(\rho),1 \Rscal_c  + \left\Lscal \frac{m^2}{\rho} ,1 \right\Rscal_c  + \hat{c}\right) .
	\end{align*}
	Setting together the estimates for $||\rho||$ and $||m||$ shows the assertion.
\end{proof}
Notably, the bound on the $\funSpace{L}^2$-norms in Theorem~\ref{theor:iso1-comred-bound} is almost independent of the discretization parameters. Only the quadrature weights of the complexity reduction step enter. To make the bound uniform, one has to require the quadrature weights to be bounded from below by a positive constant. The concluding result now reads as follows.

\begin{thrm} \label{theor:cr-iso1-wellposed}
Under the assumptions of this section, System~\ref{sys:abstract} has a unique solution $(\rho,m) \in \funSpace{C}^1([0,T);\Vsa\times \Vsb)$ for any $T >0$.
\end{thrm}
\begin{proof}
Given the solution $(\rho,m)$ exists up $t$, $t>0$, we define $\mathcal{F}(t):= \left(\max_{i\in I}{\wei_i}^{-\frac{1}{2}}\right) \left[ C_1 \HamPDE_c(\rho(t),m(t)) + C_2 \right]$ with constants $C_1, C_2$ as in Theorem~\ref{theor:iso1-comred-bound}. From the energy dissipation (Theorem~\ref{thrm:energy-diss-cr}) and the boundedness of the trace operator $\TraceOp$, it follows 
\begin{align*}
 \frac{d}{dt} \mathcal{F}(t) =	C_3 \frac{d}{dt} \HamPDE_c(\ubv{a}(t)) \leq {C}_3 \bv u(t) \cdot \TraceOp^{ }m(t) \leq  C_4 ||\bv u(t)|| \, || m(t)||
\end{align*}
for constants $C_3,C_4>0$. Further, $||\rho(t)||+|| m(t)|| \leq \mathcal{F}(t)$ holds by Theorem~\ref{theor:iso1-comred-bound}. Thus, it holds by the former estimate that ${d}/{dt} \mathcal{F}(t) \leq C_4 ||\bv u(t)|| \mathcal{F}(t)$. We follow by the Gronwall Lemma that $\mathcal{F}$ grows at most exponentially in time, and accordingly the same holds for the norm of the solution $(\rho,m)$. By that existence of a continuously differentiable solution for any $t>0$ can be deduced by the Peano existence theorem and the extension theorem for differential equations by standard arguments for ordinary differential equations, cf. \cite{book:hartman-ode}.
\end{proof}

\section{Conclusion}
We proposed a snapshot-based model reduction approach for a nonlinear flow problem on networks that is governed by the barotropic Euler equations with friction. The approach consists of a structure-preserving modification of the proper orthogonal decomposition combined with a quadrature-type complexity reduction. It yields online-efficient reduced models with remarkable structure-preserving properties, i.e., local mass conservation, an energy bound and a port-Hamiltonian structure. These properties are established using a few compatibility conditions, which we assure by appropriate adaptions in the training phase. The involved implications and appropriate efficient algorithmic implementations were a main focus of this paper. Notably, our model order reduction spaces fulfill an optimality condition under the compatibility conditions. As we demonstrated on a realistic gas transportation network benchmark, our reduced models show superior stability and overall performance compared to more generic non-structure preserving model reduction approaches. Moreover, a well-posedness result for them was established under a few additional assumptions.

\bibliographystyle{abbrv}
\bibliography{NLFlowMOR}

%

\end{document}